\pdfoutput=1

\documentclass[leqno]{prims}
\usepackage{hyperref}
\usepackage{amssymb,amsmath}
\usepackage{enumerate}

\newtheorem{thm}{Theorem}[section]
\newtheorem{cor}[thm]{Corollary}
\newtheorem{lem}[thm]{Lemma}

\newtheorem{mainthm}[thm]{Main Theorem}

\theoremstyle{definition}

\numberwithin{equation}{section}

\usepackage{mathrsfs}
\usepackage[all, pdf]{xy}

\newtheorem{dfn}[thm]{Definition}
\newtheorem{prop}[thm]{Proposition}
\newtheorem{example}[thm]{Example}

\newcommand{\C}{{\mathbb C}}

\newcommand{\Q}{{\mathbb Q}}
\newcommand{\Z}{{\mathbb Z}}

\newcommand{\normord}[1]{{:}\!\mathrel{#1}\!{:}}
\newcommand{\relmiddle}[1]{\mathrel{}\middle#1\mathrel{}}
\def\maprestrict#1#2{\left. #1 \right|_{#2}}

\VolumeNo{4x}
\YearNo{201x}
\communication{T. Arakawa. Received September 10, 2021.}

\begin{document}



\TitleHead{Quantum dilogarithm identities}
\title{Quantum dilogarithm identities arising from the product formula for universal R-matrix of quantum affine algebras}


\AuthorHead{M. Sugawara}
\author{Masaru \textsc{Sugawara}\footnote{M. Sugawara: Mathematical Institute, Tohoku University, Sendai 980-8578, Japan;
\email{masaru.sugawara.s7@dc.tohoku.ac.jp}}}


\maketitle

\begin{abstract}
In ~\cite{DGS}, four quantum dilogarithm identities containing infinitely many factors are proposed as wall-crossing formula for refined BPS invariant. We give algebraic proof of these identities using the formula for universal R-matrix of quantum affine algebra developed by K. Ito~\cite{Ito2}, which yields various product presentation of universal R-matrix by choosing various convex orders on affine root system. By the uniqueness of universal R-matrix and appropriate degeneration, we can construct various quantum dilogarithm identities including the ones proposed in \cite{DGS}, which turn out to correspond to convex orders of multiple row type.
\end{abstract}

\section{Introduction}

Dimofte, Gukov, and Soibelman proposed four remarkable identities with respect to quantum dilogarithm functions as the wall-crossing formulas for the refined BPS invariants, which they proposed in the study of type II string theory~\cite{DGS2}. In ~\cite{DGS2}, it is observed that the refined BPS invariants has very similar wall-crossing behavior to that of motivic Donaldson-Thomas invariants introduced by Kontsevich and Soibelman ~\cite{KS}, and it is conjectured that the two invariants coincide under appropriate identification of variables.

Let $q, x_1, x_2$ be indeterminate satisfying the relations $q x_1 = x_1 q$, $q x_2 = x_2 q$, $x_1 x_2 = q^2 x_2 x_1$, and let
\begin{equation}
\mathbb{E} (x) := \prod_{k = 0}^{\infty} \frac{1}{1 + q^{2k + 1}x}, \quad {\bf U}_{m, n} := \mathbb{E}(q^{-mn} x_1^m x_2^n).
\end{equation}
Then the identities they found are written down as following~\cite{DGS}. Note that the parameter $q$ in this paper corresponds to $q^{1/2}$ in ~\cite{DGS}.
\begin{align}
{\bf U}_{2, -1} {\bf U}_{0, 1} &= ({\bf U}_{0, 1} {\bf U}_{2, 1} {\bf U}_{4, 1} \dots) \label{eq:A1} \\
&\quad \times \mathbb{E} (-q x_1^2)^{-1} \mathbb{E} (-q^{-1} x_1^2)^{-1} \notag \\
&\quad \quad \times (\dots {\bf U}_{6, -1} {\bf U}_{4, -1} {\bf U}_{2, -1}), \notag \\
{\bf U}_{1, -1} {\bf U}_{1, 0} {\bf U}_{0, 1} &= ({\bf U}_{0, 1} {\bf U}_{1, 1} {\bf U}_{2, 1} {\bf U}_{3, 1} \dots) \label{eq:A2} \\
&\qquad \times {\bf U}_{1, 0}^2 \mathbb{E} (-q x_1^2)^{-1} \mathbb{E} (-q^{-1} x_1^2)^{-1} \notag \\
&\qquad \qquad \times (\dots {\bf U}_{3, -1} {\bf U}_{2, -1} {\bf U}_{1, -1}), \notag \\
{\bf U}_{1, -1}^2 {\bf U}_{0, 1}^2 &= ({\bf U}_{0, 1}^2 {\bf U}_{1, 1}^2 {\bf U}_{2, 1}^2 {\bf U}_{3, 1}^2 \dots) \label{eq:A3} \\
&\quad \times {\bf U}_{1, 0}^4 \mathbb{E} (-q x_1^2)^{-1} \mathbb{E} (-q^{-1} x_1^2)^{-1} \notag \\
&\quad \quad \times (\dots {\bf U}_{3, -1}^2 {\bf U}_{2, -1}^2 {\bf U}_{1, -1}^2), \notag \\
{\bf U}_{1, -2} {\bf U}_{0, 1}^4 &= ({\bf U}_{0, 1}^4 {\bf U}_{1, 2} {\bf U}_{1, 1}^4 {\bf U}_{3, 2} {\bf U}_{2, 1}^4 \dots) \label{eq:D4} \\
&\quad \times {\bf U}_{1, 0}^6 \mathbb{E} (-q x_1^2)^{-1} \mathbb{E} (-q^{-1} x_1^2)^{-1} \notag \\
&\quad \quad \times (\dots {\bf U}_{2, -1}^4 {\bf U}_{3, -2} {\bf U}_{1, -1}^4 {\bf U}_{1, -2}). \notag
\end{align}
The function $\mathbb{E} (x)$ is called quantum dilogarithm since
\begin{equation}
\mathbb{E} (x) = \exp \left(\mathrm{Li}_{2, q^2} \left(-qx\right)\right), \quad \mathrm{Li}_{2, q} (x) := \sum_{n = 1}^{\infty} \frac{x^n}{n(1 - q^{n})}
\end{equation}
and $(1 - q) \mathrm{Li}_{2, q} (x)$ degenerates to classical dilogarithm by $q \to 1$.

These identities, however, are derived by physical insight, and mathematically rigorous proofs for them have not been given. In this paper, we develop algebraic construction of these identities, which eventually yields mathematical proof of them as equalities of skew formal power series.

In ~\cite{KN}, Kashaev and Nakanishi established systematic construction of quantum dilogarithm identities from periods of quantum cluster algebras. Their identities, however, involve only finitely many factors, while the four identities (\ref{eq:A1}), (\ref{eq:A2}), (\ref{eq:A3}), (\ref{eq:D4}) contain infinite product. Thus, these identities belong to essentially new class of quantum dilogarithm identities.

On the other hand, K. Ito constructed the product formulas for the (quasi-) universal R-matrix of quantum affine algebra, which correspond to convex orders on affine root system~\cite{Ito2}. In the formulas, the factors corresponding to real roots are $q$-exponential function, which is in fact written as $\exp_q (x) = \mathbb{E} ((q - q^{-1}) x)$. The resemblance between the wall-crossing formulas and product formulas for the universal R-matrix implies existence of connection between wall-crossing formulas and quantum groups.

By this observation, we develop systematic construction of quantum dilogarithm identities containing infinite product, using the product formula for the universal R-matrix. As a result, we show that all four identities Dimofte et al.\ found can be derived algebraically by our method. In section \ref{sec:2}, we review general construction of convex orders on affine root systems, concrete construction of PBW type bases for the positive part $U_q^+$ of quantum affine algebra $U_q(\mathfrak{g})$ using convex order, and the explicit product formula for the quasi-universal R-matrix of $U_q(\mathfrak{g})$.

In section \ref{sec:4}, we show how to construct quantum dilogarithm identities using quasi-universal R-matrix $\Theta$ of $U_q(\mathfrak{g})$. By virtue of the uniqueness of $\Theta$, we can equate all the product presentations of $\Theta$ associated with convex orders. Thus we have infinite product identities whose parts corresponding to real roots are $q$-exponential function of root vectors. Next, we construct continuous projection of the completed quantum double algebra $U_q^+ \widehat{\otimes} U_q^-$, which contains $\Theta$, onto skew formal power series algebra $\mathcal{D}_Q$ associated with affine Dynkin quiver $Q$. By this projection, some root vectors vanish and thus their $q$-exponential become $1$ in the image. If one choose appropriate convex order and Dynkin quiver $Q$, one can make infinitely many root vectors not to vanish for the convex order, while only finitely many root vectors alive in the image for reversed convex order. Eventually one can obtain various quantum dilogarithm identities of the form ``finite product = infinite product''.

To obtain concrete identities, we have to compute the root vectors explicitly to determine whether they vanish by the projection. In section \ref{sec:2.5}, we show that every root vector can be written as ``$q$-commutator monomial'', which is finite application of $q$-bracket on the Chevalley generators. We also developed combinatorial algorithm for the computation of root vectors, which enable us to obtain concrete presentations of root vectors as $q$-commutator monomials. The computation is done as manipulations of binary trees.

We concretely found appropriate convex orders and Dynkin quivers which produce identical identities with (\ref{eq:A1}) to  (\ref{eq:D4}), which will be explicitly presented in section \ref{sec:5}. It is remarkable that the factor of $\Theta$ correspond to imaginary roots becomes $q$-exponential function by the projection, in spite of the factor itself is not $q$-exponential function. We also note that the convex orders correspond to (\ref{eq:A2}), (\ref{eq:A3}), and (\ref{eq:D4}) are of multiple row type, which was newly found by Ito~\cite{Ito1}.

\section{Product formula for the universal R-matrix of quantum affine algebras} \label{sec:2}

First we summarize Ito's works~\cite{Ito1}~\cite{Ito2} which provide explicit product presentations of the (quasi-) universal R-matrix of quantum affine algebras. 

\subsection{Quantum algebra $U_q(\mathfrak{g})$}

To begin with, we recall quantum enveloping algebra $U_q(\mathfrak{g})$ corresponding to symmetrizable Kac-Moody algebra $\mathfrak{g}$ of rank $\ell + 1$, where $q$ is an indeterminate (thus we work on generic case). We use following notations as in ~\cite{Kac}.
\begin{align*}
e_i, f_i \in \mathfrak{g}&: \text{Chevalley generators} , \\
\mathfrak{h} \subset \mathfrak{g}&: \text{Cartan subalgebra}, \\
\Check{\alpha}_i \in \mathfrak{h}&: \text{simple coroots}, \\
\alpha_i \in \mathfrak{h}^*&: \text{simple roots}, \\
s_i \in \mathrm{End}(\mathfrak{h}^*)&: \text{simple reflections} \quad (i = 0, 1, \dots, \ell), \\
\Delta \subset \mathfrak{h}^*&: \text{set of all roots},  \\
W := \langle s_0, s_1, \dots, s_{\ell} \rangle&: \text{Weyl group}, \\
\Delta_+ \subset \Delta&: \text{set of all positive roots}, \\
\Delta_- \subset \Delta&: \text{set of all negative roots}, \\
\Delta^{\mathrm{re}} \subset \Delta&: \text{set of all real roots}, \\
\Delta^{\mathrm{im}} := \Delta \setminus \Delta^{\mathrm{re}}&:  \text{set of all imaginary roots}.
\end{align*}
We also use the symbol $R_{\pm} := R \cap \Delta_{\pm}$ for every $R \subset \Delta$.

\begin{dfn}
The quantum enveloping algebra $U_q(\mathfrak{g})$ is the associative $\Q(q)$-algebra defined by following generators and relations:
\begin{align}
\text{Generators} &: E_i, F_i, K_{\lambda} \quad (i = 0, 1, \dots, \ell; \, \lambda \in P). \\
\text{Relations} &: K_{\lambda} K_{\mu} = K_{\lambda + \mu}, \quad K_0 = 1, \\
& K_{\lambda} E_i K_{\lambda}^{-1} = q^{(\lambda, \alpha_i)}E_i, \quad K_{\lambda} F_i K_{\lambda}^{-1} = q^{-(\lambda, \alpha_i)}F_i, \\
& [E_i, F_j] = \delta_{ij} \frac{K_i - K_i^{-1}}{q_i - q_i^{-1}} \quad (i = 0, 1, \dots \ell; \, \lambda, \mu \in P); \label{eq:3rd_relation_of_Uq} \\
& \sum_{k = 0}^{1 - a_{ij}} (-1)^{k} {\genfrac{[}{]}{0pt}{}{1 - a_{ij}}{k}}_{q_i} E_i^{1 - a_{ij} - k} E_j E_i^{k} = 0, \label{eq:q_Serre_rel} \\
& \sum_{k = 0}^{1 - a_{ij}} (-1)^{k} {\genfrac{[}{]}{0pt}{}{1 - a_{ij}}{k}}_{q_i} F_i^{1 - a_{ij} - k} F_j F_i^{k} = 0 \quad (i \neq j),
\end{align}
where $(\cdot, \cdot)$ is the invariant bilinear form on $\mathfrak{h}^*$,
\begin{equation*}
P := \left\{\,\lambda \in \mathfrak{h}^* \relmiddle| \langle\Check{\alpha_i}, \lambda\rangle \in \Z \, (\forall i = 0, 1, \dots, \ell)\,\right\}
\end{equation*}
is the weight lattice, and let $a_{ij} := \frac{2(\alpha_i, \alpha_j)}{(\alpha_i, \alpha_i)} \in \Z$, $q_i := q^{\frac{1}{2}(\alpha_i, \alpha_i)}$, $[n]_q := \frac{q^n - q^{-n}}{q - q^{-1}}$, $[n]_q! := [1]_q [2]_q \dots [n_q]$, ${\genfrac{[}{]}{0pt}{}{n}{k}}_q := \frac{[n]_q!}{[k]_q! [n - k]_q!} \in \Z[q, q^{-1}]$.
\end{dfn}

It is well-known that $U_q(\mathfrak{g})$ becomes a Hopf algebra with following coalgebra structure $(U_q(\mathfrak{g}), \Delta, \varepsilon)$ and antipode $S$.
\begin{align*}
\Delta(E_i) &:= E_i \otimes 1 + K_i \otimes E_i, \quad \Delta(F_i) := F_i \otimes K_i^{-1} + 1 \otimes F_i, \\
\Delta(K_{\lambda}) &:= K_{\lambda} \otimes K_{\lambda}, \quad \varepsilon(E_i) := 0, \quad \varepsilon(F_i) := 0, \quad \varepsilon(K_{\lambda}) := 1, \\
S(E_i) &:= -K_i^{-1} E_i, \quad S(F_i) := -F_i K_i, \quad S(K_{\lambda}) := K_{\lambda}^{-1}.
\end{align*}
$\Delta: U_q(\mathfrak{g}) \rightarrow U_q(\mathfrak{g}) \otimes U_q(\mathfrak{g})$ and $\varepsilon: U_q(\mathfrak{g}) \rightarrow \Q(q)$ are uniquely extended as algebra homomorphisms, and $S: U_q(\mathfrak{g}) \rightarrow U_q(\mathfrak{g})$ is also extended as an anti-automorphism.

Several subalgebras of $U_q(\mathfrak{g})$ generated by standard generators are defined as usual:
\begin{equation*}
U_q^+ := \left\langle\, E_0, E_1, \dots, E_{\ell} \,\right\rangle, \quad U_q^0 := \left\langle\, K_{\lambda} \relmiddle| \lambda \in P \,\right\rangle, \quad U_q^- := \left\langle\, F_0, F_1, \dots, F_{\ell} \,\right\rangle.
\end{equation*}
Then we have triangular decomposition of $U_q(\mathfrak{g})$~\cite[3.2.5]{Lus}.
\begin{equation} \label{eq:tri_decomp}
U_q^- \otimes U_q^0 \otimes U_q^+ \cong U_q(\mathfrak{g}), \quad x \otimes y \otimes z \mapsto xyz.
\end{equation}

Let $U_{\mu} := \left\{\, x \in U_q(\mathfrak{g}) \relmiddle| K_{\lambda} x K_{\lambda}^{-1} = q^{(\lambda, \mu)}x \, (\forall \lambda \in P) \,\right\}$ be the weight space of weight $\mu \in P$. For convenience, let $V_{\mu} := V \cap U_{\mu}$ for every subspace $V \subset U_q(\mathfrak{g})$. Then we also have the weight space decomposition
\begin{equation}
U_q(\mathfrak{g}) = \bigoplus_{\mu \in Q} U_{\mu} \quad (Q := \bigoplus_{i = 0}^{\ell} \Z \alpha_i \subset P: \text{root lattice}),
\end{equation}
and $U_q(\mathfrak{g})$ becomes a $Q$-graded algebra. Using this gradation, we introduce \textbf{$q$-bracket} $[\cdot, \cdot]_q$, which is defined on each weight spaces as follows:
\begin{equation}
[x, y]_q := xy - q^{(\mu, \nu)} yx \quad (\mu, \nu \in Q; \, x \in U_{\mu}, y \in U_{\nu}).
\end{equation}

\subsection{Convex orders on affine root system}

Next, we introduce the definition and classification of convex orders on the set of positive roots~\cite{Ito1}. We also prepare a number of notations on affine root systems.

\begin{dfn}~\cite[Definition 3.3]{Ito2} A total order $\leq$ on a set of positive roots $B \subset \Delta_+$ is called \textbf{convex} if it satisfies the following two conditions.
\begin{enumerate}
\item For any pair of positive real roots $\beta, \gamma \in B \cap \Delta_+^{\mathrm{re}}$ satisfying $\beta < \gamma$ and $\beta + \gamma \in B$, the order relation $\beta < \beta + \gamma < \gamma$ holds.

\item If $\beta \in B, \gamma \in \Delta_+ \setminus B$ and $\beta + \gamma \in B$, then $\beta < \beta + \gamma$.
\end{enumerate}
\end{dfn}

\begin{example} Set $\mathfrak{g} = \widehat{\mathfrak{sl}_2}$. Then the following order on $\Delta_+$ is convex.
\begin{equation}\label{eq:order_sl2}
\begin{aligned}
\delta - \alpha_1 &< 2\delta - \alpha_1 < 3 \delta - \alpha_1 < \dots \\
&\quad < \delta < 2 \delta < 3 \delta < \dots \\
&\quad \quad \dots < 2\delta + \alpha_1 < \delta + \alpha_1 < \alpha_1.
\end{aligned}
\end{equation}
Here $\delta := \alpha_0 + \alpha_1$ is null root. 
\end{example}

When $\mathfrak{g}$ is of untwisted affine type, convex orders on $\Delta_+$ are already classified by Ito~\cite{Ito1}. To describe convex orders in general, we have to introduce a numerous amount of symbols on affine root systems. In the rest of this section, we restrict $\mathfrak{g}$ to be an untwisted affine Lie algebra of type $X_{\ell}^{(1)}$, where $X$ is one of $A, B, C, D, E, F, G$ and $\ell$ is a positive integer. We assign indices $0, 1, \dots, \ell $ for each vertex of the Dynkin diagram corresponds to $\mathfrak{g}$ as in ~\cite{IM} so that the full subdiagram without the vertex 0 is of finite $X_{\ell}$ type.

First, let $\mathring{I} := \{1, 2, \dots, \ell\}$ be the set of indices other than 0, and $\mathring{\mathfrak{g}} \subset \mathfrak{g}$ be the Lie subalgebra generated by $\left\{\,e_i, f_i, \Check{\alpha}_i \relmiddle| i \in \mathring{I}\,\right\}$. Then $\mathring{\mathfrak{g}}$ is isomorphic to the simple Lie algebra of type $X_{\ell}$ due to our assignment of indices, and $\mathring{\mathfrak{h}} := \bigoplus_{i \in \mathring{I}} \C \Check{\alpha}_i \subset \mathfrak{h}$ is a Cartan subalgebra of $\mathring{\mathfrak{g}}$. Let $\mathring{\Delta} \subset \mathring{\mathfrak{h}}^*$ be the set of all roots of $\mathring{\mathfrak{g}}$, and $\mathring{W} = \left\langle s_i \relmiddle| i \in \mathring{I} \right\rangle \subset W$ be the finite Weyl group.

Associated to each $J \subset \mathring{I}$, we introduce several symbols below~\cite{Ito1}.
\begin{align*}
\mathring{\Pi}_{J} &:= \left\{\, \alpha_j \relmiddle| j \in J\,\right\} \subset \mathfrak{h}^*, \quad \mathring{W}_{J} := \left\langle s_j \relmiddle| j \in J \right\rangle \subset \mathring{W}, \\
\mathring{W}^{J} &:= \left\{\, w \in \mathring{W} \relmiddle| w(\alpha_j) \in \mathring{\Delta}_+ \, (\forall j \in J) \,\right\}, \\
\mathring{\Delta}_{J} &:= \mathring{W}_{J} (\mathring{\Pi}_{J}), \quad \mathring{\Delta}^{J} := \mathring{\Delta} \setminus \mathring{\Delta}_{J}, \quad \mathring{\Delta}^{J}_{\pm} := \mathring{\Delta}^{J} \cap \mathring{\Delta}_{\pm}, \\
\Delta^{J}(w, \pm) &:= \left\{\, m\delta + \varepsilon \relmiddle| m \in \Z_{\geq 0}, \, \varepsilon \in w \mathring{\Delta}^{J}_{\pm} \,\right\} \cap \Delta_+ \quad (w \in \mathring{W}), \\
\Delta_{J}(w, \pm) &:= \left\{\, m\delta + \varepsilon \relmiddle| m \in \Z_{\geq 0}, \, \varepsilon \in w \mathring{\Delta}_{J \pm} \,\right\} \cap \Delta_+ \quad (w \in \mathring{W}),
\end{align*}
where $\delta \in \Delta^{\mathrm{im}}_{+}$ is the null root. For every symbol $X_J$ (resp. $X^J$), we omit the subscript (resp. superscript) $J$ and write $X := X_J$ (resp. $X := X^J$) when $J = \mathring{I}$ (resp. $J = \emptyset$). Notice that $X^{\emptyset} = X_{\mathring{I}} = X$ for every symbols introduced above.

Since every proper full subdiagram of affine Dynkin diagram is finite direct sum of diagrams of finite type~\cite{Kac}, so is root subsystem $\mathring{\Delta}_{J}$. Thus we have the partition $J = \coprod_{c \in C} J_c$, where $C$ is the set of connected components of Dynkin diagram of $\mathring{\Delta}_{J}$ and $J_c \subset J$ is the set of vertices belong to connected component $c \in C$. Each component $\mathring{\Delta}_{J_c}$ is irreducible root system of finite type, so that there exists unique highest root $\theta_{J_c} \in \mathring{\Delta}_{J_c +}$.

Moreover, several symbols are defined for each connected component $J_c$ as follows.
\begin{align*}
\Pi_{J_c} &:=  \mathring{\Pi}_{J_c} \amalg \{\delta - \theta_{J_c}\}, \quad \Pi_{J} := \coprod_{c \in C} \Pi_{J_c}, \\
S_{J} &:= \left\{\, s_{\alpha} \relmiddle| \alpha \in \Pi_{J} \,\right\}, \quad W_{J} := \langle S_{J} \rangle \subset W, \\
\Delta_{J}^{\mathrm{re}} &:= W_{J}(\Pi_{J}), \quad \Delta_J := \Delta_{J}^{\mathrm{re}} \amalg \Delta^{\mathrm{im}}.
\end{align*}

Furthermore, we associate each $y \in W_{J}$ with a set of positive roots $\Phi_{J} (y) := y\Delta_{J-} \cap \Delta_{J+}$. Also we set
\begin{equation*}
\nabla(J, u, y) := \Delta^{J}(u, -) \amalg u \Phi_{J}(y)
\end{equation*}
for each $u \in \mathring{W}^{J}$. These infinite sets of positive real roots $\nabla(J, u, y)$ have biconvex property and play crucial role in the classification of convex orders~\cite{Ito1}.

We also need to introduce a decomposition of elements of Weyl group defined by the next lemma.
\begin{lem}
\label{lem:w_decomposition}
For every $w \in W$, there exists unique decomposition $w = w^{J} w_{J}$, where $w^{J} \in W^{J}$, $w_{J} \in \mathring{W}_J$.
\end{lem}

By definition of $W_J$, each $w \in W_J$ can be written as a finite product of elements in $S_J$. An expression $w = t_1 t_2 \dots t_m$ $(t_i \in S_J)$ is called \textbf{reduced} if the number $m$ is smallest among all the expressions of $w$ as finite product of elements in $S_J$, and the smallest number $m$ is called \textbf{length} of $w$. Let $\ell_J(w)$ denote the length of $w$. An infinite sequence of elements $u_1, u_2, \dots $in $S_J$ is called \textbf{infinite reduced word} when $\ell_J(u_1 u_2 \dots u_m) = m$ for all positive integer $m$. The set of all infinite reduced words is denoted by $\mathscr{W}_J^{\infty}$, and the $k$-th factor of ${\bf s} \in \mathscr{W}_J^{\infty}$ is denoted by ${\bf s}(k) \in S_J$. We also use a function on positive integers $\phi_{\bf s}: \Z_{\geq 1} \rightarrow \Delta_J$ defined by $\phi_{\bf s} (k) := {\bf s}(1) {\bf s}(2) \dots {\bf s}(k - 1) (\beta_k)$, where $\beta_k \in \Pi_J$ is the positive root correspond to ${\bf s}(k) = s_{\beta_k} \in S_J$. Note that 
\begin{equation*}
\Phi_J\left({\bf s}\left(1\right) {\bf s}\left(2\right) \dots {\bf s}\left(k\right)\right) = \left\{\phi_{\bf s} \left(1\right), \phi_{\bf s} \left(2\right), \dots, \phi_{\bf s} \left(k\right)\right\} \quad (k \in \Z_{\geq 1}).
\end{equation*}
We associate each ${\bf s} \in \mathscr{W}_J^{\infty}$ with a infinite set of positive roots 
\begin{equation*}
\Phi_{J} ({\bf s}) := \left\{\, \phi_{\bf s} \left(k\right) \relmiddle| k \in \Z_{\geq 1} \,\right\}.
\end{equation*}

Now, we can state the general description of convex orders. To begin with, we pick an element $w \in \mathring{W}$. Then we have the decomposition
\begin{equation}
\label{eq:positive_root_partition}
\Delta_+ = \Delta(w, -) \amalg \Delta_+^{\mathrm{im}} \amalg \Delta(w, +).
\end{equation}
Note that $\Delta(w, +) = \Delta(ww_{\circ}, -)$ with the longest element $w_{\circ} \in \mathring{W}$, since $w_{\circ}$ reverses the sign of every root in $\mathring{\Delta}$. Thus the set of positive real roots consists of two sets of the form $\Delta(w, -)$. We will construct convex orders on $\Delta(w, -)$ and connect them to construct whole order.

Convex orders on $\Delta(w, -)$ are constructed by the following procedure.
\begin{enumerate}
\item Select a positive integer $n$ and a filtration of indices

$\mathring{I} = J_0 \supsetneqq J_1 \supsetneqq J_2 \supsetneqq \dots \supsetneqq J_n = \emptyset$.
\item Select elements $y_1 \in W_{J_1}, \, y_2 \in W_{J_2}, \, \dots, \, y_n \in W_{J_n}$ and infinite reduced words ${\bf s}_0 \in \mathscr{W}_{J_0}^{\infty}, \, {\bf s}_1 \in \mathscr{W}_{J_1}^{\infty}, \, \dots, \, {\bf s}_{n - 1} \in \mathscr{W}_{J_{n - 1}}^{\infty}$ satisfying the conditions below.
\begin{multline} \label{eq:chain_of_biconvex_set}
\emptyset = \nabla(J_0, w^{J_0}, 1_W) \subsetneqq \nabla(J_1, w^{J_1}, y_1) \subsetneqq \\
\dots \subsetneqq \nabla(J_n, w^{J_n}, y_n) = \Delta(w, -),
\end{multline}
\begin{multline} \label{eq:accumulate_infinite_reduced_words}
\nabla(J_i, w^{J_i}, y_i) = \nabla(J_{i - 1}, w^{J_{i - 1}}, y_{i - 1}) \amalg w^{J_{i - 1}} y_{i - 1} \Phi_{J_{i - 1}}({\bf s}_{i - 1}) \\
\quad (i = 1, 2, \dots, n).
\end{multline}

\item Then every root $\alpha \in \Delta(w, -)$ can be uniquely written as
\begin{equation} \label{eq:presentation_of_real}
\alpha = w^{J_{k - 1}} y_{k - 1} \phi_{{\bf s}_{k - 1}}(p) \quad (1 \leq k \leq n, \, p \in \Z_{\geq 1}).
\end{equation}
Using this expression we define a total order $\leq$ on $\Delta(w, -)$ by
\begin{multline}
w^{J_{k - 1}} y_{k - 1} \left(\phi_{{\bf s}_{k - 1}}(p)\right) \leq w^{J_{l - 1}} y_{l - 1} \left(\phi_{{\bf s}_{l - 1}}(q) \right) \\
\stackrel{\mathrm{def}}{\Leftrightarrow} (k < l) \text{or} (k = l, \, p \leq q) \quad (k, l, p, q \in \Z_{\geq 1}; \, k, l \leq n).
\end{multline}
Then $\leq$ is well-ordered and its ordinal number is $n \omega$, so that this well-order $\leq$ is called \textbf{$n$-low type}.
\end{enumerate}

Using this procedure, we construct two convex orders $\leq_-, \leq_+$ on $\Delta(w, -)$, $\Delta(w, +) = \Delta(ww_{\circ}, -)$ respectively. The parameters used in the procedure can be chosen independently between $\leq_-$ and $\leq_+$. We also set a total order $\leq_0$ on $\Delta_+^{\mathrm{im}}$ arbitrarily. Then we define a total order $\leq$ on whole $\Delta_+$ as follows.
\begin{align*}
\alpha \in \Delta(w, -), \, \beta \in \Delta_+^{\mathrm{im}}, \, \gamma &\in \Delta(w, +) \Rightarrow \alpha < \beta < \gamma; \\
\alpha \leq \alpha' \stackrel{\mathrm{def}}{\Leftrightarrow} \alpha \leq_- \alpha' \quad (\alpha, \alpha' \in \Delta(w, -)); &\quad \beta \leq \beta' \stackrel{\mathrm{def}}{\Leftrightarrow} \beta \leq_0 \beta' \quad (\beta, \beta' \in  \Delta_+^{\mathrm{im}}); \\
\gamma \leq \gamma' \stackrel{\mathrm{def}}{\Leftrightarrow} \gamma' \leq_+ \gamma &\quad (\gamma, \gamma' \in \Delta(w, +)).
\end{align*}
Notice that $\leq_+$ needs to be reversed, and therefore whole $\leq$ is not well-ordered.

\begin{thm}~\cite[Theorem 7.9, Corollary 7.10]{Ito1} The total order $\leq$ on $\Delta_+$ constructed above is convex, and any convex order on $\Delta_+$ can be constructed by the above procedure.
\end{thm}

\subsection{Convex bases of $U_q^+$ constructed by convex orders}

When $\mathfrak{g}$ is of finite type, it is known that $U_q^+$ has canonical bases, which can be described concretely by using braid group action on $U_q(\mathfrak{g})$ and correspond to each reduced expressions of longest element $w_{\circ}$ of Weyl group $W$~\cite{Lus}. In the affine type case, however, a couple of difficulties arise to construct basis of $U_q^+$ due to the absence of longest element of $W$ and the existence of imaginary roots. These problems are solved by constructing certain elements correspond to imaginary roots, using extended braid group action on $U_q(\mathfrak{g})$, which is proposed by Beck~\cite{Beck}. Then Ito generalized this construction to general convex orders~\cite{Ito2}. In this subsection, we summarize the construction of PBW type bases of $U_q^+$ from convex orders. We first introduce the notion of  convex basis, which is a PBW type basis with convexity property.

\begin{dfn}
Let $U$ be a $\Q(q)$-algebra, $\Lambda \subset U$ be a subset, and $\leq$ be a total order on $\Lambda$. For every subset $\Sigma \subset \Lambda$, the set of increasing monomials consist of the elements in $\Sigma$ is denoted by
\begin{equation*}
\mathscr{E}_{<} (\Sigma) := \left\{\,E_{\lambda_1} E_{\lambda_2} \dots E_{\lambda_m} \relmiddle| E_{\lambda_k} \in \Sigma, \, E_{\lambda_1} \leq E_{\lambda_2} \leq \dots \leq E_{\lambda_m}\,\right\} \subset U.
\end{equation*}
We call a subset $I \subset \Lambda$ \textbf{interval} if $I = \Lambda$, or $I$ coincide with one of $(x, *)$, $[x, *)$, $(*, y)$, $(*, y]$, $(x, y)$, $[x, y)$, $(x, y]$, $[x, y]$ for some $x, y \in \Lambda$, where $(x, *) := \left\{\, \lambda \in \Lambda \relmiddle| x < \lambda \,\right\}$, $[x, y) := \left\{\, \lambda \in \Lambda \relmiddle| x \leq \lambda < y\,\right\}$ and so on.

$\mathscr{E}_{<} (\Lambda)$ is called \textbf{convex basis} of $U$ if it has following properties:
\begin{enumerate}
\item $\mathscr{E}_{<} (\Lambda)$ is a basis of $U$ as $\Q(q)$-linear space.
\item For every interval $I \subset \Lambda$ with respect to given order $\leq$, let $U_I$ denote the $\Q(q)$-subalgebra of $U$ generated by $I$. Then $\mathscr{E}_{<} (I)$ is a basis of $U_I$ as $\Q(q)$-linear space.
\end{enumerate}
\end{dfn}

It is known that one can construct convex bases for quantum algebra $U_q(\mathfrak{g})$ by using the braid group action on $U_q(\mathfrak{g})$, which is given explicitly by the following fundamental result.

\begin{thm} ~\cite[Chap. 37, 39]{Lus} \label{thm:braid} There exists unique $\Q(q)$-algebra automorphism $T_i \in \mathrm{Aut} \, U_q(\mathfrak{g})$ $(i = 0, 1, \dots, \ell)$ satisfying 
\begin{align}
T_i(E_i) &= -F_i K_i, \quad T_i(F_i) = -K_i^{-1} E_i, \quad T_i(K_\lambda) = K_{s_i(\lambda)} \quad (\lambda \in P), \\
T_i(E_j) &= \frac{1}{[-a_{ij}]_{q_i} !} \sum_{k = 0}^{-a_{ij}} (-1)^{k} q_i^{-k} {\genfrac{[}{]}{0pt}{}{-a_{ij}}{k}}_{q_i} E_i^{-a_{ij} - k} E_j E_i^{k}, \label{eq:braid_action} \\
T_i(F_j) &= \frac{1}{[-a_{ij}]_{q_i} !} \sum_{k = 0}^{-a_{ij}} (-1)^{k} q_i^{k} {\genfrac{[}{]}{0pt}{}{-a_{ij}}{k}}_{q_i} F_i^{k} F_j F_i^{-a_{ij} - k} \quad (j \neq i).
\end{align}
Moreover, the automorphisms $T_i$ satisfy the braid relation
\begin{equation}
\overbrace{T_i T_j T_i \dots }^{m(i, j)} = \overbrace{T_j T_i T_j \dots }^{m(i, j)} \quad (i \neq j, \, m(i, j) \neq \infty), \label{eq:braid_rel}
\end{equation}
where $m(i, j) \in \Z_{\geq 1} \cup \{\infty\}$ is the order of $s_i s_j$ in the Weyl group.
\end{thm}

Recall that the braid group $\mathcal{B}$ associated with the Weyl group $W$ is defined by generators $T_i$ and relation (\ref{eq:braid_rel}). It is well-known that $\mathcal{B}$ has the following  property with respect to reduced expressions in $W$.
\begin{prop}
\label{prop:lift_to_braid}

Let $w = s_{i_1} s_{i_2} \dots s_{i_q} = s_{j_1} s_{j_2} \dots s_{j_q}$ are two reduced expressions of $w \in W$. Then $T_{i_1} T_{i_2} \dots T_{i_q} = T_{j_1} T_{j_2} \dots T_{j_q} \in \mathcal{B}$. Therefore a map
\begin{equation} \label{eq:lifting}
f: W \rightarrow \mathcal{B}; \quad f(w) :=  T_{i_1} T_{i_2} \dots T_{i_{q - 1}} (E_{i_q}) \in \mathcal{B} \quad (w = s_{i_1} s_{i_2} \dots s_{i_q} \, \text{: reduced})
\end{equation}
is well-defined and a section of canonical surjection $\pi: \mathcal{B} \rightarrow W$; $T_i \mapsto s_i$.
\end{prop}
Thus we can define the action of $w \in W$ on $U_q(\mathfrak{g})$ by $T_w := T_{f(w)}$, where $T_b$ denotes the action of $b \in \mathcal{B}$ on $U_q(\mathfrak{g})$.

When $\mathfrak{g}$ is of finite type, set $X_q := T_{i_1} T_{i_2} \dots T_{i_{q - 1}} (E_{i_q})$ $(q = 1, 2, \dots, N)$ where $w_{\circ} = s_{i_1} s_{i_2} \dots s_{i_N}$ is a reduced expression of the longest element. Then it is known that increasing monomials $X_1^{k_1} X_2^{k_2} \dots X_N^{k_N}$ constitutes a convex basis of $U_q^+$. $X_q$ has weight $\beta_q := s_{i_1} s_{i_2} \dots s_{i_{q - 1}} (\alpha_{i_q})$ and called a \textbf{root vector} associated to the root $\beta_q$. Root vectors depend on reduced expression of $w_{\circ}$.

Reduced expression of $w_{\circ} = s_{i_1} s_{i_2} \dots s_{i_N}$ induces a convex order $\beta_1 < \beta_2 < \dots < \beta_N$. This is because if $\alpha_{i_k} + s_{i_k} s_{i_{k + 1}} \dots s_{i_{l - 1}} (\alpha_{i_l}) = s_{i_{k}} \dots s_{i_{m - 1}} (\alpha_{i_m})$ and suppose that $l < m$, then applying $s_{i_{m - 1}} s_{i_{m - 2}} \dots s_{i_k}$ both sides yields $\alpha_{i_m} \in \Delta_-$, which is absurd. Thus the given order has the convexity property. Conversely, let $w \in W$ and $\beta_1 < \beta_2 < \dots < \beta_k$ be a convex order on $\Phi(w) := w \Delta_- \cap \Delta_+$. Then $\beta_1$ must be a simple root $\alpha_{i_1}$. To see this, we suppose that $\beta_1$ is not simple. Then $\beta_1$ can be written as a sum of two positive roots, and at least one of them belong to $\Phi(w)$ due to the biconvexity of $\Phi(w)$. This contradicts the minimality of $\beta_1$, and now we conclude $\beta_1 = \alpha_{i_1}$. Since the action of $W$ preserve the addition of roots, $s_{i_1} (\beta_2) < s_{i_1} (\beta_3) < \dots < s_{i_1} (\beta_N)$ is a convex order on $\Phi(s_{i_1} w)$. By induction on the length of $w$, we can construct a reduced expression of $w$ from given convex order. Therefore, each convex order on $\Delta_+$ generates a reduced expression of $w_{\circ}$. These correspondences between convex orders and reduced expression of $w_{\circ}$ is clearly inverses of each other, the correspondences are one-to-one.

Using the correspondences above, we want to construct bases for $U_q^+$ from convex orders on $\Delta_+$. To extend the construction for affine case, we have to deal with several problems such as the definition of root vectors when given convex order has multiple lows, existence of imaginary roots, which are unreachable from simple roots by only using the braid group action. These problems have already been solved by Beck~\cite{Beck} and Ito~\cite{Ito2}.

Before introducing their construction, we need to extend the affine Weyl group properly. We now return to consider the case when $\mathfrak{g}$ is untwisted affine Lie algebra of type $X_{\ell}^{(1)}$. The linear map $t_{\lambda} \in \mathrm{End}\, \mathfrak{h}^*$, called the translation by $\lambda \in \mathring{\mathfrak{h}}^*$, is defined by
\begin{equation}
t_{\lambda}(\mu) := \mu + (\mu, \delta) \lambda - \left\{\frac{1}{2}(\lambda, \lambda) (\mu, \delta) + (\mu, \lambda) \right\} \delta \quad (\mu \in \mathfrak{h}^*),
\end{equation}
where $\delta \in \Delta_+^{\mathrm{im}}$ is the null root. Let $T := \left\{\, t_{\nu(\Check{\alpha})} \relmiddle| \Check{\alpha} \in \Check{\mathring{Q}}\,\right\}$ be the group of translations, where $\nu: \mathfrak{h} \rightarrow \mathfrak{h}^*$ is the canonical isometry and $\Check{\mathring{Q}} \subset \mathfrak{h}$ is the coroot lattice. Then it is well-known that $W = \mathring{W} \ltimes T$~\cite{Kac}. Note that in general $T$ does not contain translation $t_{\varepsilon_i}$ by fundamental coweight $\varepsilon_i \in \mathring{\mathfrak{h}}^*$,  which is characterized by $(\varepsilon_i, \alpha_j) = \delta_{ij}$ for $i, j = 1, \dots, \ell$. We extend the affine Weyl group $W$
by appending the translations $t_{\varepsilon_i}$. Let $\widehat{W}$ denote the subgroup of $\mathrm{GL} (\mathfrak{h}'^*)$ generated by $W$ and $\maprestrict{t_{\varepsilon_i}}{\mathfrak{h}'^*}$ $(i \in \mathring{I})$, where $\mathfrak{h}'^* := \bigoplus_{i = 0}^{\ell} \C \alpha_i \subset \mathfrak{h}^*$. In fact, the extended Weyl group $\widehat{W}$ coincides with a semidirect product of $W$ and a subgroup of the Dynkin automorphism group.

\begin{prop} ~\cite{IM} \cite[Proposition 2.1]{Ito2} \label{prop:extW} Let $\mathring{I}_* := \left\{\,j \in \mathring{I} \relmiddle| (\varepsilon_j, \theta_{\mathring{I}}) = 1\,\right\}$, and $\rho_{\mathring{I}j} := t_{\varepsilon_j} w_{\circ j} w_{\circ}$ for each $j \in \mathring{I}_*$, where $w_{\circ}, w_{\circ j}$ are the longest elements of $\mathring{W}$, $\mathring{W}_{\mathring{I} \setminus \{j\}}$ respectively. Then there exists an automorphism $\rho$ of the Dynkin diagram of type $X_{\ell}^{(1)}$ such that $\rho_{\mathring{I}j} (\alpha_i) = \alpha_{\rho(i)}$ for all $i = 0, 1, \dots, \ell$. The correspondence $j \mapsto \rho$ is one-to-one. Moreover, $\Omega := \left\{\,\rho_{\mathring{I}j} \relmiddle| j \in \mathring{I}_*\,\right\} \amalg \{\mathrm{id}_{\mathfrak{h}'^*}\}$ forms a subgroup of $\mathrm{GL} (\mathfrak{h}'^*)$, and
\begin{equation}
\widehat{W} = \Omega \ltimes W, \label{eq:extW}
\end{equation}
where $\rho_{\mathring{I}j} \in \Omega$ acts on $W$ by $\rho_{\mathring{I}j}.s_{i} := s_{\rho(i)}$.
\end{prop}

We define the length of $w \in \widehat{W}$ by
\begin{equation} \label{eq:length_in_extW}
\ell(w) := \ell(u) = \ell_{\mathring{I}} (u),
\end{equation}
where we use the decomposition $w = \rho u$ $(\rho \in \Omega, u \in W)$ given by (\ref{eq:extW}). Recall that $W = W_{\mathring{I}}$. We can also consider reduced expressions of $w \in \widehat{W}$. Namely, we call an expression $w = t_1 t_2 \dots t_m \in \widehat{W}$, $t_i \in S \amalg \Omega$ \textbf{reduced} if the sequence of integer $\ell(t_1), \ell(t_1 t_2), \dots, \ell(w)$ is increasing. Thus every element of $\Omega$ has length $0$, and reduced expressions of $w \in \widehat{W}$ may have different number of factors but the number of factors which belong to $S$ must coincide with the length of $w$.

Dynkin automorphism $\rho$ acts on the subalgebra $U_q'(\mathfrak{g}) := \langle E_i, F_i, K_{\alpha_i}^{\pm 1} \rangle \subset U_q(\mathfrak{g})$ as algebra automorphism by permuting indices: $E_i \mapsto E_{\rho(i)}$, $F_i \mapsto F_{\rho(i)}$, $K_{\alpha_i} \mapsto K_{\alpha_{\rho(i)}}$. Thus we have an action of the \textbf{extended braid group} $\widehat{\mathcal{B}} := \Omega \ltimes \mathcal{B}$ on $U_q'(\mathfrak{g})$ by extending the braid group action of Theorem \ref{thm:braid}. Proposition \ref{prop:lift_to_braid} also holds for $\widehat{W}$ and $\widehat{\mathcal{B}}$, and therefore every $w \in \widehat{W}$ has an action $T_w$ on $U_q'(\mathfrak{g})$.

We will define the root vectors associated to real roots by lifting the expression (\ref{eq:presentation_of_real}) to quantum algebra $U_q(\mathfrak{g})$, in which process simple reflection $s_i$ is replaced by $T_i$ and simple root $\alpha_i$ is replaced by $E_i$. In this lifting process, we also have to specify appropriate alternatives for $\delta - \theta_{J_c} \in \Pi_{J_c}$ and $s_{\delta - \theta_{J_c}} \in S_{J_c}$, where $J_c \subset \mathring{I}$ is a connected subdiagram. The simple root vector $E_{\delta - \theta_{J_c}}$ is in fact uniquely determined due to the following lemma.
\begin{lem}~\cite[Lemma 5.1]{Ito2} Let $\varepsilon \in \mathring{\Delta}_+$ and suppose that ${\bf s} := s_{i_1} s_{i_2} \dots s_{i_m}$ is a reduced expression in $W$ satisfying $\delta - \varepsilon = s_{i_1} s_{i_2} \dots s_{i_{m - 1}} (\alpha_{i_m})$ and $\Phi({\bf s}) \subset \Delta(1, -)$. Such a ${\bf s}$ exists and
\begin{equation} \label{eq:sub_simple_root_vector}
E_{\delta - \varepsilon} := T_{i_1} T_{i_2} \dots T_{i_{m - 1}} (E_{i_m}) \in U_q^+
\end{equation}
is independent of the choice of ${\bf s}$.
\end{lem}

The appropriate alternative for $s_{\delta - \theta_{J_c}}$ is given by somewhat technical manner.
\begin{dfn}~\cite[Definition 3.4]{Ito2} First, we fix an index $j_c \in J_c$ satisfying $(\varepsilon_{j_c}, \theta_{J_c}) = 1$ for every nonempty connected subdiagram $J_c \subset \mathring{I}$. Then we define a map $\widehat{\cdot}: S_J \rightarrow \widehat{W}$ by
\begin{equation} \label{eq:simple_ref_alternative}
\widehat{s_j} := s_j \quad (j \in J), \quad \widehat{s_{\delta - \theta_{J_c}}} := \left(t_{\varepsilon_{j_c}}\right)^{J_c} s_{\Bar{j_c}} \left(t_{\varepsilon_{\Bar{j_c}}}\right)^{J_c},
\end{equation}
where $\Bar{j_c} \in \mathring{I}$ is the unique index which satisfies $w_{\circ} (\alpha_{\Bar{j_c}}) = -\alpha_{j_c}$. We also define the extended map $\widehat{\cdot}: W_J \rightarrow \widehat{W}$ simply by $\widehat{w} := \widehat{t_1} \widehat{t_2} \dots \widehat{t_m}$ when $w = t_1 t_2 \dots t_m$ is a reduced expression in $W_J$.
\end{dfn}

Now, we can describe the construction of root vectors for affine case. 
\begin{dfn}~\cite[Theorem 8.4]{Ito2} Suppose that $\leq$ is a convex order on $\Delta_+$. Let $\alpha = w^{J_{k - 1}} y_{k - 1} \phi_{{\bf s}_{k - 1}}(p)$ be the expression (\ref{eq:presentation_of_real}) of a positive real root $\alpha$ determined by $\leq$. Then root vector $E_{\leq, \alpha} \in U_{\alpha}$ associated to $\alpha$ is defined by
\begin{multline}
E_{\leq, \alpha} \\
:= \begin{cases}
T_{w^{J_{k - 1}}} {T}_{\widehat{y_{k - 1}}} {T}_{\widehat{{\bf s}_{k - 1} (1)}} {T}_{\widehat{{\bf s}_{k - 1}(2)}} \dots {T}_{\widehat{{\bf s}_{k - 1} (p - 1)}} (E_{{\bf s}_{k - 1} (p)}) & (\alpha \in \Delta(w, -)) \\
\Psi T_{w^{J_{k - 1}}} {T}_{\widehat{y_{k - 1}}} {T}_{\widehat{{\bf s}_{k - 1} (1)}} {T}_{\widehat{{\bf s}_{k - 1}(2)}} \dots {T}_{\widehat{{\bf s}_{k - 1} (p - 1)}} (E_{{\bf s}_{k - 1} (p)}) & (\alpha \in \Delta(w, +))
\end{cases},
\end{multline}
where $E_{s_i} := E_i$, and $\Psi: U_q(\mathfrak{g}) \rightarrow U_q(\mathfrak{g})$ is the anti-automorphism of $\Q(q)$-algebra defined by 
$\Psi(E_i) := E_i$, $\Psi(F_i) := F_i$, $\Psi(K_\lambda) := K_{\lambda}^{-1}$. 
\end{dfn}

The root vectors for imaginary roots is constructed using the action of extended braid group, which contains coweight lattice~\cite{Beck}. Since each imaginary root has multiplicity $\ell$ in affine Lie algebra $\mathfrak{g}$ of type $X_{\ell}^{(1)}$, we will construct as many number of root vectors as the multiplicity. The construction is rather technical and we proceed step-by-step.

First, we introduce weight vectors $\mathcal{E}_{n\delta - \alpha_i}$ $(i \in \mathring{I})$, which is independent of convex order.
\begin{equation}
\mathcal{E}_{n\delta - \alpha_i} := T_{\varepsilon_i}^n T_i^{-1} (E_i) \quad (n \in \Z_{\geq 1}, i \in \mathring{I}),
\end{equation}
where $T_{\varepsilon_i} := T_{t_{\varepsilon_i}} \in \mathrm{Aut} \, U'_q(\mathfrak{g})$ was defined via the extended braid group action and lifting a reduced expression of $t_{\varepsilon_i} \in \widehat{W}$ to $\widehat{\mathcal{B}}$ by (\ref{eq:lifting}). Then we set
\begin{equation}
\varphi_{i, n} := [\mathcal{E}_{n\delta - \alpha_i}, E_i]_q = \mathcal{E}_{n\delta - \alpha_i}E_i - q_i^{-2} E_i \mathcal{E}_{n\delta - \alpha_i} \quad (n \in \Z_{\geq 1}, i \in \mathring{I}).
\end{equation}
Despite of these $\varphi_{i, n}$ have weight $n\delta \in \Delta_+^{\mathrm{im}}$, $\varphi_{i, n}$ are not yet suitable for imaginary root vectors. The genuine imaginary root vectors are constructed by modifying $\varphi_{i, n}$ through the following technical procedure. For every $i \in \mathring{I}$, let
\begin{equation}
\label{eq:gen_func_imaginary}
\varphi_i(z) := (q_i - q_i^{-1}) \sum_{n = 1}^{\infty} \varphi_{i, n} z^n \in U_q^+[[z]]
\end{equation}
be the generating function of $\varphi_{i, n}$. $U_q^+[[z]]$ has topological algebra structure by declaring that $z$ is central and $U_q^+[[z]]$ has $z$-adic topology. Then \textbf{imaginary root vectors} $I_{i, n} \in U_q^{+}$ are defined as the coefficient of the function
\begin{equation}
\label{eq:definition_imaginary}
I_i(z) := \log(1 + \varphi_i(z)) = (q_i - q_i^{-1}) \sum_{n = 1}^{\infty} I_{i, n} z^n,
\end{equation}
where the logarithm is defined by $\log (1 + x) := \sum_{m = 1}^{\infty} \frac{(-1)^{m - 1}}{m} x^m$.

It is shown that the these root vectors constitute convex bases for positive part of quantum affine algebra $U_q^+$.

\begin{thm}~\cite[Theorem 8.6]{Ito2} Let $\leq$ be a convex order on positive roots $\Delta_+$ of untwisted affine root system, and let $w \in \mathring{W}$ be the parameter determined by the decompositon (\ref{eq:positive_root_partition}) of $\Delta_+$ in accordance with the given convex order $\leq$. Let
\begin{equation}
\Lambda := \left\{\, E_{\leq, \alpha} \relmiddle| \alpha \in \Delta_+^{\mathrm{re}} \,\right\} \amalg \left\{\, T_w (I_{i, m}) \relmiddle| m \in \Z_{\geq 1}; i = 1, 2, \dots, \ell \,\right\}
\end{equation}
denote the set of root vectors constructed above, and we set the order on $\Lambda$ by using given order $\leq$ and
\begin{equation*}
T_w(I_{i, m}) \leq T_w(I_{j, m'}) \Leftrightarrow (m \leq m') \, \text{or} \, (m = m', i \leq j).
\end{equation*}
Then increasing monomials $\mathscr{E}_{<}(\Lambda)$ constitute a convex basis of $U_q^+$.
\end{thm}

Once a convex basis of $U_q^+$ is constructed, we also obtain the one for $U_q^-$ through Chevalley involution $\Omega: U_q^+ \rightarrow U_q^-;$ $E_i \mapsto F_i, q \mapsto q^{-1}$, which is anti-automorphism of $\Q$-algebra.

\subsection{Product formula for the quasi-universal R-matrix}

The convex bases for quantum affine algebra enable explicit construction of quasi-universal R-matrix. By applying Drinfeld's quantum double construction, Ito obtained the product formula for the quasi-universal R-matrix~\cite{Ito2}. Since the quasi-universal R-matrix does not lie in the algebraic tensor product $U_q(\mathfrak{g}) \otimes U_q(\mathfrak{g})$, we have to give appropriate topology on $U_q(\mathfrak{g}) \otimes U_q(\mathfrak{g})$ and complete it.

First, we set the gradation of $U_q(\mathfrak{g}) \otimes U_q(\mathfrak{g})$ by
\begin{equation}
(U_q \otimes U_q)_h := \bigoplus_{\substack{\mu, \nu \in Q_+\\ \mathrm{ht} (\mu + \nu) = h}} (U_q^- U_q^0 \otimes U_q^- U_q^0) \cdot U_{\mu}^+ \otimes U_{\nu}^+ \subset U_q(\mathfrak{g}) \otimes U_q(\mathfrak{g}) \quad (h \in \Z_{\geq 0}),
\end{equation}
that is, we only count the weight of positive part with respect to the triangular decomposition (\ref{eq:tri_decomp}). Then we set a topology which is generated by the subsets of the form
\begin{equation}
x + \bigoplus_{h = k}^{\infty} (U_q \otimes U_q)_h \quad (x \in U_q(\mathfrak{g}), k \in \Z_{\geq 0}).
\end{equation}
In short, we give $U_q(\mathfrak{g}) \otimes U_q(\mathfrak{g})$ linear topology. Let
\begin{equation}\label{eq:topology_on_U_q}
\widehat{U_q} \widehat{\otimes} \widehat{U_q} := \underset{k \geq 0}{\mathrm{\projlim}} \left(U_q\left(\mathfrak{g}\right) \otimes U_q\left(\mathfrak{g}\right) / \bigoplus_{h = k}^{\infty} \left(U_q \otimes U_q\right)_h \right)
\end{equation}
be the completion of $U_q(\mathfrak{g}) \otimes U_q(\mathfrak{g})$, and $U_q^+ \widehat{\otimes} U_q^- \subset \widehat{U_q} \widehat{\otimes} \widehat{U_q}$ denote the closure of $U_q^+ \otimes U_q^-$. The algebra structure of $U_q(\mathfrak{g}) \otimes U_q(\mathfrak{g})$ extends uniquely onto $\widehat{U_q} \widehat{\otimes} \widehat{U_q}$.

\begin{dfn} ~\cite[4.1.2]{Lus} Let $\Upsilon \in \mathrm{Aut}\, U_q(\mathfrak{g})$ be the $\Q$-algebra automorphism determined by 
\begin{equation*}
\Upsilon(E_i) := E_i, \quad \Upsilon(F_i) := F_i, \quad \Upsilon(K_{\lambda}) := K_{\lambda}^{-1}, \quad \Upsilon(q) := q^{-1},
\end{equation*}
and set $\Bar{\Delta} := (\Upsilon \otimes \Upsilon) \circ \Delta \circ \Upsilon$. The \textbf{quasi-universal R-matrix} of $U_q(\mathfrak{g})$ is the \textbf{unique} element $\Theta \in \widehat{U_q} \widehat{\otimes} \widehat{U_q}$ satisfying
\begin{enumerate}
\item $\Theta \cdot \Bar{\Delta}^{\mathrm{op}} (u) = \Delta^{\mathrm{op}} (u) \cdot \Theta \quad (\forall u \in U_q(\mathfrak{g}))$.
\item $\Theta_0 = 1 \otimes 1$,
\end{enumerate}
where $\Theta_0 \in (U_q \otimes U_q)_0$ is the image of $\Theta$ by the canonical projection
\begin{equation*}
\widehat{U_q} \widehat{\otimes} \widehat{U_q} \twoheadrightarrow U_q\left(\mathfrak{g}\right) \otimes U_q\left(\mathfrak{g}\right) / \bigoplus_{h = 1}^{\infty} \left(U_q \otimes U_q\right)_h \cong (U_q \otimes U_q)_0,
\end{equation*}
and $f^{\mathrm{op}} (u) := \sum y_i \otimes x_i$ if $f(u) = \sum x_i \otimes y_i$.
\end{dfn}

The uniqueness of $\Theta$ will be the core of the proof of identites. Finally, we introduce the product formula for quasi-universal R-matrix.
\begin{thm}~\cite{Ito0}~\cite{Ito2} Let $\leq$ be a convex order on $\Delta_+$ of affine root system, and $E_{\leq, \alpha}$, $I_{i, n}$ denote the root vectors constructed above. For every $i, j \in \mathring{I}$ and positive integer $n$, let 
\begin{equation} \label{eq:bijn}
b_{i, j; n} := \mathrm{sgn}(a_{ij})^n \frac{[a_{ij} n]_{q_i}}{n(q_j^{-1} - q_j)}, \quad \mathrm{sgn}(x) := \begin{cases}1 & x > 0 \\ 0 & x = 0 \\ -1 & x < 0\end{cases}.
\end{equation} 
Let $(c_{i, j; n})_{i, j = 1}^{\ell} \in \mathrm{Mat} (\Q(q), \ell)$ denote the inverse matrix of $(b_{i, j;n})_{i, j = 1}^{\ell}$.

We also set
\begin{align*}
F_{\leq, \alpha} := \Omega(E_{\leq, \alpha}) \in U_{-\alpha} \quad (\alpha \in \Delta_+^{\mathrm{re}}), &\quad J_{i, n} := \Omega(I_{i, n}) \in U_{-n\delta}, \\
\exp_q(x) := \sum_{n = 0}^{\infty} \frac{q^{-\frac{1}{2} n(n-1)}}{[n]_q !} x^n, &\quad q_{\alpha} := q^{\frac{1}{2} (\alpha, \alpha)} \quad (\alpha \in \Delta),
\end{align*}
\begin{align}
S_n &:= \sum_{i, j \in \mathring{I}} c_{j, i; n} I_{i, n} \otimes J_{j, n} \quad \in U_q^+ \otimes U_q^-, \label{eq:def_Sn}\\
\Theta_{\leq, \alpha} &:= \begin{cases}
\exp_{q_{\alpha}} \left\{(q_{\alpha}^{-1} - q_{\alpha}) E_{\leq, \alpha} \otimes F_{\leq, \alpha}\right\} & \alpha \in \Delta_+^{\mathrm{re}} \\
\exp \left\{T_w \otimes T_w (S_n)\right\} & \alpha = n\delta \quad (n = 1, 2, \dots)
\end{cases} \notag.
\end{align}
Then the quasi-universal R-matrix $\Theta$ has the product presentation
\begin{equation} \label{eq:prod_formula}
\Theta = \prod_{\alpha \in \Delta_+}^{>} \Theta_{\leq, \alpha} \quad \in U_q^+ \widehat{\otimes} U_q^-,
\end{equation}
where $\prod_{\alpha \in \Delta_+}^{>} X_{\alpha}$ means that if  $\alpha < \beta$, the order of multiplication is $X_{\beta} X_{\alpha}$. In short, the order of multiplication is reverse to given convex order.
\end{thm}

\section{Explicit presentation of root vectors using $q$-bracket} \label{sec:2.5}

We will construct quantum dilogarithm identities by using various presentation (\ref{eq:prod_formula}) of quasi-universal R-matrix $\Theta$, taking advantage of the uniqueness of $\Theta$. However, to obtain specific identities, we have to calculate root vectors explicitly, which is described by braid group action (Theorem \ref{thm:braid}). In this section, we show that in general quantum algebra $U_q(\mathfrak{g})$ of symmetrizable Kac-Moody algebra $\mathfrak{g}$, the element $T_w(E_i) \in U_q^+$ $(w \in W)$ can be written as "$q$-commutator monomial", that is, finite application of $q$-bracket on the generators $E_i$. We also construct concrete algorithm for getting explicit presentation of $T_w(E_i)$ as $q$-commutator monomial, which enables us direct computation.

Let $\mathfrak{g}$ be a symmetrizable Kac-Moody algebra of rank $n$.
\begin{dfn}
For every subsets $A, B \subset U_q(\mathfrak{g})$, let
\begin{equation*}
[A, B]_q := \left\{\, [x, y]_q \relmiddle| x \in A, y \in B\,\right\} \subset U_q(\mathfrak{g}).
\end{equation*}
We define subsets $P_k \subset U_q(\mathfrak{g})$ inductively by
\begin{equation*}
P_0 := \{E_1, E_2, \dots, E_n\}, \quad P_{k + 1} := \bigcup_{i, j = 0}^{k} [P_i, P_j]_q \quad (k \in \Z_{\geq 0}).
\end{equation*}
We call the elements of the form $cM \in U_q(\mathfrak{g})$ for some $c \in \Q(q), M \in \bigcup_{k = 0}^{\infty} P_k$ \textbf{$q$-commutator monomial}.
\end{dfn}

Our claim is that $T_w(E_i) \in U_q^+$ is $q$-commutator monomial for all $i = 1, \dots, n$ and $w \in W$. To prove it, several formulas have to be prepared. First we recall
\begin{equation} \label{eq:braid_action2}
T_i (E_j) = \frac{1}{[-a_{ij}]_{q_i} !} \overbrace{[E_i, [E_i, \cdots, [E_i, }^{-a_{ij}} E_j]_q ]_q \dots ]_q \quad (i \neq j)
\end{equation}
by definition of the braid group action, and
\begin{equation} \label{eq:braid_preserve_q_comm}
T_i \left(\left[x, y\right]_q\right) = \left[T_i \left(x\right), T_i \left(y\right)\right]_q \quad (x, y \in U_q(\mathfrak{g}), \, i = 1, 2, \dots, n),
\end{equation}
since the Weyl group action preserves the invariant bilinear form. The basic process of calculation for $T_w(E_j)$ $(w \in W, 1 \leq j \leq n)$ is as follows: choose a reduced expression $w = s_{i_1} s_{i_2} \dots s_{i_m}$, and expand every $T_{i_k}$ of $T_w = T_{i_1} \dots T_{i_m}$ from the tail using (\ref{eq:braid_action2}) and (\ref{eq:braid_preserve_q_comm}). However, there is a problem that $T_k (E_k) = -F_k K_k$ may be appeared in the process of expansion. To resolve it, we use the following formula.

\begin{lem}\label{lem:q_adjoint_power_commutation}
For every $1 \leq i \neq j \leq n$ and positive integer $m$,
\begin{equation}\label{eq:q_adjoint_power_commutation}
\left[\left(\stackrel{\rightarrow}{\mathrm{ad}}\, E_i \right)^m (E_j), \, T_i(E_i)\right]_q = [m]_{q_i} [1 - a_{ij} - m]_{q_i} \left(\stackrel{\rightarrow}{\mathrm{ad}}\, E_i \right)^{m - 1} (E_j),
\end{equation}
where $\stackrel{\rightarrow}{\mathrm{ad}}\, x (y) := [x, y]_q$.
\end{lem}

\begin{proof}
By defining relation of $U_q(\mathfrak{g})$,
\begin{align*}
E_i K_i &=  q_i^{-2} K_i E_i, \quad E_j K_i = q_i^{-a_{ij}} K_i E_j, \\
E_i F_i &= F_i E_i + \frac{K_i - K_i^{-1}}{q_i - q_i^{-1}}, \quad E_j F_i = F_i E_j.
\end{align*}
Thus the commutation relations of $F := T_i(E_i) = -F_i K_i$ and $E_i, E_j$ are as follows.
\begin{equation*}
E_i F = q_i^{-2} F E_i - \frac{K_i^2 - 1}{q_i - q_i^{-1}}, \quad E_j F = q_i^{-a_{ij}} F E_j.
\end{equation*}
Since the weight of $F = -F_i K_i$ is $-\alpha_i$, we have
\begin{equation} \label{eq:m=0case}
[E_j, F]_q = E_j F - q_i^{-a_{ij}} F E_j = q_i^{-a_{ij}} F E_j - q_i^{-a_{ij}} F E_j = 0.
\end{equation}

Now, we begin the proof by induction on $m$. Suppose that (\ref{eq:q_adjoint_power_commutation}) holds for some positive integer $m$. Let
\begin{equation*}
C_m := [m]_{q_i} [1 - a_{ij} - m]_{q_i}, \quad X_m := \left(\stackrel{\rightarrow}{\mathrm{ad}}\, E_i \right)^m (E_j).
\end{equation*}
Since the weight of $X_m$ is $m \alpha_i + \alpha_j$, we have $X_m K_i = q_i^{-2m - a_{ij}} K_i X_m$. Then by the induction hypothesis,
\begin{equation*}
X_m F = q_i^{-2m - a_{ij}} F X_m + C_m X_{m - 1}.
\end{equation*}
Using these commutation relations, we obtain the equation for the case of $m + 1$.
\begin{align*}
&\left[X_{m + 1}, \, F\right]_q \\
&= \left[\left[E_i, X_m \right]_q, \, F\right]_q \\
&= \left[E_i X_m - q_i^{2m + a_{ij}} X_m E_i, \, F\right]_q \\
&= E_i X_m F - q_i^{-2(m + 1) - a_{ij}} F E_i X_m \\
&\qquad - q_i^{2m + a_{ij}} \left\{X_m E_i F - q_i^{-2(m + 1) - a_{ij}} F X_m E_i \right\} \\
&= E_i \left\{q_i^{-2m - a_{ij}} F X_m + C_m X_{m - 1}\right\} - q_i^{-2(m + 1) - a_{ij}} F E_i X_m \\
&\qquad -q_i^{2m + a_{ij}} \left\{X_m \left(q_i^{-2} F E_i - \frac{K_i^2 - 1}{q_i - q_i^{-1}} \right) -q_i^{-2(m + 1) - a_{ij}} F X_m E_i \right\} \\
&= q_i^{-2m - a_{ij}} \left(q_i^{-2} F E_i - \frac{K_i^2 - 1}{q_i - q_i^{-1}} \right) X_m  + C_m E_i X_{m - 1} - q_i^{2(m + 1) + a_{ij}} F E_i X_m \\
&\qquad -q_i^{2m + a_{ij}} q_i^{-2} \left(q_i^{-2m - a_{ij}}F X_m + C_m X_{m - 1} \right) E_i + q_i^{2m + a_{ij}} X_m \frac{K_i^2 - 1}{q_i - q_i^{-1}} \\
&\qquad \qquad + q_i^{-2} FX_m E_i \\
&= -q_i^{-2m - a_{ij}} \frac{K_i^2 - 1}{q_i - q_i^{-1}} X_m + C_m E_i X_{m - 1} \\
&\qquad - q_i^{2(m - 1) + a_{ij}} C_m X_{m - 1} E_i + q_i^{2m + a_{ij}} \frac{q_i^{2(-2m - a_{ij})}K_i^2 - 1}{q_i - q_i^{-1}} X_m \\
&= \frac{q_i^{-2m - a_{ij}} - q_i^{2m + a_{ij}}}{q_i - q_i^{-1}} X_m + C_m \left[E_i, X_{m - 1}\right]_q \\
&= \left(\left[-a_{ij} - 2m\right]_{q_i} + C_m\right) X_m.
\end{align*}
Thus we obtain the recursion formula
\begin{equation*}
C_{m + 1} = C_m + \left[-a_{ij} - 2m\right]_{q_i} \quad (m \geq 1).
\end{equation*}
It is easy to verify $C_m := [m]_{q_i} [1 - a_{ij} - m]_{q_i}$ satisfies this recurrence relation. Therefore, (\ref{eq:q_adjoint_power_commutation}) holds for $m + 1$.

For the case of $m = 1$, above calculation works if one uses (\ref{eq:m=0case}) in place of induction hypothesis and lets $C_0 := 0$.
\end{proof}

\begin{prop}\label{prop:q_pol_finite}
Suppose that the root subsystem spanned by $\alpha_i, \alpha_j$ $(i \neq j)$ is of finite type and $w = s_i s_j s_i s_j \dots$ is a reduced expression. Then $T_w(E_k)$ ($k = i$ if $\ell(w)$ is even, $k = j$ otherwise) is a $q$-commutator monomial consists of $E_i$ and $E_j$.
\end{prop}

\begin{proof}
Since length of reduced expression of the form $s_i s_j s_i s_j \dots$ is at most $5$ when finite type case, our task is just compute $T_w(E_k)$ directly for all cases. Using the formula (\ref{eq:q_adjoint_power_commutation}), the computation is easily accomplished. For example, when $a_{ij} = a_{ji} = -1$ we have
\begin{align*}
T_i T_j (E_i) &\stackrel{\text{(\ref{eq:braid_action2})}}{=} T_i \left([E_j, E_i]_q\right) \stackrel{\text{(\ref{eq:braid_preserve_q_comm})}}{=} \left[T_i(E_j), T_i(E_i)\right]_q \stackrel{\text{(\ref{eq:braid_action2})}}{=} \left[[E_i, E_j]_q, T_i(E_j)\right]_q \\
&\stackrel{\text{Lemma \ref{lem:q_adjoint_power_commutation}}}{=} [1]_{q_i} [1 - (-1) - 1]_{q_i} E_j = E_j.
\end{align*}
Thus we have a reduction formula
\begin{equation} \label{eq:reduction_A2}
T_i T_j (E_i) = E_j \quad \text{if} \  (a_{ij}, a_{ji}) = (-1, -1).
\end{equation}
\end{proof}

For infinite type case, we use the following formulas, which can be verified by direct computation using (\ref{eq:q_adjoint_power_commutation}).
\begin{lem}\label{lem:braid_rec_formulas}
For indices $i, j$ ($i \neq j$) and nonnegative integers $p, k$, let
\begin{align}
{\bf s}_{i,j;p} &:= \overbrace{\dots s_j s_i s_j s_i s_j}^{p}, \\
V_{i, j; p}^{(k)} &:= T_{{\bf s}_{i,j;p}} \left(\left(\stackrel{\leftarrow}{\mathrm{ad}}\, E_j\right)^k (E_i)\right) \\
&=  \overbrace{\dots T_i T_j T_i T_j }^{p} \left([[\dots [E_i, \overbrace{E_j]_q, E_j]_q, \dots, E_j}^{k}]_q\right).
\end{align}
These $V_{i, j; p}^{(k)}$ satisfy the following recurrence relations.
\begin{align}
V_{i, j; p + 1}^{(1)} &= \frac{1}{[-a_{ji} - 1]_{q_j} !} \left(\stackrel{\rightarrow}{\mathrm{ad}}\, E_{{\bf s}_{i, j; p + 1}}\right)^{-a_{ji} - 2} \left(V_{j, i; p}^{(1)}\right) \quad (a_{ji} \leq -2), \label{eq:V_ijp_formula1}\\
V_{i, j; p + 2}^{(1)} &= \frac{1}{[-a_{ji} - 1]_{q_j} !} \left(\stackrel{\rightarrow}{\mathrm{ad}}\, V_{i, j; p}^{(1)}\right)^{-a_{ji} - 3} \left(V_{i, j; p}^{(2)}\right) \quad (a_{ij} = -1, a_{ji} \leq -3), \label{eq:V_ijp_formula2}\\
V_{i, j; p + 2}^{(2)} &= \frac{[2]_{q_j}}{[-a_{ji} - 2]_{q_j} !} \left(\stackrel{\rightarrow}{\mathrm{ad}}\, V_{i, j; p}^{(1)}\right)^{-a_{ji} - 4} \left(V_{i, j; p}^{(2)}\right) \quad (a_{ij} = -1, a_{ji} \leq -4). \label{eq:V_ijp_formula3}
\end{align}

\begin{thm}\label{thm:q_comm}
For every $w \in W$ and index $j$ satisfying $w(\alpha_j) \in \Delta_+$, $T_w(E_j) \in U_q(\mathfrak{g})$ is a $q$-commutator monomial.
\end{thm}

\begin{proof}
The proof is by induction on $\ell(w)$. The case $\ell(w) = 1$ is immediate by (\ref{eq:braid_action2}). Suppose that there exists integer $m \geq 2$ such that $T_w(E_j)$ is a $q$-commutator monomial if  $w(\alpha_j) \in \Delta_+$ and $\ell(w) < m$. Let $w \in W$ satisfies $w(\alpha_j) \in \Delta_+$ and $\ell(w) = m$. Take a reduced expression of $w$ and let $s_i$ be its suffix. Then $i \neq j$ due to the assumption. Let $w^{\{i, j\}} \in W$ be the shortest element satisfying $w = w^{\{i, j\}} \dots s_j s_i s_j s_i$. Then $w^{\{i, j\}} (\alpha_i), w^{\{i, j\}} (\alpha_j) \in \Delta_+$ due to the minimality of $w^{\{i, j\}}$. By the induction hypothesis, $T_{w^{\{i, j\}}} (E_i)$ and $T_{w^{\{i, j\}}} (E_j)$ are $q$-commutator monomials. Thus using (\ref{eq:braid_preserve_q_comm}), the proof completes if $\dots T_j T_i (E_j)$ turns out to be a $q$-commutator monomial consists of only $E_i$ and $E_j$.

Let $\mathscr{F}_{ij}$ $(i \neq j)$ denotes the set of $q$-commutator monomials consist of only $E_i$ and $E_j$. We are going to prove that if ${\bf s}_{i,j;p}$ $(p \in \Z_{\geq 1})$ is a reduced expression, then $E_{{\bf s}_{i,j;p}} := \dots T_j T_i (E_j) \in \mathscr{F}_{ij}$ by induction on $p$. The cases when $p = 1, 2$ are immediate by (\ref{eq:braid_action2}). When $\alpha_i$ and $\alpha_j$ span a finite root system, ${\bf s}_{i,j;p}$ is reduced only for finitely many $p \in \Z_{\geq 1}$. Thus when $(a_{ij}, a_{ji}) = (-1, -1)$, $(-1, -2)$, $(-1, -3)$, $(-2, -1)$, $(-3, -1)$, we can verify $E_{{\bf s}_{i,j;p}} \in \mathscr{F}_{ij}$ by direct computation since there exists only finitely many cases. The computation is easily accomplished using the formula (\ref{eq:q_adjoint_power_commutation}). For example when $(a_{ij}, a_{ji}) = (-1, -1)$ and $p = 3$, we have
\begin{align*}
E_{{\bf s}_{i,j;3}} &:= T_j T_i (E_j) \stackrel{\text{(\ref{eq:braid_action2})}}{=} T_j \left([E_i, E_j]_q\right) \stackrel{\text{(\ref{eq:braid_preserve_q_comm})}}{=} \left[T_j(E_i), T_j(E_j)\right]_q \\
&\stackrel{\text{(\ref{eq:braid_action2})}}{=} \left[[E_j, E_i]_q, T_j(E_j)\right]_q \stackrel{\text{Lemma \ref{lem:q_adjoint_power_commutation}}}{=} [1]_{q_j} [1 - (-1) - 1]_{q_j} E_i = E_i.
\end{align*}

When $\alpha_i$ and $\alpha_j$ span infinite root system, then $a_{ij} a_{ji} \geq 4$. Now we suppose that $p \geq 2$ and $E_{{\bf s}_{i,j;r}} \in \mathscr{F}_{ij}$ for all $r \leq p$. First, $E_{{\bf s}_{i, j; p + 1}}$ can be written as follows.
\begin{equation}\label{eq:expansion}
E_{{\bf s}_{i, j; p + 1}} = \frac{1}{[-a_{ij}]_{q_i} !} \left(\stackrel{\rightarrow}{\mathrm{ad}}\, E_{{\bf s}_{j, i; p}} \right)^{-a_{ij} - 1} \left(V_{i, j; p - 1}^{(1)}\right).
\end{equation}
Since $E_{{\bf s}_{j, i; p}} \in \mathscr{F}_{ij}$, which is the induction hypothesis, we are reduced to verify $V_{i, j; p - 1}^{(1)} \in \mathscr{F}_{ij}$. When $a_{ij} \leq -2, a_{ji} \leq -2$, the fact $V_{i, j; 1}^{(1)} \in \mathscr{F}_{ij}$ derived from (\ref{eq:q_adjoint_power_commutation}) and inductive use of the formula (\ref{eq:V_ijp_formula1}) show $V_{i, j; p - 1}^{(1)} \in \mathscr{F}_{ij}$. When $a_{ij} = -1, a_{ji} \leq -4$, we can directly verify $V_{i, j; 0}^{(1)}$, $V_{i, j; 1}^{(1)}$, $V_{i, j; 0}^{(2)}$, $V_{i, j; 1}^{(2)} \in \mathscr{F}_{ij}$ using (\ref{eq:q_adjoint_power_commutation}), and the recurrence formulas (\ref{eq:V_ijp_formula2}), (\ref{eq:V_ijp_formula3}) show $V_{i, j; p - 1}^{(1)}, V_{i, j; p - 1}^{(2)} \in \mathscr{F}_{ij}$ for all $p \geq 1$.

When $a_{ij} \leq -4, a_{ji} = -1$, we need to continue the calculation of (\ref{eq:expansion}) slightly. Using the formula (\ref{eq:q_adjoint_power_commutation}), we have
\begin{equation*}
T_j \left(\left[E_i, [E_i, E_j]_q\right]_q\right) = \left[[E_j, E_i]_q, E_i\right]_q.
\end{equation*}
Thus $E_{{\bf s}_{i, j; p + 1}}$ can be written as
\begin{equation*}
E_{{\bf s}_{i, j; p + 1}} = \frac{1}{[-a_{ij}]_{q_i} !} \left(\stackrel{\rightarrow}{\mathrm{ad}}\, E_{{\bf s}_{j, i; p}} \right)^{-a_{ij} - 2}  \left(V_{j, i; p - 2}^{(2)}\right).
\end{equation*}
Since $V_{j, i; p - 2}^{(2)} \in \mathscr{F}_{ij}$ due to the discussion of the case $a_{ij} = -1, a_{ji} \leq -4$, we conclude $E_{{\bf s}_{i, j; p + 1}} \in \mathscr{F}_{ij}$.
\end{proof}
\end{lem}

By the proof of the theorem \ref{thm:q_comm}, we can easily construct an algorithm for describing $T_w(E_j)$ as a concrete $q$-commutator monomial once formulas for the elements of the form $\dots T_i T_j T_i (E_j)$ are prepared. In particular for simply laced case, we have a simple graphical algorithm for the calculation of $T_w(E_j)$, which we are going to describe below.

First, we introduce a graphical notation of $q$-bracket, which is convenient to write down $q$-commutator monomials.
\begin{equation} \label{eq:tree_notation}
\begin{xy}
 (-3, 3) *{X}, (3, 3) *{Y},
\ar @{-} (-3, 1); (0, -2), \ar @{-} (0, -2) ; (3, 1)
\end{xy} := [X, Y]_q \quad (X, Y \in U_q(\mathfrak{g})).
\end{equation}
We also abbreviate $E_i$ to $i$ in the schematic notation. For example, the $q$-Serre relation (\ref{eq:q_Serre_rel}) can be written as the following binary tree.
\begin{equation*}
\begin{xy}
(-15, 4) * {i}, (-12, 4) * {i}, (-9, 4) * {i}, (-6, 4) * {i}, (0, 4) * {\cdots}, (6, 4) * {i}, (9, 4) * {i}, (12, 4) * {j}, (-3, 9) * {1 - a_{ij}},
(-17, 5) ; (-3, 7) **\crv{(-16, 5.99) & (-9, 6) & (-4, 6.01)},
(11, 5) ; (-3, 7) **\crv{(10, 5.99) & (4, 6) & (0, 6.01)},
\ar @{-} (-15, 2) ; (0, -7),
\ar @{-} (-12, 2) ; (1.33, -6),
\ar @{-} (-9, 2) ; (2.66, -5),
\ar @{-} (-6, 2) ; (4, -4),
\ar @{-} (6, 2) ; (9.33, 0),
\ar @{-} (9, 2) ; (10.66, 1),
\ar @{-} (0, -7) ; (12, 2),
\end{xy}
= \left(\stackrel{\rightarrow}{\mathrm{ad}}\, E_i \right)^{1-a_{ij}}(E_j) = 0 \quad (i \neq j).
\end{equation*}
Using these notation, we can describe every $q$-commutator monomial as a binary tree, whose node represents $q$-bracket and whose leaf denotes a Chevalley generator $E_i$. Now we can describe the algorithm for simply laced case.
\begin{prop} \label{prop:algorithm}
Let $\mathfrak{g}$ be a simply laced Kac-Moody algebra, $w \in W$ and $j$ be an index satisfying $w(\alpha_j) \in \Delta_+$. Let $w = s_{i_1} s_{i_2} \dots s_{i_m}$ be a reduced expression. Then the binary tree constructed by the following procedure represents a $q$-commutator monomial equal to $T_w(E_j)$.

\begin{enumerate}
\item In this procedure, we manipulate a binary tree, each of whose leaf holds a pair of a reduced expression $s_{j_1} s_{j_2} \dots s_{j_k}$ and an index $p$ such that $s_{j_1} s_{j_2} \dots s_{j_k} (\alpha_p)$ $\in \Delta_+$.

\item At the beginning we have a binary tree consists of only the root, whose reduced expression is $s_{i_1} s_{i_2} \dots s_{i_m}$ and whose index is $j$. The procedure terminates immediately when $m = 0$.

\item For each leaf of the binary tree, the following manipulations are applied recursively. Let $s_{j_1} s_{j_2} \dots s_{j_k}$ and $p$ be the reduced expression and the index
of the leaf we are working on respectively.

\begin{enumerate}
\item We are done for the leaf if $k = 0$.

\item If $k \geq 1$, then $j_k \neq p$. If $a_{j_k p} = 0$, then delete the factor $s_{j_k}$ in the reduced expression since $T_{j_k} (E_p) = E_p$. Repeat this deletion until $a_{j_k p} = -1$.

\item If $s_{j_1} s_{j_2} \dots s_{j_{k - 1}} (\alpha_p) \in \Delta_-$, then there exists a number $l$ such that
\begin{equation*}
s_{j_1} s_{j_2} \dots s_{j_{k - 1}} = s_{j_1} s_{j_2} \dots s_{j_{l - 1}} s_{j_{l + 1}} \dots s_{j_{k - 1}} s_p,
\end{equation*}
due to the exchange condition ~\cite{Kac}. By $a_{j_k p} = a_{p j_k} = -1$ and (\ref{eq:reduction_A2}), we have
\begin{align*}
T_{j_1} T_{j_2} \dots T_{j_{k - 1}} T_{j_k} (E_p) &= T_{j_1} T_{j_2} \dots T_{j_{l - 1}} T_{j_{l + 1}} \dots T_{j_{k - 1}} T_p T_{j_k} (E_p) \\
&= T_{j_1} T_{j_2} \dots T_{j_{l - 1}} T_{j_{l + 1}} \dots T_{j_{k - 1}} (E_{j_k}).
\end{align*}
According to this calculation, replace the reduced expression with \\ $s_{j_1} s_{j_2} \dots s_{j_{l - 1}} s_{j_{l + 1}} \dots s_{j_{k - 1}}$ and replace the index with $j_k$. Repeat this replacement until $s_{j_1} s_{j_2} \dots s_{j_{k - 1}} (\alpha_p) \in \Delta_+$.

\item Finally, when $s_{j_1} s_{j_2} \dots s_{j_{k - 1}} (\alpha_p) \in \Delta_+$, then $a_{j_k p} = a_{p j_k} = -1$ by the manipulations so far. Thus $T_{j_k} (E_p) = [E_{j_k}, E_p]_q$ by (\ref{eq:braid_action2}), and \\ $s_{j_1} s_{j_2} \dots s_{j_{k - 1}} s_{j_k}$, $s_{j_1} s_{j_2} \dots s_{j_{k - 1}} s_p$ are reduced. Therefore, create new branch at the current leaf and generate two leaves as below, where ${\bf s}' := s_{j_1} s_{j_2} \dots s_{j_{k - 1}}$. The new two leaves have index $j_k, p$ respectively, and both reduced expression is ${\bf s}'$. Figure \ref{fig:branching} shows this branching procedure, where ${\bf s}[p]$ denotes reduced expression ${\bf s}$ and index $p$.
\begin{figure}[htbp]
\centering
\begin{equation*}
\begin{xy}
(0, 5) * {{\bf s}'s_{j_k}[p]}, (10, 0) * {\cdots},
\ar @{-} (0, 2) ; (7, 0)
\end{xy}
\quad \rightsquigarrow \quad
\begin{xy}
(0, 8) * {{\bf s}'[j_k]}, (10, 8) * {{\bf s}'[p]}, (15, 0) * {\cdots},
\ar @{-} (0, 5) ; (5, 2),
\ar @{-} (5, 2) ; (10, 5),
\ar @{-} (5, 2) ; (12, 0)
\end{xy}
\end{equation*}
\caption{Branching rule}
\label{fig:branching}
\end{figure}
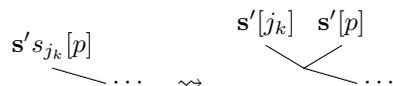
\end{enumerate}

\item Repeat above procedure until all reduced expressions in the leaves has length 0. This algorithm terminates within finite steps because each manipulation shortens the length of reduced expression of target leaf.
\end{enumerate}
\end{prop}

\section{Construction of quantum dilogarithm identities} \label{sec:4}

In this section, we show how to construct quantum dilogarithm identities using the product formula (\ref{eq:prod_formula}) of the quasi-universal R-matrix $\Theta$.

First, we introduce certain projections of the algebra $U_q^+ \widehat{\otimes} U_q^-$ onto skew formal power series algebras determined by Dynkin quivers. Through the projections, most of the elements of the form $\Theta_{\leq, \alpha} \in U_q^+ \widehat{\otimes} U_q^-$ $(\alpha \in \Delta_+^{\mathrm{re}})$ become the unit of image, while several factors survive and retain their form as $q$-exponential function, which can also be seen as quantum dilogarithm functions. Moreover, under appropriate setting of parameters, the product of factors in $\Theta$ associated with imaginary roots can be written using quantum dilogarithm functions in the image of the projection. Thus the image of $\Theta$ will be written as certain product of quantum dilogarithm functions. Choosing various convex orders, one can obtain various product presentations of the image of $\Theta$, which have finitely or infinitely many factors depending on selected order. Eventually, we can construct quantum dilogarithm identities of the form ``finite product = infinite product'', which will exactly coincide with the identities proposed in ~\cite{DGS} after suitable change of variables.

\subsection{Projections of $U_q^+ \widehat{\otimes} U_q^-$ onto skew formal power series algebras}
\label{sec:representation}
Let $\mathfrak{g}$ be symmetrisable Kac-Moody algebra of rank $n$ and $A = (a_{ij})_{i, j = 1}^n$ be its Cartan matrix. Let $d_1, d_2, \dots, d_n$ be coprime positive integers such that $c_{ij} := d_i a_{ij} = d_j a_{ji}$ $(1 \leq i, j \leq n)$. Then $C := (c_{ij})_{i, j = 1}^n$ is symmetrized matrix of $A$. We normalize the invariant bilinear form $(\cdot, \cdot)$ so that $(\alpha_i, \alpha_j) = c_{ij}$ for $1 \leq i, j \leq n$. Choose $\sigma_{ij} \in \{\pm1\}$ for each pair of indices $i < j$ such that $a_{ij} \neq 0$, and set
\begin{equation}
b_{ij} := \begin{cases}
\sigma_{ij} c_{ij} & i < j \\
0 & i = j \\
-\sigma_{ij} c_{ij} & i > j \\
\end{cases}.
\end{equation}
Then the matrix $B = (b_{ij})_{i, j = 1}^n$ is skew-symmetric matrix, and this data can be interpreted as the Dynkin quiver which has an arrow from $i$ to $j$ if $\sigma_{ij} = +1$. Let $\{\cdot, \cdot\}_B: Q \times Q \rightarrow \Z$ be the skew-symmetric form satisfying $\{\alpha_i, \alpha_j\}_B = b_{ij}$ for all indices $i, j$.

Let $\mathcal{P}_B$ be a $\Q(q)$-algebra defined by the generators and relations below.
\begin{align*}
\text{Generators} &: e_1, e_2, \dots, e_n. \\
\text{Relations} &: e_i e_j = q^{b_{ij}} e_j e_i \quad (i, j = 1, 2, \dots, n).
\end{align*}
$\mathcal{P}_B$ has natural $Q$-graded algebra structure if each $e_i$ is supposed to have weight $\alpha_i$, hence $q$-bracket makes sense on $\mathcal{P}_B$. Each weight space of $\mathcal{P}_B$ is one dimensional subspace spanned by a monomial of the form $e_1^{k_1} e_2^{k_2} \dots e_n^{k_n}$. The degree of each monomial
\begin{equation}
\deg e_1^{k_1} e_2^{k_2} \dots e_n^{k_n} := k_1 + k_2 + \dots + k_n
\end{equation}
coincides with the height of its weight. Let $\mathcal{P}_{B, m}$ be the subspace of $\mathcal{P}_B$ spanned by monomials of degree $m$. Also note that
\begin{equation}
[e_i, e_j]_q = e_i e_j - q^{(\alpha_i, \alpha_j)} e_j e_i = \left(q^{b_{ij}} - q^{c_{ij}}\right) e_j e_i = 0 \quad \text{if} \  \sigma_{ij} = +1.
\end{equation}

Due to the following well-known fact, $\mathcal{P}_B$ turns out to be a quotient of $U_q^+ \subset U_q(\mathfrak{g})$.
\begin{prop}~\cite{Lus}[33.1.3] The positive part of the quantum enveloping algebra $U_q^+$ is isomorphic to the $\Q(q)$-algebra whose generators are $E_1, E_2, \dots,$ $E_n$ and whose relation is given by the quantum Serre relation (\ref{eq:q_Serre_rel}).
\end{prop}

\begin{prop}\label{prop:projection}
There exists unique $Q$-graded algebra surjection
\begin{equation}
\pi_B: U_q^+ \rightarrow \mathcal{P}_B
\end{equation}
such that $\pi_B(E_i) = e_i$ for all $i = 1, \dots, n$. 
\end{prop}

Let $\mathcal{P}_B^+ := \mathcal{P}_B$ and $\mathcal{P}_B^-$ be a copy of $\mathcal{P}_B$ but whose generators $e_i$ are replaced by $f_i$. Recall that $U_q^-$ is isomorphic to $U_q^+$ as algebra~\cite{Lus}[3.2.6]. Let $\pi_B^+ := \pi_B: U_q^+ \rightarrow \mathcal{P}_B^+$, $\pi_B^-: U_q^- \rightarrow \mathcal{P}_B^-$ $(\pi_B^-(F_i) := f_i)$ are $Q$-graded algebra surjections given by the Proposition \ref{prop:projection}. Then we have an algebra surjection
\begin{equation}
\pi_B^+ \otimes \pi_B^-: U_q^+ \otimes U_q^- \twoheadrightarrow \mathcal{P}_B^+ \otimes \mathcal{P}_B^-.
\end{equation}

We want to construct a completion of this surjection to define image of $\Theta \in U_q^+ \widehat{\otimes} U_q^-$. To define the completion, we give a topology on $\mathcal{D}_B := \mathcal{P}_B^+ \otimes \mathcal{P}_B^-$ so that the surjection $\pi_B^+ \otimes \pi_B^-$ becomes continuous.

For every nonnegative integer $m$, set
\begin{equation}
\mathcal{D}_m := \bigoplus_{k = m}^{\infty} \mathcal{P}_{B, m}^+ \otimes \mathcal{P}_B^- \quad \subset \mathcal{D}_B,
\end{equation}
and define the completion of $\mathcal{D}_B$ by
\begin{equation}
\widehat{\mathcal{D}_B}  := \underset{m \geq 0}{\mathrm{\projlim}} \ \mathcal{D}_B / \mathcal{D}_m.
\end{equation}
Recalling the definition (\ref{eq:topology_on_U_q}) of $\widehat{U_q} \widehat{\otimes} \widehat{U_q}$, composition of the surjection $\pi_B^+ \otimes \pi_B^-: U_q^+ \otimes U_q^- \twoheadrightarrow \mathcal{D}_B$ and inclusion $\iota: \mathcal{D}_B \hookrightarrow \widehat{\mathcal{D}_B}$ is continuous with respect to the relative topology on $U_q^+ \otimes U_q^- \subset \widehat{U_q} \widehat{\otimes} \widehat{U_q}$. Hence this map induces an unique continuous map
\begin{equation}
\pi_B^+ \widehat{\otimes} \pi_B^-: U_q^+ \widehat{\otimes} U_q^- \rightarrow \widehat{\mathcal{D}_B}
\end{equation}
due to the completeness of $\widehat{\mathcal{D}_B}$.

\subsection{Skew formal power series algebras}

Let $y_i := (q_i - q_i^{-1}) e_i \otimes (q_i^{-1} - q_i) f_i \in \widehat{\mathcal{D}_B}$, $\mathcal{S}_B \subset \widehat{\mathcal{D}_B}$ be $\Q(q)$-subalgebra generated by $y_1, y_2, \dots, y_n$, and $\widehat{\mathcal{S}_B} \subset \widehat{\mathcal{D}_B}$ be its closure. Since the increasing monomials 
\begin{equation*}
e_1^{m_1} e_2^{m_2} \dots e_n^{m_n} \otimes f_1^{m'_1} f_2^{m'_2} \dots f_n^{m'_n} \in \widehat{\mathcal{D}_B}
\end{equation*}
form a topological basis of $\widehat{\mathcal{D}_B}$, the increasing monomials $y_1^{m_1} y_2^{m_2} \dots y_n^{m_n}$ form a topological basis of $\widehat{\mathcal{S}_B}$. Therefore $\widehat{\mathcal{S}_B}$ is isomorphic to formal power series algebra $\Q(q) [[y_1, y_2, \dots, y_n]]$ as $\Q(q)$-linear space. This isomorphism endows the $\Q(q)$-linear space $\Q(q) [[y_1, y_2, \dots, y_n]]$ with a complete topological $\Q(q)$-algebra structure, whose multiplication is uniquely determined by the commutation relations $y_i y_j = q^{2b_{ij}} y_j y_i$.

In the same way, skew Laurent polynomial algebra $\mathcal{L}_B$ can be defined. Namely, $\mathcal{L}_B$ is Laurent polynomial algebra $\Q(q)[y_1^{\pm 1}, y_2^{\pm 2}, \dots, y_n^{\pm 1}]$ as $\Q(q)$-linear space, and multiplication in $\mathcal{L}_B$ is uniquely defined by the commutation relations $y_i y_j = q^{2b_{ij}} y_j y_i$. $\mathcal{S}_B$ can be naturally considered as a subalgebra of $\mathcal{L}_B$.

Let $L$ be the lower triangular part of $B$. Since $B$ is skew-symmetric, $B = L - {}^t\!L$. We define normal ordered product in $\mathcal{L}_B$ by
\begin{equation}
\normord{y^{\bf m}} := q^{{}^t\!{\bf m} L {\bf m}} y^{\bf m} \quad ({\bf m} = {}^t\!(m_1, m_2, \dots, m_n) \in \Z^n),
\end{equation}
where $y^{\bf m} := y_1^{m_1} y_2^{m_2} \dots y_n^{m_n}$.

Let $B' = (b'_{kl})_{k, l = 1}^{n'} \in M_{n'} (\Z)$ be another skew-symmetric matrix. We shall consider algebra homomorphism $\psi_R: \mathcal{L}_B \rightarrow \mathcal{L}_{B'} \cong \Q(q) [y_1'^{\pm 1}, y_2'^{\pm 1}, \dots, y_{n'}'^{\pm 1}]$ which is determined by a $n' \times n$-matrix $R \in M_{n',n}(\Z)$ and
\begin{equation}
\psi_R (y_i) := \normord{y'^{R v_i}},
\end{equation}
where $v_i \in \Z^n$ is the $i$-th unit vector. $\psi_R$ is well-defined if and only if it preserves the commutation relation $y_i y_j = q^{2b_{ij}} y_j y_i$ for all $i, j = 1, 2, \dots, n$, in other words
\begin{equation*}
\normord{y'^{R v_i}} \normord{y'^{R v_j}} = q^{2 {}^t\!v_i B v_j } \normord{y'^{R v_j}} \normord{y'^{R v_i}} \quad (i, j = 1, 2, \dots, n).
\end{equation*}
On the other hand, 
\begin{equation*}
y'^{R v_i} y'^{R v_j} = q^{2{}^t\!(R v_i) B' R v_j} y'^{R v_j} y'^{R v_i}\quad \in \mathcal{L}_{B'}.
\end{equation*}
Thus $\psi_R$ is well-defined if and only if ${}^t\!v_i {}^t\!R B' R v_j =  {}^t\!v_i B v_j$ for all $i, j$. This shows

\begin{prop} \label{prop:S_B_hom}
Let $B \in M_n (\Z)$, $B' \in M_{n'} (\Z)$ be skew-symmetric matrices, and $R$ be integer-valued $n' \times n$-matrix. There exists unique algebra homomorphism $\psi_R: \mathcal{L}_B \rightarrow \mathcal{L}_{B'}$ satisfying $\psi_R (y_i) = \normord{y'^{R v_i}}$ $(i = 1, 2, \dots, n)$ if and only if
\begin{equation} \label{eq:cond_S_B_hom}
{}^t\!R B' R = B.
\end{equation}
\end{prop}

Moreover, $\psi_R$ preserves the normal ordered product.
\begin{prop} \label{prop:preserve_normord}
Suppose that $R$ satisfies (\ref{eq:cond_S_B_hom}). Then
\begin{equation} \label{eq:preserve_normord}
\psi_R (\normord{y^{\bf m}}) = \normord{y'^{R {\bf m}}} \quad ({\bf m} \in \Z^n).
\end{equation}
\end{prop}

\begin{proof}
We prove by induction on $\mathrm{deg}\,{\bf m} := |m_1| + |m_2| + \dots + |m_n|$. The case $\mathrm{deg}\,{\bf m} = 1$ is trivial since ${}^t\!v_i L v_i = 0$ for any unit vector $v_i$. Suppose that (\ref{eq:preserve_normord}) holds if $\mathrm{deg}\,{\bf m} < N$ for some integer $N \geq 2$. Let ${\bf m}$ be of degree $N$. Choose ${\bf m}_1, {\bf m}_2 \in \Z^n$ so that ${\bf m} = {\bf m}_1 + {\bf m}_2$ and $\mathrm{deg}\,{\bf m}_k < N$ $(k = 1, 2)$. Since $y^{{\bf m}_1} y^{{\bf m}_2}= q^{2{}^t\!{\bf m}_1 L {\bf m}_2} y^{{\bf m}_1 + {\bf m}_2}$, we have
\begin{align*}
\normord{y^{{\bf m}_1 + {\bf m}_2}} &= q^{{}^t\!({\bf m}_1 + {\bf m}_2) L ({\bf m}_1 + {\bf m}_2)} y^{{\bf m}_1 + {\bf m}_2} \\
&= q^{{}^t\!{\bf m}_2 L {\bf m}_1 - {}^t\!{\bf m}_1 L {\bf m}_2} \normord{y^{{\bf m}_1}} \normord{y^{{\bf m}_2}} \\
&= q^{{}^t\!{\bf m}_2 B {\bf m}_1} \normord{y^{{\bf m}_1}} \normord{y^{{\bf m}_2}}.
\end{align*}
Then by the induction hypothesis,
\begin{align*}
\psi_R(\normord{y^{{\bf m}_1 + {\bf m}_2}}) &= q^{{}^t\!{\bf m}_2 B {\bf m}_1} \psi_R (\normord{y^{{\bf m}_1}} \normord{y^{{\bf m}_2}}) \\
&= q^{{}^t\!{\bf m}_2 B {\bf m}_1} \normord{y'^{R {\bf m}_1}} \normord{y'^{R {\bf m}_2}} \\
&= q^{{}^t\!{\bf m}_2 B {\bf m}_1 - {}^t\!(R{\bf m}_2) B' R {\bf m}_1} \normord{y'^{{\bf m}_1 + {\bf m}_2}} \\
&= \normord{y'^{{\bf m}_1 + {\bf m}_2}}. \quad (\because \text{(\ref{eq:cond_S_B_hom})})
\end{align*}
Thus (\ref{eq:preserve_normord}) holds for arbitrary ${\bf m}$ of degree $N$.
\end{proof}

When all the component of matrix $R$ satisfying (\ref{eq:cond_S_B_hom}) are nonnegative, we have the restricted homomorphism $\psi_R: \mathcal{S}_B \rightarrow \mathcal{S}_{B'}$.
\begin{prop} \label{prop:R_conti}
Suppose that $R \in M_{n',n}(\Z_{\geq 0})$ satisfies (\ref{eq:cond_S_B_hom}). The algebra homomorphism $\psi_R: \mathcal{S}_B \rightarrow \mathcal{S}_{B'}$ is continuous with respect to the relative topology in $\widehat{\mathcal{S}_B}$ and $\widehat{\mathcal{S}_{B'}}$ if and only if each column of $R$ contains nonzero component.
\end{prop}
\begin{proof}
By the definition of $\psi_R$, it is continuous if and only if:
\begin{quote}
For any $N \in \Z_{\geq 0}$, there exists some $M \in \Z_{\geq 0}$ such that ${\bf m} = {}^t\!(m_1, \dots, m_n) \in \Z_{\geq 0}^n$ of total degree $\deg {\bf m} \geq M \Rightarrow \deg R{\bf m} \geq N$.
\end{quote}
Since $\deg R{\bf m} = \sum_{i = 1}^n r_i m_i$ where $r_i$ is the sum of components in $i$-th column of $R$, this condition holds if and only if all the $r_i$ are positive.
\end{proof}
Thus, when $R$ satisfies these conditions, $\psi_R$ uniquely extends to continuous algebra homomorphism $\widehat{\psi_R}: \widehat{\mathcal{S}_B} \rightarrow \widehat{\mathcal{S}_{B'}}$.

\subsection{Several formulas related to quantum dilogarithm}

To compute the image of the quasi-universal R-matrix $\Theta$, we briefly prepare for a couple of formulas related to quantum dilogarithm function $\mathrm{Li}_{2, q} (x)$. Let
\begin{align}
\log (1 - x) &:= -\sum_{n = 1}^{\infty} \frac{x^n}{n}, \label{eq:def_log} \\
\mathrm{Li}_{2, q} (x) &:= \sum_{n = 1}^{\infty} \frac{x^n}{n(1 - q^{n})}, \\
\mathbb{E} (x) &:= \exp \left(\mathrm{Li}_{2, q^2} \left(-qx\right)\right), \\
(x; q)_{\infty} &:= \prod_{n = 0}^{\infty} (1 - q^n x).
\end{align}
In this paper we consider these functions just as formal power series. The function $(x; q)_{\infty}$ is characterized by the recurrence relation
\begin{equation}
(1 - x) (qx; q)_{\infty} = (x; q)_{\infty}.
\end{equation}

Since $-\mathrm{Li}_{2, q} (qx) + \log (1 - x) = -\mathrm{Li}_{2, q} (x)$, $\exp (-\mathrm{Li}_{2, q} (x))$ satisfies the recurrence relation. Therefore
\begin{align}
\exp \mathrm{Li}_{2, q} (x) &= (x; q)_{\infty}^{-1}, \\
\mathbb{E} (x) &= (-qx; q^2)_{\infty}^{-1},
\end{align}
which coincides with the product presentation of $\mathbb{E} (x)$ in the introduction. By this presentation, $\mathbb{E} (x)$ is characterized by the recurrence relation
\begin{equation}
(1 + qx) \mathbb{E} (x) = \mathbb{E} (q^2x).
\end{equation}

Recall that $q$-exponential function was defined by
\begin{equation}
\exp_q(x) := \sum_{n = 0}^{\infty} \frac{q^{-\frac{1}{2} n(n-1)}}{[n]_q !} x^n.
\end{equation}
Then it can be directly verified that $\exp_q(x)$ satisfies
\begin{equation}
(1 + q(q - q^{-1}) x) \exp_q(x) = \exp_q(q^2x)
\end{equation}
and we conclude that
\begin{equation}
\exp_q(x) = \mathbb{E} \left(\left(q-q^{-1}\right) x\right). \label{eq:qexp_bb}
\end{equation}

In the same way, we can also prove another presentation of $(x; q)_{\infty}$ ~\cite{Tera}, which we will use in the computation of imaginary root vectors.
\begin{equation}
\exp \left(\sum_{m = 1}^{\infty} -\frac{1}{m(1 - q^{m})} x^m\right) = (x; q)_{\infty}. \label{eq:q_dilog}
\end{equation}

\subsection{Computation of the image of quasi-universal R-matrix $\Theta$}

Now, we suppose that $\mathfrak{g}$ is untwisted affine Lie algebra. We shall compute \\ $\pi_B^+ \widehat{\otimes} \pi_B^- (\Theta) \in \widehat{\mathcal{D}_B}$ for various product presentation of $\Theta$ (\ref{eq:prod_formula}) and equate them to obtain concrete identities. 

First we remark that $q$-commutator monomials degenerate to ordinary monomials. Let $X_{\alpha}, X_{\beta} \in \mathcal{P}_B^+$ have weight $\alpha, \beta \in Q$ respectively. Then by definition of $\mathcal{P}_B^+$, $X_{\alpha}$ is a linear combination of monomials $e_{i_1} e_{i_2} \dots e_{i_m}$, where $\alpha_{i_1} + \alpha_{i_2} + \dots + \alpha_{i_m} = \alpha$. Hence $X_{\alpha} X_{\beta} = q^{\{\alpha, \beta\}_B} X_{\beta} X_{\alpha}$ and we have
\begin{equation} \label{eq:weight_vec_comm}
\left[X_{\alpha}, X_{\beta}\right]_q := X_{\alpha} X_{\beta} - q^{(\alpha, \beta)} X_{\beta} X_{\alpha} = \left(1 - q^{(\beta, \alpha) + \{\beta, \alpha\}_B}\right) X_{\alpha} X_{\beta}.
\end{equation}
It is convenient to introduce the bilinear form $\langle \alpha, \beta \rangle_B := (\alpha, \beta) - \{\alpha, \beta\}_B$. Then the values of the bilinear form $\langle \cdot, \cdot \rangle_B$ are even integers, since
\begin{equation*}
\langle \alpha_i, \alpha_j \rangle_B = c_{ij} - b_{ij} = \begin{cases}
(1 \pm 1) c_{ij} & i \neq j \\
2d_i & i = j
\end{cases}.
\end{equation*}
The formula (\ref{eq:weight_vec_comm}) shows that $\left[X_{\alpha}, X_{\beta}\right]_q$ vanishes if and only if $\langle \alpha, \beta \rangle_B = 0$, otherwise it is nonzero multiple of $X_{\alpha} X_{\beta}$. Therefore we have the following vanishing criteria for $q$-commutator monomials.

\begin{prop} A $q$-commutator monomial $M \in U_q^+$ lies in the kernel of $\pi_B^+$ if and only if there exists an application of $q$-bracket $\left[E_{\alpha}, E_{\beta}\right]_q$ in $M$ for some $E_{\alpha}, E_{\beta} \in U_q^+$ of weight $\alpha, \beta \in Q$ satisfying $\langle \alpha, \beta \rangle_B = 0$. If there are no such application of $q$-bracket in $M$, $\pi_B^+(M) \in \mathcal{P}_B^+$ is a nonzero monomial.
\end{prop}

Recall that the root vectors for $U_q^-$ were defined by $F_{\leq, \alpha} := \Omega(E_{\leq, \alpha})$, where $\alpha \in \Delta_+^{\mathrm{re}}$ and $\Omega: U_q^+ \rightarrow U_q^-$; $\Omega(E_i) := F_i$, $\Omega(q) := q^{-1}$ was Chevalley involution, which is anti-automorphism of $\Q$-algebra. Notice that $\Omega$ preserves $q$-bracket except for multiple of a power of $q$.
\begin{equation}
\Omega\left(\left[E_{\alpha}, E_{\beta}\right]_q\right) = -q^{-(\alpha, \beta)} \left[F_{\alpha}, F_{\beta}\right]_q.
\end{equation}
Thus $F_{\leq, \alpha}$ are also $q$-commutator monomials, and they coincide with $E_{\leq, \alpha}$ except for multiple of $\pm q^k$ and replacing $E_i$ with $F_i$.

There exists unique anti-isomorphism of $\Q$-algebra $\overline{\Omega_B}: \mathcal{P}_B^+ \rightarrow \mathcal{P}_B^-$ which sends $e_i$ to $f_i$ and $q$ to $q^{-1}$. This is useful to compute $\pi_B^- (F_{\leq, \alpha})$ because
\begin{equation} \label{eq:omega_bar_commute}
\overline{\Omega_B} \circ \pi_B^+ = \pi_B^- \circ \Omega.
\end{equation}

Since $\mathrm{ht} \, \alpha - 1$ times of $q$-brackets occur in $E_{\leq, \alpha}$ for $\alpha \in \Delta_+^{\mathrm{re}}$, its image takes the form below.
\begin{equation}
\pi_B^+ (E_{\leq, \alpha}) = C q^u \left(\prod_{i = 1}^{\mathrm{ht} \, \alpha - 1} \left(1 - q^{k_i}\right) \right) e_0^{m_0} e_1^{m_1} \dots e_n^{m_n},
\end{equation}
where $C \in \Q(q)$ is the coefficient of $E_{\leq, \alpha}$ as $q$-commutator monomial, $\alpha = \sum_{i = 0}^n m_i \alpha_i$ and $u, k_i \in \Z$. When simply-laced case, $C = 1$ because no non-trivial scalar multiple occur in the algorithm of section \ref{sec:2.5} (Proposition \ref{prop:algorithm}). We also note that each $k_i$ is a value of the bilinear form $\langle \cdot, \cdot \rangle_B$ and thus even integer. Using $\overline{\Omega_B}$, the image of $F_{\leq, \alpha}$ is
\begin{equation}
\pi_B^- (F_{\leq, \alpha}) = \overline{\Omega_B} \circ \pi_B^+(E_{\leq, \alpha}) = \Omega(C) q^{-u} \left(\prod_{i = 1}^{\mathrm{ht} \, \alpha - 1} \left(1 - q^{-k_i}\right) \right) f_n^{m_n} f_{n - 1}^{m_{n - 1}} \dots f_0^{m_0}.
\end{equation}

Recall that the subalgebra $\mathcal{S}_B \subset \mathcal{D}_B$, which is generated by $y_i := (q_i - q_i^{-1}) (q_i^{-1} - q_i) e_i \otimes f_i$ for $i = 0, 1, \dots, n$. The normal ordered product of monomial $y^{\bf m} := y_0^{m_0} y_1^{m_1} \dots y_n^{m_n}$ $({\bf m} = {}^t\!(m_0, m_1, \dots, m_n) \in \Z_{\geq 0}^{n + 1})$ was defined as $\normord{y^{\bf m}} := q^{{}^t\!{\bf m} L {\bf m}} y^{\bf m}$, where $L$ was the lower triangular part of $B$.

Using these notations,  when simply-laced case we have
\begin{multline}
\pi_B^+ \otimes \pi_B^- \left( \left(q - q^{-1}\right) \left(q^{-1} - q\right) E_{\leq, \alpha} \otimes F_{\leq, \alpha}\right) \\
= \left(\prod_{i = 1}^{\mathrm{ht} \, \alpha - 1} \left[\frac{k_i}{2}\right]_q^2 \right) \normord{y_0^{m_0} y_1^{m_1} \dots y_n^{m_n}}.
\end{multline}

Recall that $\exp_q (x) = \mathbb{E} ((q - q^{-1})x)$ (\ref{eq:qexp_bb}). Therefore if $k_i = \pm 2$ for all $i$, we have simple description of the image of $\Theta_{\leq, \alpha}$ for positive real root $\alpha$.
\begin{prop} \label{prop:image_real_root_vec}
Let $\Delta$ be simply-laced affine root system, $\alpha = \sum_{i = 0}^n m_i \alpha_i \in \Delta_+^{\mathrm{re}}$ and $\leq$ be a convex order. We suppose that when a presentation of $E_{\leq, \alpha}$ as $q$-commutator monomial is given, $\langle \alpha, \beta\rangle_B = \pm 2$ for each application of $q$-bracket $[X_{\alpha}, X_{\beta}]_q$, where $\alpha, \beta \in \Delta^+$ and $X_{\alpha} \in U_{\alpha}^+, X_{\beta} \in U_{\beta}^+$. Then
\begin{equation}
\pi_B^+ \widehat{\otimes} \pi_B^- (\Theta_{\leq, \alpha}) = \mathbb{E}\left(\normord{y_0^{m_0} y_1^{m_1} \dots y_n^{m_n}}\right).
\end{equation}
\end{prop}

While we have the general simple description of the images of real root vectors, the computation of the images of imaginary root vectors requires some ingenuity. Recall that $\varphi_{i, n} \in U_{n\delta}^+$ were defined as
\begin{equation}
\varphi_{i, n} := \left[T_{\varepsilon_i}^n T_i^{-1} \left(E_i\right), E_i\right]_q \quad (n \in \Z_{\geq 1}, i \in \mathring{I}),
\end{equation}
and imaginary root vectors $I_{i, n}$ were polynomials consists of $\varphi_{i, n}$. Since $T_{\varepsilon_i}^n T_i^{-1} \left(E_i\right)$ can be written as a $q$-commutator monomial by using the algorithm in section \ref{sec:2.5}, $\varphi_{i, n}$ themselves are $q$-commutator monomials. But we need to compute $T_w (I_{i, n})$ $(w \in \mathring{W})$ for general convex order, and we cannot apply the algorithm to $T_w (\varphi_{i, n})$ when $w(\alpha_i) \in \Delta_-$ because $T_w (E_i)$ no longer lies in $U_q^+$.

First, we compute $T_i (\varphi_{i, n})$ using the following fact.
\begin{prop}~\cite{Ito0}
For every $i, j \in \mathring{I}$ and positive integer $n$,
\begin{equation}
T_{\varepsilon_j} (\varphi_{i, n}) = \varphi_{i, n}.
\end{equation}
\end{prop}
We will also use the property that for $u, v \in \widehat{W}$, $T_{uv} = T_{u} T_{v}$ if $\ell(uv) = \ell(u) + \ell(v)$.

Since $t_{\varepsilon_i} (\alpha_i) = -\delta + \alpha_i \in \Delta_-$, the length of $u := t_{\varepsilon_i} s_i$ in $\widehat{W}$ is $\ell(u) = \ell(t_{\varepsilon_i}) - 1$. Thus $T_{\varepsilon_i} = T_u T_i$ and we have
\begin{equation} \label{eq:invert}
T_i (\varphi_{i, n}) = T_u^{-1} (\varphi_{i, n}) = \left[T_i T_{\varepsilon_i}^{n - 1} T_i^{-1} \left(E_i\right), T_u^{-1} (E_i)\right]_q.
\end{equation}

Now we can compute arbitrary $T_w (\varphi_{i, n})$ for $w \in \mathring{W}$. When $w(\alpha_i) \in \Delta_+$, we have
\begin{equation}
T_w (\varphi_{i, n}) = \left[T_w T_{\varepsilon_i}^n T_i^{-1} \left(E_i\right), T_w \left(E_i\right)\right]_q
\end{equation}
and thus simply applying the algorithm to $T_w T_{\varepsilon_i}^{n - 1} T_u (E_i)$ and $T_w (E_i)$ yields explicit presentation of $T_w (\varphi_{i, n})$ as a $q$-commutator monomial.

When $w(\alpha_i) \in \Delta_-$, let $w' := w s_i$. Then $T_w = T_{w'} T_i$ and using (\ref{eq:invert}) we have
\begin{equation}
T_w (\varphi_{i, n}) = T_{w'} \left[T_i T_{\varepsilon_i}^{n - 1} T_i^{-1} \left(E_i\right), T_u^{-1} \left(E_i\right)\right]_q.
\end{equation}
Let $\lambda := \varepsilon_1 + \varepsilon_2 + \dots + \varepsilon_{\ell} \in \mathring{\mathfrak{h}}^*$, which is strictly dominant weight. Then $T_{\lambda} := T_{t_{\lambda}} = T_{\varepsilon_1} T_{\varepsilon_2} \dots T_{\varepsilon_{\ell}}$ and $t_{\lambda} (\alpha_i) = -\delta + \alpha_i \in \Delta_-$ $(i \in \mathring{I})$. Thus $t_{\lambda}$ inverts every positive roots in $\mathring{\Delta}_+$, which implies that there is an expression $t_{\lambda} = v w_{\circ}$ where $v \in \widehat{W}$ and $w_{\circ} \in \mathring{W}$ is the longest element satisfying $\ell(t_{\lambda}) = \ell(v) + \ell(w_{\circ})$. Notice that
\begin{equation*}
T_{w'} T_u^{-1} (E_i) = T_w T_i^{-1} T_u^{-1} (E_i) = T_w T_{\varepsilon_i}^{-1} (E_i),
\end{equation*}
and $T_{\varepsilon_i}^{-1} (E_i) = T_{\lambda}^{-1} (E_i)$ since $T_{\varepsilon_j} (E_i) = E_i$ if $i \neq j$. We also note that $T_{w_{\circ}} = T_{w_{\circ} w^{-1}} T_w$ for every $w \in \mathring{W}$ due to the maximality of $w_{\circ}$.

Recall the anti-automorphism of $\Q(q)$-algebra $\Psi: U_q^+ \rightarrow U_q^+$ defined by $\Psi(E_i) := E_i$. Since $\Psi$ preserves weights, it reverses $q$-bracket: 
\begin{equation}
\Psi([x, y]_q) = [\Psi(y), \Psi(x)]_q \quad (x, y \in U_q^+).
\end{equation}
It is also easy to verify that $\Psi T_i = T_i^{-1} \Psi$ for $i \in I$.

Let $t_{\lambda} = \tau v' w_{\circ}$ $(\tau \in \Omega, v' \in W)$ and $w_{\circ} = \Tilde{w} w$. Then
\begin{align*}
T_w T_{\lambda}^{-1} (E_i) &= T_w T_w^{-1} T_{\Tilde{w}}^{-1} T_{v'}^{-1} T_{\tau}^{-1} (E_i) \\
&= T_{\Tilde{w}}^{-1} T_{v'}^{-1} T_{\tau}^{-1} \Psi (E_i) \\
&= \Psi T_{(v'\Tilde{w})^{-1}} (E_{\tau^{-1} (i)}).
\end{align*}

Finally we have
\begin{equation}
T_w (\varphi_{i, n}) = \left[T_w T_{\varepsilon_i}^{n - 1} T_i^{-1} \left(E_i\right), \Psi T_{(v'\Tilde{w})^{-1}} \left(E_{\tau^{-1} (i)}\right)\right]_q.
\end{equation}
We can apply the algorithm to $T_{\Tilde{w} v'} \left(E_{\tau^{-1} (i)}\right)$. Since $\Psi$ just reverses the directions of $q$-brackets, $\Psi T_{\Tilde{w} v'} \left(E_{\tau^{-1} (i)}\right)$ is a $q$-commutator monomial. 

\begin{prop}
$T_w (\varphi_{i, n})$ is a $q$-commutator monomial for every $w \in \mathring{W}$, $i \in \mathring{I}$, and positive integer $n$.
\end{prop}

\section{Examples of quantum dilogarithm identities} \label{sec:5}

In this final section, we give specific convex orders and Dynkin quivers, which eventually induce the identities proposed in ~\cite{DGS}.

Recall that affine positive root system $\Delta_+$ is decomposed as
\begin{equation*}
\Delta_+ = \Delta(w, -) \amalg \Delta_+^{\mathrm{im}} \amalg \Delta(w, +)
\end{equation*}
and convex orders on $\Delta_+$ consists of convex orders on each $\Delta(w, \pm)$ (the order on $\Delta_+^{\mathrm{im}}$ is not significant since any total order can be chosen). Convex order on $\Delta(w, -)$ was determined by the following parameters with several restrictions (\ref{eq:chain_of_biconvex_set}) (\ref{eq:accumulate_infinite_reduced_words}):
\begin{enumerate}
\item A positive integer $n$ and a filtration of indices

$\mathring{I} = J_0 \supsetneqq J_1 \supsetneqq J_2 \supsetneqq \dots \supsetneqq J_n = \emptyset$.
\item $y_1 \in W_{J_1}, \, y_2 \in W_{J_2}, \, \dots, \, y_n \in W_{J_n}$.
\item Infinite reduced words ${\bf s}_0 \in \mathscr{W}_{J_0}^{\infty}, \, {\bf s}_1 \in \mathscr{W}_{J_1}^{\infty}, \, \dots, \, {\bf s}_{n - 1} \in \mathscr{W}_{J_{n - 1}}^{\infty}$.
\end{enumerate}

We have to specify not only the parameters for $\Delta(w, -)$, but also for  $\Delta(w, +) = \Delta(w w_{\circ}, -)$ to construct whole convex order on $\Delta_+$. In the examples below, let $\Check{\cdot}$ denote the parameters for $\Delta(w, +)$. For instance, $\Check{w} = w w_{\circ}$.

Fortunately, the parameters $y_i$ are all $1$ in our examples below, so we omit the value of $y_i$ in the examples. Also, the numbers of rows $n$ are same for both $\Delta(w, -)$ and $\Delta(w, +)$, in other words $\Check{n} := n$ for all the examples below.

\subsection{Type $A_1^{(1)}$}

Let $\mathfrak{g} = \widehat{\mathfrak{sl}_2}$ be affine algebra of type $A_1^{(1)}$. In this case, there are only two convex orders on $\Delta_+$ except for the order on $\Delta_+^{\mathrm{im}}$, and one of them is just the reverse of the another one. The corresponding parameters are as follows:
\begin{align*}
w &:= 1, n := 1; \quad \mathring{I} = \{1\} = J_0 \supsetneqq J_1 = \emptyset, \\
{\bf s}_0 &:= (s_0 s_1)^{\infty}, \\
\mathring{I} &= \{1\} = \Check{J}_0 \supsetneqq \Check{J}_1 = \emptyset, \\
\Check{\bf s}_0 &:= (s_1 s_0)^{\infty},
\end{align*}
where $({\bf s})^{\infty} := {\bf s}{\bf s}{\bf s}\dots$ denote infinite repetition of ${\bf s}$. Then corresponding convex order $\leq$ turns out to coincide with (\ref{eq:order_sl2}).

Next, we compute the root vectors from this convex order. Since $A_1^{(1)}$ is not simply-laced, we cannot use the algorithm of section \ref{sec:2.5}. But the following formula is sufficient to accomplish the computation.
\begin{equation}
T_1\left(\left[E_0, E_1\right]_q\right) = [E_1, E_0]_q, \quad T_0\left(\left[E_1, E_0\right]_q\right) = [E_0, E_1]_q.
\end{equation}
By definition of root vectors and the chosen order, one can verify
\begin{equation}
\begin{aligned}
E_{\leq, (2n + 1) \delta - \alpha_1} &= (T_0 T_1)^n (E_0) \\
&= \frac{1}{[2]_q^{2n}} \left(\left(\stackrel{\leftarrow}{\mathrm{ad}}\, \left[E_0, E_1\right]_q\right)^{2n} (E_0)\right), \\
E_{\leq, (2n + 2) \delta - \alpha_1} &= (T_0 T_1)^n E_0 (E_1) \\
&= \frac{1}{[2]_q^{2n + 1}} \left(\left(\stackrel{\leftarrow}{\mathrm{ad}}\, \left[E_0, E_1\right]_q\right)^{2n + 1} (E_0)\right), \\
E_{\leq, 2n \delta + \alpha_1} &= \Psi (T_1 T_0)^n (E_1) \\
&= \frac{1}{[2]_q^{2n}} \left(\left(\stackrel{\rightarrow}{\mathrm{ad}}\, \left[E_0, E_1\right]_q\right)^{2n} (E_1)\right), \\
E_{\leq, (2n + 1) \delta + \alpha_1} &= \Psi (T_1 T_0)^n T_1 (E_0) \\
&= \frac{1}{[2]_q^{2n + 1}} \left(\left(\stackrel{\rightarrow}{\mathrm{ad}}\, \left[E_0, E_1\right]_q\right)^{2n + 1} (E_1)\right),
\end{aligned}
\end{equation}
for all $n = 0, 1, 2, \dots$.

Using the reduced expression $t_{\varepsilon_1} = \rho s_1$, where $\rho \in \Omega$ is the transposition of $0$ and $1$, one can show that for any positive integer $m$,
\begin{equation} \label{eq:im_A_1_calc}
\begin{aligned}
\varphi_{1, m} &= \left[\left(T_\rho T_1\right)^{m - 1} \left(E_0\right), E_1\right]_q \\
&= \left[E_{\leq, m\delta - \alpha_1}, E_1\right]_q, \\
T_1 (\varphi_{1, m}) &= T_{\rho^{-1}} (\varphi_{1, m}) \\
&= \left[\Psi E_{\leq, (m - 1)\delta + \alpha_1}, E_0\right]_q.
\end{aligned}
\end{equation}

We set the projection of section \ref{sec:representation} by $\sigma_{01} := +1$. Then corresponding skew-symmetric matrix is $B = \bigl(\begin{smallmatrix} 0 & -2 \\ 2 & 0 \end{smallmatrix}\bigr)$, and matrix presentation of the bilinear form $\langle \cdot, \cdot\rangle_B$ is $(\langle \alpha_i, \alpha_j \rangle_B)_{i, j = 0}^{1} = \bigl(\begin{smallmatrix} 2 & 0 \\ -4 & 2 \end{smallmatrix}\bigr)$.
Since $\langle \alpha_0, \alpha_1 \rangle_B = 0$, the projection $\pi_B^+: U_q^+ \rightarrow \mathcal{P}_B^+$ annihilates $[E_0, E_1]_q$. Thus all the root vectors vanish in $\mathcal{P}_B^+$ except for simple root vectors $E_{\leq, \alpha_i} = E_i$ $(i = 0, 1)$. Therefore, the image of quasi-universal R-matrix $\Theta$ is
\begin{equation} \label{eq_A1_image_theta}
\pi_B^+ \widehat{\otimes} \pi_B^- (\Theta) = \mathbb{E}(y_1) \mathbb{E}(y_0) \quad \in \widehat{\mathcal{D}_B}.
\end{equation}
Beware that the order of product is reverse to given convex order (\ref{eq:prod_formula}).

Now we consider the reversed order $\leq'$, which is in fact obtained by just swapping $\Delta(w, -)$ for $\Delta(w, +)$. The corresponding parameters are also just swapping every parameter $\cdot$ for $\Check{\cdot}$. Thus $E_{\leq', \alpha} = \Psi E_{\leq, \alpha}$ for every real root $\alpha \in \Delta_+^{\mathrm{re}}$. One can verify that all the real root vectors for $\leq'$ satisfy the condition of Proposition \ref{prop:image_real_root_vec} and thus do not vanish.

Since $\Check{w} = s_1$, $T_1 (I_{1, m})$ $(m \geq 1)$ are used as imaginary root vectors. Using (\ref{eq:im_A_1_calc}) we have
\begin{equation}
\pi_B^+ (T_1 \varphi_{1, m}) = [m + 1]_q (q - q^{-1})^{2m - 1} (e_0 e_1)^m.
\end{equation}
Let $D := (q - q^{-1})^2 e_0 e_1$. Then the image of generating function $T_1 \varphi_1 (z) \in U_q^+ [[z]]$ is
\begin{equation}
\pi_B^+ \left(T_1 \left(1 + \varphi_1 \left(z\right)\right)\right) = \sum_{m = 0}^{\infty} [m + 1]_q (Dz)^m = \frac{1}{(1 - qDz) (1 - q^{-1}Dz)},
\end{equation}
where $\pi_B^+: U_q^+ [[z]] \rightarrow \mathcal{P}_B^+[[z]]$ is defined degreewise. Recall that $I_1 (z) = (q - q^{-1}) \sum_{m = 1}^{\infty} I_{1, m} z^m := \log(1 + \varphi_1(z))$. Since $\log (1 + x) = -\sum_{m = 1}^{\infty} (-1)^m x^m / m$ (\ref{eq:def_log}),
\begin{equation}
\pi_B^+ (T_1 I_1(z)) = \sum_{m = 1}^{\infty} \frac{q^m + q^{-m}}{m} D^m z^m
\end{equation}
and therefore
\begin{equation}
\pi_B^+ (I_{1, m}) = \frac{q^m + q^{-m}}{m(q - q^{-1})} D^m.
\end{equation}
Now we compute the image of $S'_m := T_1 \otimes T_1 (S_m)$ (\ref{eq:def_Sn}). By definition $b_{1,1;m} = [2m]_q / (m(q^{-1} - q))$ and $c_{1, 1; m} = b_{1, 1; m}^{-1}$. Thus
\begin{equation*}
S_m = c_{1, 1; m} I_{1, m} \otimes J_{1, m} = \frac{m(q^{-1} - q)}{[2m]_q} I_{1, m} \otimes \Omega I_{1, m} \quad \in U_q^+ \otimes U_q^-.
\end{equation*}
Let $D' := \overline{\Omega_B} (D) = (q^{-1} - q)^2 f_1 f_0$. By virtue of (\ref{eq:omega_bar_commute}) and $\Omega T_i = T_i \Omega$ $(i = 0, 1, \dots, \ell)$, we can compute as follows.
\begin{align*}
\pi_B^+ \otimes \pi_B^- (S'_m) &= \frac{m(q^{-1} - q)}{[2m]_q} \frac{q^m + q^{-m}}{m(q - q^{-1})} \frac{q^m + q^{-m}}{m(q^{-1} - q)} D^m \otimes D'^m \\
&= \frac{1}{m} \frac{q^m + q^{-m}}{q^m - q^{-m}} (D \otimes D')^m = -\frac{q^m (q^m + q^{-m})}{m(1 - q^{2m})} (D \otimes D')^m.
\end{align*}
By (\ref{eq:q_dilog}), we obtain the image of $\Theta_{\mathrm{im}} := \prod_{m = 1}^{\infty} \Theta_{m \delta}$.
\begin{equation}
\begin{aligned}
\pi_B^+ \widehat{\otimes} \pi_B^- (\Theta_{\mathrm{im}}) &= \mathbb{E} (-q D\otimes D')^{-1} \mathbb{E} (-q^{-1} D\otimes D')^{-1} \\
&= \mathbb{E} (-q \normord{y_0 y_1})^{-1} \mathbb{E} (-q^{-1} \normord{y_0 y_1})^{-1}.
\end{aligned}
\end{equation}

Finally, we have
\begin{equation*}
\begin{aligned}
\pi_B^+ \widehat{\otimes} \pi_B^- (\Theta) &= \mathbb{E} (\normord{y_0}) \mathbb{E} (\normord{y_0^2 y_1}) \mathbb{E} (\normord{y_0^3 y_1^2}) \dots \\
&\quad \times \mathbb{E} (-q \normord{y_0 y_1})^{-1} \mathbb{E} (-q^{-1} \normord{y_0 y_1})^{-1} \\
&\quad \quad \times \dots \mathbb{E} (\normord{y_0^2 y_1^3}) \mathbb{E} (\normord{y_0 y_1^2}) \mathbb{E} (\normord{y_1}).
\end{aligned}
\end{equation*}
Comparing with (\ref{eq_A1_image_theta}), we eventually obtain the following quantum dilogarithm identity, which was first found by Terasaki~\cite{Tera}.
\begin{thm}~\cite{Tera}
Let $y_0, y_1$ be indeterminate. Then the following identity holds in skew formal power series algebra $\widehat{\mathcal{S}_1} := \Q(q)[[y_0, y_1]]$ with commutation relation $y_0 y_1 = q^{-4} y_1 y_0$.
\begin{equation} \label{eq:A1_result}
\begin{aligned}
\mathbb{E}(\normord{y_1}) \mathbb{E}(\normord{y_0}) &= \mathbb{E} (\normord{y_0}) \mathbb{E} (\normord{y_0^2 y_1}) \mathbb{E} (\normord{y_0^3 y_1^2}) \dots \\
&\quad \times \mathbb{E} (-q \normord{y_0 y_1})^{-1} \mathbb{E} (-q^{-1} \normord{y_0 y_1})^{-1} \\
&\quad \quad \times \dots \mathbb{E} (\normord{y_0^2 y_1^3}) \mathbb{E} (\normord{y_0 y_1^2}) \mathbb{E} (\normord{y_1}),
\end{aligned}
\end{equation}
where $\normord{y_0^{m_0} y_1^{m_1}} = q^{2 m_0 m_1} y_0^{m_0} y_1^{m_1}$.
\end{thm}

Let $B' = \bigl(\begin{smallmatrix} 0 & 1 \\ -1 & 0 \end{smallmatrix}\bigr)$, $\widehat{\mathcal{S}} := \widehat{\mathcal{S}_{B'}} \cong \Q(q)[[x_1, x_2]]$. Then $x_1 x_2 = q^2 x_2 x_1$, which coincides with the commutation relation in the introduction. If we set $R := \bigl(\begin{smallmatrix} 0 & 2 \\ 1 & 3 \end{smallmatrix}\bigr)$, $R$ satisfies ${}^t\!R B' R = B$. Thus by Proposition \ref{prop:S_B_hom} and \ref{prop:R_conti}, there exists unique continuous algebra homomorphism $\widehat{\psi_1}: \widehat{\mathcal{S}_1} \rightarrow \widehat{\mathcal{S}}$ satisfying 
\begin{equation}
\widehat{\psi_1}(y_0) = x_2, \quad \widehat{\psi_1}(y_1) = \normord{x_1^2 x_2^3} = q^{-6} x_1^2 x_2^3.
\end{equation}

Let $\mathcal{L} := \mathcal{L}_{B'} \cong \Q(q)[x_1^{\pm 1}, x_2^{\pm 1}]$. Since $S := \bigl(\begin{smallmatrix} 1 & 0 \\ -2 & 1 \end{smallmatrix}\bigr)$ satisfies ${}^t\!S B' S = B'$, we have algebra automorphism $\psi_S \in \mathrm{Aut}\, \mathcal{L}$. $\psi_S$ transforms variables $x_1, x_2$ as $\psi_S(x_1) = \normord{x_1 x_2^{-2}} = q^2 x_1 x_2^{-2}$, $\psi_S(x_2) = x_2$.

Applying $\widehat{\psi_1}$ on (\ref{eq:A1_result}) and transforming variables $x_1, x_2$ by $\psi_S$, we obtain (\ref{eq:A1}) in $\Q(q)[[\frac{x_1}{x_2^2}, x_2]]$. Recall that $\widehat{\psi_1}$ and $\psi_S$ preserves the normal ordered product (Proposition \ref{prop:preserve_normord}).

\begin{cor}
The identity (\ref{eq:A1}) holds in $\Q(q)[[\frac{x_1}{x_2^2}, x_2]]$.
\end{cor}

\textbf{Remark.} We shall call a group homomorphism $Z: Q \rightarrow \C$ central charge. When $Z$ is injective and $Z(\Delta_+)$ lies in the (closure of) upper half plane $\overline{\mathcal{H}} := \left\{\, z \in \C \relmiddle| \mathrm{Im}\, z \geq 0\,\right\}$,
\begin{equation} \label{eq:order_by_central_charge}
\alpha \leq_Z \beta \stackrel{\mathrm{def}}{\Leftrightarrow} \arg Z(\alpha) \leq \arg Z(\beta) \quad (\alpha, \beta \in \Delta_+^{\mathrm{re}})
\end{equation}
defines a convex order on positive real roots, where we choose principal value of argument so that $0 \leq \arg z < 2\pi$.

Setting $Z(\alpha_0) := 1$, $Z(\alpha_1) := 1 + \sqrt{-1}$ yields the convex order (\ref{eq:order_sl2}). Notice that \textbf{$\leq_Z$ yields only convex order of single row}, because every root is of the form $m\delta + \alpha$ $(m \in \Z, \alpha \in \mathring{\Delta} \cup \{0\})$ and thus $Z(\Delta)$ lies in finite number of lines parallel to $Z(\delta)$.

\subsection{Type $A_2^{(1)}$}

Let $\mathfrak{g} = \widehat{\mathfrak{sl}_3}$ be affine algebra of type $A_2^{(1)}$. We choose a convex order by setting
\begin{align*}
w &:= s_1, \quad n := 2; \\
\mathring{I} &= \{1, 2\} = J_0 \supsetneqq J_1 := \{1\} \supsetneqq J_2 = \emptyset, \\
{\bf s}_0 &:= (s_0 s_1 s_2)^{\infty}, \quad {\bf s}_1 := (s_1 s_{\delta - \alpha_1})^{\infty}, \\
\mathring{I} &= \Check{J}_0 \supsetneqq \Check{J}_1 := \{2\} \supsetneqq \Check{J}_2 = \emptyset, \\
\Check{\bf s}_0 &:= (s_2 s_1 s_0)^{\infty}, \quad \Check{\bf s}_1 := (s_{\delta - \alpha_2} s_2)^{\infty}. 
\end{align*}
Then the corresponding convex order $\leq$ is
\begin{equation*}
\begin{aligned}
\delta - \alpha_1 - \alpha_2 &< \delta - \alpha_2 < 2\delta - \alpha_1 - \alpha_2 < 2\delta - \alpha_2 < \dots \\
< \alpha_1 &< \delta + \alpha_1 < 2\delta + \alpha_1 < 3\delta + \alpha_1 < \dots \\
< \delta &< 2\delta < 3\delta < 4\delta < \dots \\
\dots &< 3\delta - \alpha_1 < 2\delta - \alpha_1 < \delta - \alpha_1 \\
\dots &< \delta + \alpha_1 + \alpha_2 < \delta + \alpha_2 < \alpha_1 + \alpha_2 < \alpha_2,
\end{aligned}
\end{equation*}
where the null root $\delta = \alpha_0 + \alpha_1 + \alpha_2$.

Using the algorithm of Proposition \ref{prop:algorithm} and notation (\ref{eq:tree_notation}), real root vectors in the first row of $\Delta(w, -)$ are computed as follows.
\begin{equation}
E_{\leq, m\delta - \alpha_1 - \alpha_2} = \begin{xy}
(-22, 3) ; (-18, 4) **\crv{(-21, 3.49) & (-20, 3.5) & (-19, 3.51)},
(-14, 3) ; (-18, 4) **\crv{(-15, 3.49) & (-16, 3.5) & (-17, 3.51)},
(-13, 3) ; (-9, 4) **\crv{(-12, 3.49) & (-11, 3.5) & (-10, 3.51)},
(-5, 3) ; (-9, 4) **\crv{(-6, 3.49) & (-7, 3.5) & (-8, 3.51)},
(5, 3) ; (9, 4) **\crv{(6, 3.49) & (7, 3.5) & (8, 3.51)},
(13, 3) ; (9, 4) **\crv{(12, 3.49) & (11, 3.5) & (10, 3.51)},
(14, 3) ; (18, 4) **\crv{(15, 3.49) & (16, 3.5) & (17, 3.51)},
(22, 3) ; (18, 4) **\crv{(21, 3.49) & (20, 3.5) & (19, 3.51)},
(-22, 4) ; (0, 5) **\crv{(-21, 4.49) & (-11, 4.5) & (-1, 4.51)},
(22, 4) ; (0, 5) **\crv{(21, 4.49) & (11, 4.5) & (1, 4.51)},
(-24, 2) * {0}, (-21, 2) * {1}, (-18, 2) * {0}, (-15, 2) * {2}, (-12, 2) * {1}, (-9, 2) * {0}, (-6, 2) * {2}, (0, 2) * {\cdots},
(  6, 2) * {1}, (  9, 2) * {0}, ( 12, 2) * {2}, ( 15, 2) * {1}, ( 18, 2) * {0}, (21, 2) * {2}, (0, 6) * {m - 1},
\ar @{-} (-24, 0) ; (-1.5, -7.5),
\ar @{-} (-21, 0) ; (-22.5,-0.5),
\ar @{-} (-18, 0) ; (-16.5, -0.5),
\ar @{-} (-15, 0) ; (-19.5, -1.5),
\ar @{-} (-12, 0) ; (-18, -2),
\ar @{-} (-9, 0) ; (-7.5, -0.5),
\ar @{-} (-6, 0) ; (-15, -3),
\ar @{-} (6, 0) ; (-9, -5),
\ar @{-} (9, 0) ; (10.5, -0.5),
\ar @{-} (12, 0) ; (-6, -6),
\ar @{-} (15, 0) ; (-4.5, -6.5),
\ar @{-} (18, 0) ; (19.5, -0.5),
\ar @{-} (21, 0) ; (-1.5, -7.5)
\end{xy},
\end{equation}
\begin{equation}
E_{\leq, m\delta - \alpha_2} = \begin{xy}
(-22, 3) ; (-18, 4) **\crv{(-21, 3.49) & (-20, 3.5) & (-19, 3.51)},
(-14, 3) ; (-18, 4) **\crv{(-15, 3.49) & (-16, 3.5) & (-17, 3.51)},
(-13, 3) ; (-9, 4) **\crv{(-12, 3.49) & (-11, 3.5) & (-10, 3.51)},
(-5, 3) ; (-9, 4) **\crv{(-6, 3.49) & (-7, 3.5) & (-8, 3.51)},
(5, 3) ; (9, 4) **\crv{(6, 3.49) & (7, 3.5) & (8, 3.51)},
(13, 3) ; (9, 4) **\crv{(12, 3.49) & (11, 3.5) & (10, 3.51)},
(14, 3) ; (18, 4) **\crv{(15, 3.49) & (16, 3.5) & (17, 3.51)},
(22, 3) ; (18, 4) **\crv{(21, 3.49) & (20, 3.5) & (19, 3.51)},
(-22, 4) ; (0, 5) **\crv{(-21, 4.49) & (-11, 4.5) & (-1, 4.51)},
(22, 4) ; (0, 5) **\crv{(21, 4.49) & (11, 4.5) & (1, 4.51)},
(-24, 2) * {0}, (-21, 2) * {1}, (-18, 2) * {0}, (-15, 2) * {2}, (-12, 2) * {1}, (-9, 2) * {0}, (-6, 2) * {2}, (0, 2) * {\cdots},
(  6, 2) * {1}, (  9, 2) * {0}, ( 12, 2) * {2}, ( 15, 2) * {1}, ( 18, 2) * {0}, (21, 2) * {2}, (24, 2) * {1}, (0, 6) * {m - 1},
\ar @{-} (-24, 0) ; (0,-8),
\ar @{-} (-21, 0) ; (-22.5,-0.5),
\ar @{-} (-18, 0) ; (-16.5, -0.5),
\ar @{-} (-15, 0) ; (-19.5, -1.5),
\ar @{-} (-12, 0) ; (-18, -2),
\ar @{-} (-9, 0) ; (-7.5, -0.5),
\ar @{-} (-6, 0) ; (-15, -3),
\ar @{-} (6, 0) ; (-9, -5),
\ar @{-} (9, 0) ; (10.5, -0.5),
\ar @{-} (12, 0) ; (-6, -6),
\ar @{-} (15, 0) ; (-4.5, -6.5),
\ar @{-} (18, 0) ; (19.5, -0.5),
\ar @{-} (21, 0) ; (-1.5, -7.5),
\ar @{-} (24, 0) ; (0, -8)
\end{xy} \quad (m \geq 1). 
\end{equation}
This can be directly proven by induction on $m$, noting that 
\begin{equation}
T_0 T_1 T_2([E_0, E_2]_q) = E_1, \quad T_0 T_1 T_2(E_1) = [E_0, E_2]_q,
\end{equation}
and $T_0 T_1 T_2([E_0, E_1]_q) = E_{\leq, 3\delta - \alpha_1 - \alpha_2}$, which means that applying $T_0 T_1 T_2$ on $[E_0, E_1]_q$ adds 3 branches from right.

The computation of real root vectors in the second row requires more preparation. First, we have to compute $E_{\delta - \alpha_1}$ of (\ref{eq:sub_simple_root_vector}). Since $\Phi(s_0 s_2) = \{\alpha_0, \alpha_0 + \alpha_2\} \subset \Delta(1, -)$, we have $E_{\delta - \alpha_1} = T_0 (E_2) = [E_0, E_2]_q$. Next, we need to compute $\widehat{s_{\delta - \alpha_1}} \in \widehat{W}$, which is appropriate extension of $s_{\delta - \alpha_1} \in W_{J_1}$. By the definition (\ref{eq:simple_ref_alternative}), $\widehat{s_{\delta - \alpha_1}} = t_{\varepsilon_1}^{J_1} s_1 t_{\varepsilon_1}^{J_1}$. To compute this, we require reduced expression of $t_{\varepsilon_1} \in \widehat{W}$. Let $\rho \in \Omega$ denote the Dynkin automorphism which acts on the indices as $\rho(0) = 1$, $\rho(1) = 2$, $\rho(2) = 0$. Then
\begin{equation}
t_{\varepsilon_1} = \rho s_2 s_1, \quad t_{\varepsilon_2} = \rho^2 s_1 s_2 \quad \in \widehat{W} \subset \mathrm{GL}(\mathfrak{h'}^*)
\end{equation}
are reduced expressions. Note that the length of $w \in \widehat{W}$ defined by (\ref{eq:length_in_extW}) coincides with the number of positive roots $\alpha \in \Delta_+$ such that $w(\alpha) \in \Delta_-$. More generally, 

\begin{prop} \label{prop:translation_reduced_rep}
Let $\mathfrak{g} = \widehat{\mathfrak{sl}_{\ell + 1}}$ be affine algebra of type $A_{\ell}^{(1)}$. We set the indices $0, 1, \dots, \ell$ so that $a_{i i + 1} \neq 0$ for $i = 0, 1, \dots, \ell - 1$. Let $\rho \in \Omega$ be the Dynkin automorphism which acts on indices as $\rho(i) = i + 1$ ($0 \leq i \leq \ell - 1$), $\rho(\ell) = 0$. Then
\begin{equation} \label{eq:translation_reduced_rep}
t_{\varepsilon_i} = \left(\rho^{-1} s_1 s_2 \dots s_i\right)^{\ell + 1 - i} \quad (i = 1, 2, \dots, \ell)
\end{equation}
are reduced expressions in $\widehat{W}$.
\end{prop}
\begin{proof}
Firstly, we have to check the equality. It is enough to compare the action of both sides on simple roots since $\widehat{W} \subset \mathrm{GL}\,(\mathfrak{h}'^*)$. Moreover, since every element of $\widehat{W}$ fixes the null root $\delta = \alpha_0 + \alpha_1 + \dots + \alpha_{\ell}$, it is enough to check the action on $\alpha_j$ for $j = 1, 2, \dots, \ell$. On the one hand, $t_{\varepsilon_i} (\alpha_j) = \alpha_j - \delta_{ij} \delta$. On the other hand, let $R_i := \rho^{-1} s_1 s_2 \dots s_i$. When $\ell = 1$, $R_1(\alpha_1) = \rho^{-1}(-\alpha_1) = - \alpha_0 = \alpha_1 - \delta$ and hence (\ref{eq:translation_reduced_rep}) holds. Next, we assume $\ell \geq 2$. Recall that for $1 \leq i, j \leq \ell$,
\begin{equation}
s_i(\alpha_j) = \begin{cases}
\alpha_j & |i - j| > 1 \\
\alpha_i + \alpha_j & j = i \pm 1 \\
- \alpha_i & j = i
\end{cases}
\end{equation}
in the root system of type $A_{\ell}$. By direct calculation we have
\begin{equation}
R_i(\alpha_j) = \begin{cases}
\alpha_j & 1 \leq j \leq i - 1\\
-\delta + \alpha_i + \alpha_{i + 1} + \dots + \alpha_{\ell} & j = i \\
\delta - \alpha_{i + 1} - \alpha_{i + 2} - \dots - \alpha_{\ell} & j = i + 1 \\
\alpha_{j - 1} & j > i + 1
\end{cases}.
\end{equation}
When $1 \leq j < i$, it is clear that $R_i^{\ell + 1 - i} (\alpha_j) = \alpha_j$.

When $j = i$, the case $i = \ell$ is the above formula. When $i < \ell$, notice that $R_i(\alpha_i + \alpha_{i + 1}) = \alpha_i$. Using this inductively,
\begin{align*}
R_i^{\ell + 1 - i}(\alpha_i) &= R_i^{\ell - i}(-\delta + \alpha_i + \alpha_{i + 1} + \dots + \alpha_{\ell}) \\
&= R_i^{\ell - i - 1}(-\delta + \alpha_i + \alpha_{i + 1} + \dots + \alpha_{\ell - 1}) \\
&\dots \\
&= R_i(-\delta + \alpha_i + \alpha_{i + 1}) = -\delta + \alpha_i
\end{align*}
and thus $R_i^{\ell + 1 - i}(\alpha_i) = -\delta + \alpha_i$.

When $j = i + 1$, notice that $R_i^2(\alpha_{i + 1}) = \alpha_{\ell}$. Thus
\begin{equation*}
R_i^{\ell + 1 - i}(\alpha_{i + 1}) = R_i^{\ell - i - 1} (\alpha_{\ell}) = \alpha_{i + 1}.
\end{equation*}

When $j > i + 1$, $R_i^{\ell + 1 - i}(\alpha_j) = R_i^{\ell - j + 2} (\alpha_{i + 1}) = R_i^{\ell - j} (\alpha_{\ell}) = \alpha_j$.

By above calculation, we conclude $t_{\varepsilon_i} = R_i^{\ell + 1 - i}$.

To verify that $R_i^{\ell + 1 - i}$ is reduced expression, it is enough to show that the length of $t_{\varepsilon_i} \in \widehat{W}$ is $i(\ell + 1 - i)$. The length of $t_{\varepsilon_i}$ coincides with the number of positive roots which $t_{\varepsilon_i}$ sends to negative roots. Recall that 
\begin{equation}
\mathring{\Delta} = \left\{\,\pm\left(\alpha_i + \alpha_{i + 1} + \dots + \alpha_j\right) \relmiddle| 1 \leq i \leq j \leq \ell\,\right\}
\end{equation}
in finite root system of type $A_{\ell}$. $t_{\varepsilon_i}$ translates the roots containing $\pm\alpha_i$ by $\mp\delta$. Thus, if $\alpha = m\delta + \varepsilon \in \Delta_+$ ($m \in \Z_{\geq 0}$, $\varepsilon \in \mathring{\Delta}$) satisfies $t_{\varepsilon_i}(\alpha) \in \Delta_-$, then $m = 0$ and $\varepsilon$ must be a positive root containing $\alpha_i$. Such $\varepsilon$ takes the form $\alpha_j + \alpha_{j + 1} + \dots + \alpha_k$ ($j \leq i \leq k$), and the number of such $(j, k)$ is $i(\ell - i + 1)$. This shows that the length of $t_{\varepsilon_i}$ is $i(\ell + 1 - i)$.
\end{proof}

By Proposition \ref{prop:translation_reduced_rep},  $t_{\varepsilon_1}^{J_1} = (\rho^{-1} s_1 \rho^{-1} s_1)^{J_1} = (\rho s_2 s_1)^{J_1} = \rho s_2$ and we have
\begin{equation}
\widehat{s_{\delta - \alpha_1}} = \rho s_2 s_1 \rho s_2.
\end{equation}

Now we can compute real root vectors in the second row. Since $w^{J_1} = s_1^{J_1} = 1$ and $E_{\delta - \alpha_1} = T_0(E_1) = T_{\rho} T_2(E_1)$,
\begin{equation*}
\begin{aligned}
E_{\leq, m\delta + \alpha_1} &= \begin{cases}
(T_1 T_{\widehat{s_{\delta - \alpha_1}}})^{m/2} (E_1) & m: \text{even} \\
(T_1 T_{\widehat{s_{\delta - \alpha_1}}})^{(m - 1)/2} T_1(E_{\delta - \alpha_1}) & m: \text{odd}
\end{cases} \\
&= (T_1 T_{\rho} T_2)^m (E_1) \quad (m \geq 0).
\end{aligned}
\end{equation*}
Observing $T_1 T_{\rho} T_2 (E_0) = E_0$ and $T_1 T_{\rho} T_2([E_1, E_2]_q) =[E_1, E_2]_q$, we have
\begin{equation}
E_{\leq, m\delta + \alpha_1} = \begin{xy}
(-22, 3) ; (-18, 4) **\crv{(-21, 3.49) & (-20, 3.5) & (-19, 3.51)},
(-14, 3) ; (-18, 4) **\crv{(-15, 3.49) & (-16, 3.5) & (-17, 3.51)},
(-13, 3) ; (-9, 4) **\crv{(-12, 3.49) & (-11, 3.5) & (-10, 3.51)},
(-5, 3) ; (-9, 4) **\crv{(-6, 3.49) & (-7, 3.5) & (-8, 3.51)},
(5, 3) ; (9, 4) **\crv{(6, 3.49) & (7, 3.5) & (8, 3.51)},
(13, 3) ; (9, 4) **\crv{(12, 3.49) & (11, 3.5) & (10, 3.51)},
(14, 3) ; (18, 4) **\crv{(15, 3.49) & (16, 3.5) & (17, 3.51)},
(22, 3) ; (18, 4) **\crv{(21, 3.49) & (20, 3.5) & (19, 3.51)},
(-22, 4) ; (0, 5) **\crv{(-21, 4.49) & (-11, 4.5) & (-1, 4.51)},
(22, 4) ; (0, 5) **\crv{(21, 4.49) & (11, 4.5) & (1, 4.51)},
(-24, 2) * {1}, (-21, 2) * {0}, (-18, 2) * {1}, (-15, 2) * {2}, (-12, 2) * {0}, (-9, 2) * {1}, (-6, 2) * {2}, (0, 2) * {\cdots},
(  6, 2) * {0}, (  9, 2) * {1}, ( 12, 2) * {2}, ( 15, 2) * {0}, ( 18, 2) * {1}, (21, 2) * {2}, (0, 6) * {m},
\ar @{-} (-24, 0) ; (-1.5, -7.5),
\ar @{-} (-21, 0) ; (-22.5,-0.5),
\ar @{-} (-18, 0) ; (-16.5, -0.5),
\ar @{-} (-15, 0) ; (-19.5, -1.5),
\ar @{-} (-12, 0) ; (-18, -2),
\ar @{-} (-9, 0) ; (-7.5, -0.5),
\ar @{-} (-6, 0) ; (-15, -3),
\ar @{-} (6, 0) ; (-9, -5),
\ar @{-} (9, 0) ; (10.5, -0.5),
\ar @{-} (12, 0) ; (-6, -6),
\ar @{-} (15, 0) ; (-4.5, -6.5),
\ar @{-} (18, 0) ; (19.5, -0.5),
\ar @{-} (21, 0) ; (-1.5, -7.5)
\end{xy} \quad (m \geq 0).
\end{equation}

Similarly, root vectors for real roots in $\Delta(w, +) = \Delta(s_2 s_1, -)$ are computed as follows. Comparing ${\bf s}_0 = (s_0 s_1 s_2)^{\infty}$ and $\Check{\bf s}_0 = (s_2 s_1 s_0)^{\infty}$, $E_{\leq, m\delta + \alpha_2}$, $E_{\leq, m\delta + \alpha_1 + \alpha_2}$ are obtained by swapping all the index $0$ for $2$ in $E_{\leq, m\delta - \alpha_1 - \alpha_2}$, $E_{\leq, m\delta - \alpha_2}$ and applying $\Psi$, which just reverses all the directions of $q$-bracket. As a result,
\begin{equation}
E_{\leq, m\delta + \alpha_2} = \begin{xy}
(-22, 3) ; (-18, 4) **\crv{(-21, 3.49) & (-20, 3.5) & (-19, 3.51)},
(-14, 3) ; (-18, 4) **\crv{(-15, 3.49) & (-16, 3.5) & (-17, 3.51)},
(-13, 3) ; (-9, 4) **\crv{(-12, 3.49) & (-11, 3.5) & (-10, 3.51)},
(-5, 3) ; (-9, 4) **\crv{(-6, 3.49) & (-7, 3.5) & (-8, 3.51)},
(5, 3) ; (9, 4) **\crv{(6, 3.49) & (7, 3.5) & (8, 3.51)},
(13, 3) ; (9, 4) **\crv{(12, 3.49) & (11, 3.5) & (10, 3.51)},
(14, 3) ; (18, 4) **\crv{(15, 3.49) & (16, 3.5) & (17, 3.51)},
(22, 3) ; (18, 4) **\crv{(21, 3.49) & (20, 3.5) & (19, 3.51)},
(-22, 4) ; (0, 5) **\crv{(-21, 4.49) & (-11, 4.5) & (-1, 4.51)},
(22, 4) ; (0, 5) **\crv{(21, 4.49) & (11, 4.5) & (1, 4.51)},
(-21, 2) * {0}, (-18, 2) * {2}, (-15, 2) * {1}, (-12, 2) * {0}, (-9, 2) * {2}, (-6, 2) * {1}, (0, 2) * {\cdots},
(  6, 2) * {0}, (  9, 2) * {2}, ( 12, 2) * {1}, ( 15, 2) * {0}, ( 18, 2) * {2}, (21, 2) * {1}, (24, 2) * {2}, (0, 6) * {m},
\ar @{-} (-21, 0) ; (1.5,-7.5),
\ar @{-} (-18, 0) ; (-19.5, -0.5),
\ar @{-} (-15, 0) ; (4.5, -6.5),
\ar @{-} (-12, 0) ; (6, -6),
\ar @{-} (-9, 0) ; (-10.5, -0.5),
\ar @{-} (-6, 0) ; (9, -5),
\ar @{-} (6, 0) ; (15, -3),
\ar @{-} (9, 0) ; (7.5, -0.5),
\ar @{-} (12, 0) ; (18, -2),
\ar @{-} (15, 0) ; (19.5, -1.5),
\ar @{-} (18, 0) ; (16.5, -0.5),
\ar @{-} (21, 0) ; (22.5, -0.5),
\ar @{-} (24, 0) ; (1.5, -7.5)
\end{xy},
\end{equation}
\begin{equation}
E_{\leq, m\delta + \alpha_1 + \alpha_2} = \begin{xy}
(-22, 3) ; (-18, 4) **\crv{(-21, 3.49) & (-20, 3.5) & (-19, 3.51)},
(-14, 3) ; (-18, 4) **\crv{(-15, 3.49) & (-16, 3.5) & (-17, 3.51)},
(-13, 3) ; (-9, 4) **\crv{(-12, 3.49) & (-11, 3.5) & (-10, 3.51)},
(-5, 3) ; (-9, 4) **\crv{(-6, 3.49) & (-7, 3.5) & (-8, 3.51)},
(5, 3) ; (9, 4) **\crv{(6, 3.49) & (7, 3.5) & (8, 3.51)},
(13, 3) ; (9, 4) **\crv{(12, 3.49) & (11, 3.5) & (10, 3.51)},
(14, 3) ; (18, 4) **\crv{(15, 3.49) & (16, 3.5) & (17, 3.51)},
(22, 3) ; (18, 4) **\crv{(21, 3.49) & (20, 3.5) & (19, 3.51)},
(-22, 4) ; (0, 5) **\crv{(-21, 4.49) & (-11, 4.5) & (-1, 4.51)},
(22, 4) ; (0, 5) **\crv{(21, 4.49) & (11, 4.5) & (1, 4.51)},
(-24, 2) * {1}, (-21, 2) * {0}, (-18, 2) * {2}, (-15, 2) * {1}, (-12, 2) * {0}, (-9, 2) * {2}, (-6, 2) * {1}, (0, 2) * {\cdots},
(  6, 2) * {0}, (  9, 2) * {2}, ( 12, 2) * {1}, ( 15, 2) * {0}, ( 18, 2) * {2}, (21, 2) * {1}, (24, 2) * {2}, (0, 6) * {m},
\ar @{-} (-24, 0) ; (0,-8),
\ar @{-} (-21, 0) ; (1.5,-7.5),
\ar @{-} (-18, 0) ; (-19.5, -0.5),
\ar @{-} (-15, 0) ; (4.5, -6.5),
\ar @{-} (-12, 0) ; (6, -6),
\ar @{-} (-9, 0) ; (-10.5, -0.5),
\ar @{-} (-6, 0) ; (9, -5),
\ar @{-} (6, 0) ; (15, -3),
\ar @{-} (9, 0) ; (7.5, -0.5),
\ar @{-} (12, 0) ; (18, -2),
\ar @{-} (15, 0) ; (19.5, -1.5),
\ar @{-} (18, 0) ; (16.5, -0.5),
\ar @{-} (21, 0) ; (22.5, -0.5),
\ar @{-} (24, 0) ; (0, -8)
\end{xy} \quad (m \geq 0).
\end{equation}
Since $\Check{w}^{\Check{J}_1} = s_2 s_1$ and $\Check{\bf s}_1 = (s_{\delta - \alpha_2} s_2)^{\infty}$, root vectors for $m\delta - \alpha_1$ are defined as
\begin{equation}
E_{\leq, m\delta - \alpha_1} = \begin{cases}
\Psi T_2 T_1 (T_{\widehat{s_{\delta - \alpha_2}}} T_2)^{(m - 1) / 2} (E_{\delta - \alpha_2}) & m: \text{odd} \\
\Psi T_2 T_1 (T_{\widehat{s_{\delta - \alpha_2}}} T_2)^{(m - 2) / 2} T_{\widehat{s_{\delta - \alpha_2}}} (E_2) & m: \text{even}
\end{cases} \quad (m \geq 1),
\end{equation}
where $E_{\delta - \alpha_2} = T_0(E_1)$. By Proposition \ref{prop:translation_reduced_rep}, $t_{\varepsilon_2}^{\Check{J}_1} = (\rho^{-1} s_1 s_2)^{\Check{J}_1} = \rho^2 s_1$ and thus $\widehat{s_{\delta - \alpha_2}} := t_{\varepsilon_2}^{\Check{J}_1} s_2 t_{\varepsilon_2}^{\Check{J}_1}= \rho^2 s_1 s_2 \rho^2 s_1$. Rewriting $E_{\delta - \alpha_2} = T_0(E_1) = T_{\rho^2} T_1 (E_2)$,
\begin{equation}
E_{\leq, m\delta - \alpha_1} = \Psi T_2 T_1 (T_{\rho^2} T_1 T_2)^{m - 1} T_{\rho^2} T_1 (E_2) \quad (m \geq 1).
\end{equation}
Moreover, one can deform this presentation to the form $\Psi (T_u)^m T_v (E_i)$ by realizing $T_{\varepsilon_1}^{-1} (E_2) = E_2$. Since $T_{\varepsilon_1}^{-1} = T_1^{-1} T_2^{-1} T_{\rho}^{-1}$, we can replace $E_2$ with $T_1^{-1} T_2^{-1} T_{\rho}^{-1} (E_2)$. In the extended braid group $\widehat{B}$, 
\begin{align*}
T_{\rho^2} T_1 T_{\varepsilon_1}^{-1} &= T_1^{-1} T_{\rho}, \\
T_{\rho^2} T_1 T_2 \cdot T_1^{-1} T_{\rho} &= T_1^{-1} T_{\rho} \cdot T_{\rho^2} T_0 T_1,
\end{align*}
where we used $T_1^{-1} T_{\rho^2} = T_{\rho^2} T_2^{-1}$ and $T_1 T_2 T_1^{-1} = T_2^{-1} T_1 T_2$. Therefore we have
\begin{equation*}
T_2 T_1 (T_{\rho^2} T_1 T_2)^{m - 1} T_{\rho^2} T_1 (E_2) = T_2 T_{\rho} (T_{\rho^2} T_0 T_1)^{m - 1} (E_2).
\end{equation*}
Finally, rewriting $E_2 = T_{\rho}^{-1} (E_0)$ yields
\begin{align*}
T_2 T_{\rho} (T_{\rho^2} T_0 T_1)^{m - 1} (E_2) &= T_2 T_{\rho} (T_{\rho^2} T_0 T_1)^{m - 1} T_{\rho}^{-1} (E_0) \\
&= T_2 (T_{\rho^2} T_1 T_2)^{m - 1} (E_0) \\
&= (T_2 T_{\rho^2} T_1)^{m - 1} T_2 (E_0).
\end{align*}
As a result, we obtain
\begin{equation}
E_{\leq, m\delta - \alpha_1} = \Psi (T_2 T_{\rho^2} T_1)^{m - 1} T_2 (E_0) \quad (m \geq 1).
\end{equation}

The advantage of this presentation is that inductive computation becomes easy. In fact,
\begin{equation}
T_2 T_{\rho^2} T_1 T_2 (E_0) = \begin{xy}
(-6, 2) * {2}, (-3, 2) * {0}, (0, 2) * {2}, (3, 2) * {1}, (6, 2) * {0},
\ar @{-} (-6, 0) ; (0, -4),
\ar @{-} (-3, 0) ; (-4.5, -1),
\ar @{-} (0, 0) ; (1.5, -1),
\ar @{-} (3, 0) ; (-1.5, -3),
\ar @{-} (6, 0) ; (0, -4)
\end{xy}
\end{equation}
and thus applying $T_2 T_{\rho^2} T_1$ to $[E_2, E_0]_q$ adds 2 branches from right. By virtue of $[E_2, E_1]_q$ and $E_0$ being invariant by $T_2 T_{\rho^2} T_1$, finally we have
\begin{equation}
E_{\leq, m\delta - \alpha_1} = \begin{xy}
(-22, 3) ; (-18, 4) **\crv{(-21, 3.49) & (-20, 3.5) & (-19, 3.51)},
(-14, 3) ; (-18, 4) **\crv{(-15, 3.49) & (-16, 3.5) & (-17, 3.51)},
(-13, 3) ; (-9, 4) **\crv{(-12, 3.49) & (-11, 3.5) & (-10, 3.51)},
(-5, 3) ; (-9, 4) **\crv{(-6, 3.49) & (-7, 3.5) & (-8, 3.51)},
(5, 3) ; (9, 4) **\crv{(6, 3.49) & (7, 3.5) & (8, 3.51)},
(13, 3) ; (9, 4) **\crv{(12, 3.49) & (11, 3.5) & (10, 3.51)},
(14, 3) ; (18, 4) **\crv{(15, 3.49) & (16, 3.5) & (17, 3.51)},
(22, 3) ; (18, 4) **\crv{(21, 3.49) & (20, 3.5) & (19, 3.51)},
(-22, 4) ; (0, 5) **\crv{(-21, 4.49) & (-11, 4.5) & (-1, 4.51)},
(22, 4) ; (0, 5) **\crv{(21, 4.49) & (11, 4.5) & (1, 4.51)},
(-21, 2) * {0}, (-18, 2) * {1}, (-15, 2) * {2}, (-12, 2) * {0}, (-9, 2) * {1}, (-6, 2) * {2}, (0, 2) * {\cdots},
(  6, 2) * {0}, (  9, 2) * {1}, ( 12, 2) * {2}, ( 15, 2) * {0}, ( 18, 2) * {1}, (21, 2) * {2}, (24, 2) * {0}, (27, 2) * {2}, (0, 6) * {m - 1},
\ar @{-} (-21, 0) ; (3,-8),
\ar @{-} (-18, 0) ; (4.5, -7.5),
\ar @{-} (-15, 0) ; (-16.5, -0.5),
\ar @{-} (-12, 0) ; (7.5, -6.5),
\ar @{-} (-9, 0) ; (9, -6),
\ar @{-} (-6, 0) ; (-7.5, -0.5),
\ar @{-} (6, 0) ; (16.5, -3.5),
\ar @{-} (9, 0) ; (18, -3),
\ar @{-} (12, 0) ; (10.5, -0.5),
\ar @{-} (15, 0) ; (21, -2),
\ar @{-} (18, 0) ; (22.5, -1.5),
\ar @{-} (21, 0) ; (19.5, -0.5),
\ar @{-} (24, 0) ; (25.5, -0.5),
\ar @{-} (27, 0) ; (3, -8)
\end{xy} \quad (m \geq 1).
\end{equation}

Next, we compute imaginary root vectors. Since $w = s_1$, $T_1 (I_{i, m})$ $(i = 1, 2; m \in \Z_{\geq 1})$ are used as imaginary root vectors. We use (\ref{eq:invert}) to compute $T_1(\varphi_{1, m})$. Since $t_{\varepsilon_1} = \rho s_2 s_1 = u s_1$, $u = \rho s_2$. Thus
\begin{align*}
T_1(\varphi_{1, m}) &= \left[T_1\left(T_{\rho} T_2 T_1\right)^{m - 1} T_1^{-1}\left(E_1\right), T_u^{-1} \left(E_1\right)\right]_q \\
&= \left[\left(T_1 T_{\rho} T_2\right)^{m - 1} \left(E_1\right), \Psi T_2 T_{\rho}^{-1} \left(E_1\right)\right]_q.
\end{align*}
Since $T_1 T_{\rho} T_2 (E_1) = \begin{xy}
(-6, 2) * {1}, (-3, 2) * {0}, (0, 2) * {1}, (3, 2) * {2},
\ar @{-} (-6, 0) ; (-1.5, -3),
\ar @{-} (-3, 0) ; (-4.5, -1),
\ar @{-} (0, 0) ; (1.5, -1),
\ar @{-} (3, 0) ; (-1.5, -3)
\end{xy}
$ and $T_1 T_{\rho} T_2$ fixes $E_0$ and $[E_1, E_2]_q$,
\begin{equation}
T_1(\varphi_{1, m}) = \begin{xy}
(-22, 3) ; (-18, 4) **\crv{(-21, 3.49) & (-20, 3.5) & (-19, 3.51)},
(-14, 3) ; (-18, 4) **\crv{(-15, 3.49) & (-16, 3.5) & (-17, 3.51)},
(-13, 3) ; (-9, 4) **\crv{(-12, 3.49) & (-11, 3.5) & (-10, 3.51)},
(-5, 3) ; (-9, 4) **\crv{(-6, 3.49) & (-7, 3.5) & (-8, 3.51)},
(5, 3) ; (9, 4) **\crv{(6, 3.49) & (7, 3.5) & (8, 3.51)},
(13, 3) ; (9, 4) **\crv{(12, 3.49) & (11, 3.5) & (10, 3.51)},
(14, 3) ; (18, 4) **\crv{(15, 3.49) & (16, 3.5) & (17, 3.51)},
(22, 3) ; (18, 4) **\crv{(21, 3.49) & (20, 3.5) & (19, 3.51)},
(-22, 4) ; (0, 5) **\crv{(-21, 4.49) & (-11, 4.5) & (-1, 4.51)},
(22, 4) ; (0, 5) **\crv{(21, 4.49) & (11, 4.5) & (1, 4.51)},
(-24, 2) * {1}, (-21, 2) * {0}, (-18, 2) * {1}, (-15, 2) * {2}, (-12, 2) * {0}, (-9, 2) * {1}, (-6, 2) * {2}, (0, 2) * {\cdots},
(  6, 2) * {0}, (  9, 2) * {1}, ( 12, 2) * {2}, ( 15, 2) * {0}, ( 18, 2) * {1}, (21, 2) * {2}, (24, 2) * {0}, (27, 2) * {2}, (0, 6) * {m - 1},
\ar @{-} (-24, 0) ; (1.5,-8.5),
\ar @{-} (-21, 0) ; (-22.5,-0.5),
\ar @{-} (-18, 0) ; (-16.5, -0.5),
\ar @{-} (-15, 0) ; (-19.5, -1.5),
\ar @{-} (-12, 0) ; (-18, -2),
\ar @{-} (-9, 0) ; (-7.5, -0.5),
\ar @{-} (-6, 0) ; (-15, -3),
\ar @{-} (6, 0) ; (-9, -5),
\ar @{-} (9, 0) ; (10.5, -0.5),
\ar @{-} (12, 0) ; (-6, -6),
\ar @{-} (15, 0) ; (-4.5, -6.5),
\ar @{-} (18, 0) ; (19.5, -0.5),
\ar @{-} (21, 0) ; (-1.5, -7.5),
\ar @{-} (24, 0) ; (25.5, -0.5)
\ar @{-} (27, 0) ; (1.5, -8.5)
\end{xy} \quad (m \geq 1).
\end{equation}
The computation of $T_1(\varphi_{2, m})$ is more easy. By definition, 
\begin{equation*}
T_1(\varphi_{2, m}) = \left[T_1 T_{\varepsilon_2}^m T_2^{-1}\left(E_2\right), T_1 \left(E_2\right)\right]_q 
\end{equation*}
and in the computation of $E_{\leq, m\delta - \alpha_1}$ we have already computed
\begin{equation*}
T_1 T_{\varepsilon_2}^m T_2^{-1}\left(E_2\right) = T_{\rho} (T_{\rho^2} T_0 T_1)^{m - 1} (E_2) = (T_{\rho^2} T_1 T_2)^{m - 1} (E_0) = T_{\varepsilon_2}^{m - 1} (E_0).
\end{equation*}
Since $T_{\rho^2} T_1 T_2 (E_0) = \begin{xy}
(-6, 2) * {0}, (-3, 2) * {1}, (0, 2) * {0}, (3, 2) * {2},
\ar @{-} (-6, 0) ; (-1.5, -3),
\ar @{-} (-3, 0) ; (-4.5, -1),
\ar @{-} (0, 0) ; (1.5, -1),
\ar @{-} (3, 0) ; (-1.5, -3)
\end{xy}
$ and $T_{\varepsilon_2}$ fixes $E_1$ and $[E_0, E_2]_q$,
\begin{equation}
T_1(\varphi_{2, m}) = \begin{xy}
(-22, 3) ; (-18, 4) **\crv{(-21, 3.49) & (-20, 3.5) & (-19, 3.51)},
(-14, 3) ; (-18, 4) **\crv{(-15, 3.49) & (-16, 3.5) & (-17, 3.51)},
(-13, 3) ; (-9, 4) **\crv{(-12, 3.49) & (-11, 3.5) & (-10, 3.51)},
(-5, 3) ; (-9, 4) **\crv{(-6, 3.49) & (-7, 3.5) & (-8, 3.51)},
(5, 3) ; (9, 4) **\crv{(6, 3.49) & (7, 3.5) & (8, 3.51)},
(13, 3) ; (9, 4) **\crv{(12, 3.49) & (11, 3.5) & (10, 3.51)},
(14, 3) ; (18, 4) **\crv{(15, 3.49) & (16, 3.5) & (17, 3.51)},
(22, 3) ; (18, 4) **\crv{(21, 3.49) & (20, 3.5) & (19, 3.51)},
(-22, 4) ; (0, 5) **\crv{(-21, 4.49) & (-11, 4.5) & (-1, 4.51)},
(22, 4) ; (0, 5) **\crv{(21, 4.49) & (11, 4.5) & (1, 4.51)},
(-24, 2) * {0}, (-21, 2) * {1}, (-18, 2) * {0}, (-15, 2) * {2}, (-12, 2) * {1}, (-9, 2) * {0}, (-6, 2) * {2}, (0, 2) * {\cdots},
(  6, 2) * {1}, (  9, 2) * {0}, ( 12, 2) * {2}, ( 15, 2) * {1}, ( 18, 2) * {0}, (21, 2) * {2}, (24, 2) * {1}, (27, 2) * {2}, (0, 6) * {m - 1},
\ar @{-} (-24, 0) ; (1.5,-8.5),
\ar @{-} (-21, 0) ; (-22.5,-0.5),
\ar @{-} (-18, 0) ; (-16.5, -0.5),
\ar @{-} (-15, 0) ; (-19.5, -1.5),
\ar @{-} (-12, 0) ; (-18, -2),
\ar @{-} (-9, 0) ; (-7.5, -0.5),
\ar @{-} (-6, 0) ; (-15, -3),
\ar @{-} (6, 0) ; (-9, -5),
\ar @{-} (9, 0) ; (10.5, -0.5),
\ar @{-} (12, 0) ; (-6, -6),
\ar @{-} (15, 0) ; (-4.5, -6.5),
\ar @{-} (18, 0) ; (19.5, -0.5),
\ar @{-} (21, 0) ; (-1.5, -7.5),
\ar @{-} (24, 0) ; (25.5, -0.5)
\ar @{-} (27, 0) ; (1.5, -8.5)
\end{xy} \quad (m \geq 1).
\end{equation}

We also require $T_{\Check{w}} (\varphi_{i, m}) = T_2 T_1 (\varphi_{i, m})$ for the reversed order $\leq'$. Fortunately, in this case we can simply apply $T_2$ on every leaf of the presentation of $T_1(\varphi_{1, m})$. As a result,
\begin{equation} \label{eq:A2_reversed_im_root_vectors}
T_2 T_1(\varphi_{1, m}) = \begin{xy}
(-22, 3) ; (-18, 4) **\crv{(-21, 3.49) & (-20, 3.5) & (-19, 3.51)},
(-14, 3) ; (-18, 4) **\crv{(-15, 3.49) & (-16, 3.5) & (-17, 3.51)},
(-13, 3) ; (-9, 4) **\crv{(-12, 3.49) & (-11, 3.5) & (-10, 3.51)},
(-5, 3) ; (-9, 4) **\crv{(-6, 3.49) & (-7, 3.5) & (-8, 3.51)},
(5, 3) ; (9, 4) **\crv{(6, 3.49) & (7, 3.5) & (8, 3.51)},
(13, 3) ; (9, 4) **\crv{(12, 3.49) & (11, 3.5) & (10, 3.51)},
(14, 3) ; (18, 4) **\crv{(15, 3.49) & (16, 3.5) & (17, 3.51)},
(22, 3) ; (18, 4) **\crv{(21, 3.49) & (20, 3.5) & (19, 3.51)},
(-22, 4) ; (0, 5) **\crv{(-21, 4.49) & (-11, 4.5) & (-1, 4.51)},
(22, 4) ; (0, 5) **\crv{(21, 4.49) & (11, 4.5) & (1, 4.51)},
(-27, 2) * {2}, (-24, 2) * {1}, (-21, 2) * {2}, (-18, 2) * {0}, (-15, 2) * {1}, (-12, 2) * {2}, (-9, 2) * {0}, (-6, 2) * {1}, (0, 2) * {\cdots},
(  6, 2) * {2}, (  9, 2) * {0}, ( 12, 2) * {1}, ( 15, 2) * {2}, ( 18, 2) * {0}, (21, 2) * {1}, (24, 2) * {0}, (0, 6) * {m - 1},
\ar @{-} (-27, 0) ; (-1.5,-8.5),
\ar @{-} (-24, 0) ; (-25.5, -0.5),
\ar @{-} (-21, 0) ; (-19.5, -0.5),
\ar @{-} (-18, 0) ; (-22.5, -1.5),
\ar @{-} (-15, 0) ; (-21, -2),
\ar @{-} (-12, 0) ; (-10.5, -0.5),
\ar @{-} (-9, 0) ; (-18, -3),
\ar @{-} (-6, 0) ; (-16.5, -3.5),
\ar @{-} (6, 0) ; (7.5, -0.5),
\ar @{-} (9, 0) ; (-9, -6),
\ar @{-} (12, 0) ; (-7.5, -6.5),
\ar @{-} (15, 0) ; (16.5, -0.5),
\ar @{-} (18, 0) ; (-4.5, -7.5),
\ar @{-} (21, 0) ; (-3, -8),
\ar @{-} (24, 0) ; (-1.5, -8.5)
\end{xy},
\end{equation}
\begin{equation} \label{eq:A2_reversed_im_root_vectors2}
T_2 T_1(\varphi_{2, m}) = \begin{xy}
(-22, 3) ; (-18, 4) **\crv{(-21, 3.49) & (-20, 3.5) & (-19, 3.51)},
(-14, 3) ; (-18, 4) **\crv{(-15, 3.49) & (-16, 3.5) & (-17, 3.51)},
(-13, 3) ; (-9, 4) **\crv{(-12, 3.49) & (-11, 3.5) & (-10, 3.51)},
(-5, 3) ; (-9, 4) **\crv{(-6, 3.49) & (-7, 3.5) & (-8, 3.51)},
(5, 3) ; (9, 4) **\crv{(6, 3.49) & (7, 3.5) & (8, 3.51)},
(13, 3) ; (9, 4) **\crv{(12, 3.49) & (11, 3.5) & (10, 3.51)},
(14, 3) ; (18, 4) **\crv{(15, 3.49) & (16, 3.5) & (17, 3.51)},
(22, 3) ; (18, 4) **\crv{(21, 3.49) & (20, 3.5) & (19, 3.51)},
(-22, 4) ; (0, 5) **\crv{(-21, 4.49) & (-11, 4.5) & (-1, 4.51)},
(22, 4) ; (0, 5) **\crv{(21, 4.49) & (11, 4.5) & (1, 4.51)},
(-27, 2) * {2}, (-24, 2) * {0}, (-21, 2) * {2}, (-18, 2) * {1}, (-15, 2) * {0}, (-12, 2) * {2}, (-9, 2) * {1}, (-6, 2) * {0}, (0, 2) * {\cdots},
(  6, 2) * {2}, (  9, 2) * {1}, ( 12, 2) * {0}, ( 15, 2) * {2}, ( 18, 2) * {1}, (21, 2) * {0}, (24, 2) * {1}, (0, 6) * {m - 1},
\ar @{-} (-27, 0) ; (-1.5,-8.5),
\ar @{-} (-24, 0) ; (-25.5, -0.5),
\ar @{-} (-21, 0) ; (-19.5, -0.5),
\ar @{-} (-18, 0) ; (-22.5, -1.5),
\ar @{-} (-15, 0) ; (-21, -2),
\ar @{-} (-12, 0) ; (-10.5, -0.5),
\ar @{-} (-9, 0) ; (-18, -3),
\ar @{-} (-6, 0) ; (-16.5, -3.5),
\ar @{-} (6, 0) ; (7.5, -0.5),
\ar @{-} (9, 0) ; (-9, -6),
\ar @{-} (12, 0) ; (-7.5, -6.5),
\ar @{-} (15, 0) ; (16.5, -0.5),
\ar @{-} (18, 0) ; (-4.5, -7.5),
\ar @{-} (21, 0) ; (-3, -8),
\ar @{-} (24, 0) ; (-1.5, -8.5)
\end{xy} \quad (m \geq 1).
\end{equation}

Set the projection $\pi_B^+: U_q^+ \rightarrow \mathcal{P}_B^+$ of section \ref{sec:representation} by $\sigma_{01} := +1, \sigma_{02} := +1, \sigma_{12} := +1$. Then corresponding skew-symmetric matrix is $B = \bigl(\begin{smallmatrix} 0 & -1 & -1 \\ 1 & 0 & -1 \\ 1 & 1 & 0\end{smallmatrix}\bigr)$, which corresponds to the Dynkin quiver 
\begin{xy}
\ar (5, 8.66) *+[Fo]{0}="A"; (0, 0) *+[Fo]{1}="B",
\ar "A" ; (10, 0) *+[Fo]{2}="C",
\ar "B" ; "C"
\end{xy}
. The matrix presentation of the bilinear form $\langle \cdot, \cdot\rangle_B$ is $(\langle \alpha_i, \alpha_j \rangle_B)_{i, j = 0}^{2} = \bigl(\begin{smallmatrix} 2 & 0 & 0 \\ -2 & 2 & 0 \\ -2 & -2 & 2\end{smallmatrix}\bigr)$. Thus $[E_0, E_1]_q$, $[E_0, E_2]_q$, and $[E_1, E_2]_q$ lie in the kernel of $\pi_B^+$. Examining the computed presentations of root vectors above, all the root vectors except for simple root vectors vanish by $\pi_B^+$. Therefore, the image of quasi-universal R-matrix $\Theta$ is
\begin{equation} \label{eq_A2_image_theta}
\pi_B^+ \widehat{\otimes} \pi_B^- (\Theta) = \mathbb{E}(y_2) \mathbb{E}(y_1) \mathbb{E}(y_0) \quad \in \widehat{\mathcal{D}_B}.
\end{equation}

Now we reverse the given order. Recall that the real root vectors for the reversed order $\leq'$ are obtained just by reversing all the direction of $q$-bracket. Then one can verify that the real root vectors in the first row, namely, $E_{\leq', m\delta - \alpha_1} = \Psi (E_{\leq, m\delta - \alpha_1})$,  $E_{\leq', m\delta - \alpha_1 - \alpha_2} = \Psi (E_{\leq, m\delta - \alpha_1 - \alpha_2})$ $(m \geq 1)$ satisfy the condition of Proposition \ref{prop:image_real_root_vec}.

However, by contrast, real root vectors in second row behave differently. $E_{\leq', \delta - \alpha_1} = \Psi (E_{\leq, \delta - \alpha_1}) = [E_2, E_0]_q$ does not vanish by $\pi_B^+$ since $\langle \alpha_2, \alpha_0\rangle_B = -2 \neq 0$ and satisfies the condition of Proposition \ref{prop:image_real_root_vec}. In the same way $E_{\leq', \alpha_1} = E_1$ does not vanish. But $E_{\leq', m\delta - \alpha_1}$ and $E_{\leq', (m - 1)\delta + \alpha_1}$ vanish for $m > 1$ because $\langle \alpha_0 + \alpha_2, \alpha_1 + \alpha_2\rangle_B = \langle \alpha_0, \alpha_1\rangle_B = 0$. Thus the real root vectors in the second row vanish except for $E_{\leq', \delta - \alpha_1}$ and $E_{\leq', \alpha_1}$.

Next, we have to compute the image of imaginary root vectors $T_2 T_1 (I_{i, m})$. Using (\ref{eq:weight_vec_comm}) and the presentations (\ref{eq:A2_reversed_im_root_vectors}) (\ref{eq:A2_reversed_im_root_vectors2}), we have
\begin{align*}
\pi_B^+ T_2 T_1(\varphi_{1, m}) &= (1 - q^{-2})^{3m - 2} (1 - q^{-2(m + 1)}) e_2 e_1 (e_2 e_0 e_1)^{m - 1} e_0 \\
&= (q - q^{-1})^{3m - 1} [m + 1]_q (e_0 e_1 e_2)^m, \\
\pi_B^+ T_2 T_1(\varphi_{2, m}) &= \begin{cases}
(q - q^{-1}) e_0 e_2 & m = 1 \\
0 & m > 1
\end{cases},
\end{align*}
since $\langle \alpha_2 + \alpha_0, \alpha_2 + \alpha_1\rangle_B = 0$. The second equality is due to the following calculation. Recall that the commutation relations in $\mathcal{P}_B^+$ become $e_1 e_0 = q e_0 e_1$, $e_2 e_0 = q e_0 e_2$, $e_2 e_1 = q e_1 e_2$. Thus
\begin{equation*}
e_2 e_1 (e_2 e_0 e_1)^{m - 1} e_0 = q e_1 e_2 (q^2 e_0 e_1 e_2)^{m - 1} e_0 = q^{4m - 1} (e_0 e_1 e_2)^m \quad (m \geq 1).
\end{equation*}

Let $D := (q - q^{-1})^3 e_0 e_1 e_2$. Then by definition of generating function $\varphi_i(z) \in U_q^+[[z]]$,
\begin{align*}
\pi_B^+ (T_2 T_1 \varphi_1(z)) &= \sum_{m = 1}^{\infty} [m + 1]_q D^m z^m = \frac{1}{(1 - qDz) (1 - q^{-1}Dz)} \quad \in \mathcal{P}_B^+ [[z]], \\
\pi_B^+ (T_2 T_1 \varphi_2(z)) &= Dz.
\end{align*}
Thus the image of imaginary root vectors are
\begin{align}
\pi_B^+ (T_2 T_1 I_{1, m}) &= \frac{q^m + q^{-m}}{m(q - q^{-1})} D^m, \\
\pi_B^+ (T_2 T_1 I_{2, m}) &= \frac{(-1)^{m - 1}}{m(q - q^{-1})} D^m \quad (m \geq 1).
\end{align}

At last, we compute the image of $S'_m := (T_2 T_1 \otimes T_2 T_1) (S_m)$ $(m \geq 1)$ in (\ref{eq:def_Sn}). By definition,
\begin{align*}
\begin{bmatrix} b_{1, 1; m} & b_{1, 2; m} \\ b_{2, 1; m} & b_{2, 2; m} \end{bmatrix} &= \frac{1}{m(q^{-1} - q)} \begin{bmatrix} [2m]_q & (-1)^{m - 1}[m]_q \\ (-1)^{m - 1}[m]_q & [2m]_q \end{bmatrix}, \\
\begin{bmatrix} c_{1, 1; m} & c_{1, 2; m} \\ c_{2, 1; m} & c_{2, 2; m} \end{bmatrix} &= \frac{m(q^{-1} - q)}{[2m]_q^2 - [m]_q^2} \begin{bmatrix} [2m]_q & (-1)^{m}[m]_q \\ (-1)^{m}[m]_q & [2m]_q \end{bmatrix}
\end{align*}
and thus $S_m$ is written down as
\begin{multline}
\label{eq:A_2(1)_S_m}
S_m = \frac{m(q^{-1} - q)}{[2m]_q^2 - [m]_q^2} \Bigl\{[2m]_q(I_{1, m} \otimes J_{1, m} + I_{2, m} \otimes J_{2, m}) \\
+ (-1)^m [m]_q (I_{1, m} \otimes J_{2, m} + I_{2, m} \otimes J_{1, m})\Bigr\}.
\end{multline}
Let $D' := \overline{\Omega_B} (D) = (q^{-1} - q)^2 f_2 f_1 f_0$. Then
\begin{align}
\pi_B^- (T_2 T_1 J_{1, m}) &= \frac{q^m + q^{-m}}{m(q^{-1} - q)} D'^m, \\
\pi_B^- (T_2 T_1 J_{2, m}) &= \frac{(-1)^{m - 1}}{m(q^{-1} - q)} D'^m \quad (m \geq 1).
\end{align}
Therefore the image of $S'_m$ is computed as follows.
\begin{align*}
&(\pi_B^+ \otimes \pi_B^-) (S'_m) \\
&= \frac{m(q^{-1} - q)}{[2m]_q^2 - [m]_q^2} \Biggl\{[2m]_q\left(-\frac{(q^m - q^{-m})^2}{m^2(q - q^{-1})^2} - \frac{1}{m^2(q - q^{-1})^2}\right) D^m \otimes {D'}^m \\
&\quad  + (-1)^m [m]_q \cdot 2 \frac{(-1)^m(q^m + q^{-m})}{m^2(q - q^{-1})^2} D^m \otimes {D'}^m \Biggr\} \\
&= \frac{m(q^{-1} - q)}{[2m]_q^2 - [m]_q^2} \cdot \frac{-[2m]_q \left(1 + (q^m + q^{-m})^2\right) + 2[m]_q (q^m + q^{-m})}{m^2(q - q^{-1})^2} (D \otimes D')^m \\
&= \frac{[2m]_q(q^{2m} + 3 + q^{-2m}) - 2[m]_q (q^m + q^{-m})}{\left([2m]_q^2 - [m]_q^2\right) \cdot m(q - q^{-1})} (D \otimes D')^m \\
&= \frac{1}{m} \frac{(q^{2m} - q^{-2m}) (q^{2m} + 3 + q^{-2m}) - 2(q^m - q^{-m}) (q^m + q^{-m})}{(q^{2m} - q^{-2m})^2 - (q^m - q^{-m})^2} (D \otimes D')^m \\
&= \frac{1}{m} \frac{(q^m - q^{-m}) (q^m + q^{-m}) (q^{2m} + 1 + q^{-2m})}{(q^m - q^{-m})^2 (q^{2m} + 1 + q^{-2m})} (D \otimes D')^m \\
&= -\frac{1}{m} \frac{q^m (q^m + q^{-m})}{1 - q^{2m}} (D \otimes D')^m.
\end{align*}
This result coincides with the case of type $A_1^{(1)}$. Therefore
\begin{equation}
\begin{aligned}
\pi_B^+ \widehat{\otimes} \pi_B^- (\Theta_{\mathrm{im}}) &= \mathbb{E} (-q D\otimes D')^{-1} \mathbb{E} (-q^{-1} D\otimes D')^{-1} \\
&= \mathbb{E} (-q \normord{y_0 y_1 y_2})^{-1} \mathbb{E} (-q^{-1} \normord{y_0 y_1 y_2})^{-1}.
\end{aligned}
\end{equation}
Comparing (\ref{eq_A2_image_theta}), finally we attain the following identity.
\begin{thm}
Let $y_0, y_1, y_2$ be indeterminate with commutation relation $y_0 y_1 = q^{-2} y_1 y_0$, $y_0 y_2 = q^{-2} y_2 y_0$, $y_1 y_2 = q^{-2} y_2 y_1$. Then the following identity holds in skew formal power series algebra $\widehat{\mathcal{S}_2} := \Q(q)[[y_0, y_1, y_2]]$.
\begin{equation} \label{eq:A2_result}
\begin{aligned}
&\mathbb{E}(\normord{y_2}) \mathbb{E}(\normord{y_1}) \mathbb{E}(\normord{y_0}) \\
&= \quad \left\{\prod_{m \geq 0}^{\rightarrow} \mathbb{E}(\normord{y_0^{m + 1} y_1^m y_2^m}) \mathbb{E}(\normord{y_0^{m + 1} y_1^{m + 1} y_2^m})\right\} \mathbb{E}(\normord{y_0 y_2}) \\
&\quad \qquad \times \mathbb{E} (-q \normord{y_0 y_1 y_2})^{-1} \mathbb{E} (-q^{-1} \normord{y_0 y_1 y_2})^{-1} \\
&\quad \qquad \qquad \times \mathbb{E}(\normord{y_1}) \left\{\prod_{m \geq 0}^{\leftarrow} \mathbb{E}(\normord{y_0^m y_1^{m + 1} y_2^{m + 1}}) \mathbb{E}(\normord{y_0^m y_1^m y_2^{m + 1}})\right\}.
\end{aligned}
\end{equation}
where $\prod_{m \geq 0}^{\rightarrow} a_m := a_0 a_1 a_2 \dots$, $\prod_{m \geq 0}^{\leftarrow} a_m := \dots a_2 a_1 a_0$, and normal ordered product is $\normord{y_0^{m_0} y_1^{m_1} y_2^{m_2}} = q^{m_0 m_1 + m_0 m_2 + m_1 m_2} y_0^{m_0} y_1^{m_1} y_2^{m_2}$.
\end{thm}

Set $R := \bigl(\begin{smallmatrix} 0 & 1 & 1 \\ 1 & 2 & 1 \end{smallmatrix}\bigr)$. Then the condition of Proposition \ref{prop:S_B_hom} holds for $B' = \bigl(\begin{smallmatrix} 0 & 1 \\ -1 & 0 \end{smallmatrix}\bigr)$. Thus we have continuous algebra homomorphism $\widehat{\psi_2}: \widehat{\mathcal{S}_2} \rightarrow \widehat{\mathcal{S}}$ which satisfies 
\begin{equation}
\widehat{\psi_2}(y_0) = x_2, \quad \widehat{\psi_2}(y_1) = \normord{x_1 x_2^2}, \quad \widehat{\psi_2}(y_2) = \normord{x_1 x_2}.
\end{equation}
Applying $\widehat{\psi_2}$ on (\ref{eq:A2_result}) and transforming variables using $\psi_S$, we obtain (\ref{eq:A2}). This proves that the identity (\ref{eq:A2}) holds in $\Q(q)[[\frac{x_1}{x_2^2}, x_2]]$.

\textbf{Remark.} To derive (\ref{eq:A2}), convex order of multiple row is mandatory because the factor ${\bf U}_{1, 0}^2$ in the middle cannot be appear by only using convex orders of single row. We also note that convex orders of multiple row never appear as the form of $\leq_Z$ determined by central charge $Z$ (\ref{eq:order_by_central_charge}).

\subsection{Type $A_3^{(1)}$}

Let $\mathfrak{g} = \widehat{\mathfrak{sl}_4}$ be affine algebra of type $A_3^{(1)}$. We set a convex order by
\begin{align*}
w &:= s_2 s_1 s_3, \quad n := 3; \\
\mathring{I} &= J_0 \supsetneqq J_1 := \{1, 3\} \supsetneqq J_2 := \{3\} \supsetneqq J_3 = \emptyset, \\
{\bf s}_0 &:= (s_0 s_2 s_1 s_3 s_2 s_0 s_3 s_1)^{\infty}, \quad {\bf s}_1 := (s_1 s_{\delta - \alpha_1})^{\infty}, \\
{\bf s}_2 &:= (s_3 s_{\delta - \alpha_3})^{\infty}, \\
& \\
\Check{w} &= ww_{\circ} = s_1 s_3 s_2; \\
\mathring{I} &= \Check{J}_0 \supsetneqq \Check{J}_1 := \{1, 3\} \supsetneqq \Check{J}_2 := \{1\} \supsetneqq \Check{J}_3 = \emptyset, \\
\Check{\bf s}_0 &:= (s_1 s_3 s_2 s_0 s_3 s_1 s_0 s_2)^{\infty}, \quad \Check{\bf s}_1 := (s_{\delta - \alpha_3} s_3)^{\infty}, \\
\Check{\bf s}_2 &:= (s_{\delta - \alpha_1} s_1)^{\infty}.
\end{align*}
Then the corresponding convex order $\leq$ is as follows.
\begin{equation*}
\begin{aligned}
\delta - \alpha_1 - \alpha_2 - \alpha_3 < \alpha_2 < &\delta - \alpha_3 < \delta - \alpha_1 \\
\quad  < 2\delta - \alpha_1 - \alpha_2 - \alpha_3 < &\delta + \alpha_2 < 2\delta - \alpha_3 < 2\delta - \alpha_1 \\
&\dots \\
< \alpha_1 + \alpha_2 < \delta + \alpha_1 + &\alpha_2 < 2\delta + \alpha_1 + \alpha_2 < 3\delta + \alpha_1 + \alpha_2 < \dots \\
< \alpha_2 + \alpha_3 < \delta + \alpha_2 + &\alpha_3 < 2\delta + \alpha_2 + \alpha_3 < 3\delta + \alpha_2 + \alpha_3 < \dots \\
< \delta < 2\delta < &3\delta < 4\delta < \dots \\
\dots < 3\delta - \alpha_2 - \alpha_3 &< 2\delta - \alpha_2 - \alpha_3 < \delta - \alpha_2 - \alpha_3 \\
\dots < 3\delta - \alpha_1 - \alpha_2 &< 2\delta - \alpha_1 - \alpha_2 < \delta - \alpha_1 - \alpha_2 \\
&\dots \\
\quad &< 2\delta - \alpha_2 < \delta + \alpha_1 + \alpha_2 + \alpha_3 < \delta + \alpha_3 < \delta + \alpha_1 \\
\quad &\quad < \delta - \alpha_2 < \alpha_1 + \alpha_2 + \alpha_3 < \alpha_3 < \alpha_1.
\end{aligned}
\end{equation*}
where the null root $\delta = \alpha_0 + \alpha_1 + \alpha_2 + \alpha_3$.

Set
\begin{equation}
B := \begin{pmatrix}
0 & -1 & 0 & -1 \\
1 &  0 &  1 &  0 \\
0 & -1 & 0 & -1 \\
1 &  0 &  1 &  0
\end{pmatrix} \quad (\text{corresponding quiver:}
\begin{xy}
(0, 5) *+[Fo]{0}="0", (0, -5) *+[Fo]{1}="1", (10, -5) *+[Fo]{2}="2", (10, 5) *+[Fo]{3}="3",
\ar "0" ; "1",
\ar "0" ; "3",
\ar "2" ; "1",
\ar "2" ; "3"
\end{xy}
).
\end{equation}
The matrix presentation of the bilinear form $\langle \cdot, \cdot\rangle_B$ is 
\begin{equation*}
(\langle \alpha_i, \alpha_j \rangle_B)_{i, j = 0}^{3} = \begin{pmatrix}
 2 & 0 & 0 & 0 \\
 -2 & 2 & -2 & 0 \\
 0 & 0 & 2 & 0 \\
 -2 & 0 & -2 & 2
\end{pmatrix}.
\end{equation*}
In the same way in the examples so far, one can verify that all the root vectors except for simple root vectors vanish by the projection $\pi_B^+$. On the other hand, real root vectors for the reversed order $\leq'$, namely, $E_{\leq', m\delta + \alpha}$ for $\alpha = \pm \alpha_1, \pm \alpha_2, \pm \alpha_3, \pm (\alpha_1 + \alpha_2 + \alpha_3)$ ($m \geq 0$ if $\alpha \in \Delta_+$, $m \geq 1$ if $\alpha \in \Delta_-$) do not vanish by $\pi_B^+$ and satisfy the condition of Proposition \ref{prop:image_real_root_vec}. The real root vectors in the second and third row $E_{\leq', m\delta + \alpha_1 + \alpha_2}$, $E_{\leq', m\delta + \alpha_2 + \alpha_3}$, $E_{\leq', (m + 1) \delta - \alpha_1 - \alpha_2}$, $E_{\leq', (m + 1) \delta - \alpha_2 - \alpha_3}$ vanish if and only if $m > 0$, and satisfy the condition of Proposition \ref{prop:image_real_root_vec} when $m = 0$. This behavior in the second and third row resembles that of the case of type $A_2^{(1)}$.

Although the computation of imaginary root vectors also has resemblance to the previous examples and in fact we will obtain identical presentation of $\pi_B^+ \widehat{\otimes} \pi_B^- (\Theta_{\mathrm{im}})$, the process of computation is far from obvious. After somewhat lengthy computation (we used $T_{\varepsilon_2}^{-1} (E_i) = E_i$ for $i = 1, 3$ in the process), one will obtain
\begin{align*}
T_{\Check{w}} (\varphi_{1, m}) &= [(T_1 T_{\rho} T_3 T_2)^{m - 1} T_1 (E_0), [E_3, E_2]_q]_q, \\
T_{\Check{w}} (\varphi_{2, m}) &= [(T_1 T_3 T_2T_{\rho^2} T_2)^{m - 1} T_3 T_1 (E_2), E_0]_q, \\
T_{\Check{w}} (\varphi_{3, m}) &= [(T_3 T_{\rho^3} T_1 T_2)^{m - 1} T_3 (E_0), [E_1, E_2]_q]_q \quad (m \geq 1),
\end{align*}
\begin{align*}
\pi_B^+ T_{\Check{w}} (\varphi_{1, m}) &= \pi_B^+ T_{\Check{w}} (\varphi_{3, m}) = \begin{cases}
(q - q^{-1}) e_0 e_2 e_1 e_3 & m = 1 \\
0 & m > 1
\end{cases}, \\
\pi_B^+ T_{\Check{w}} (\varphi_{2, m}) &= (q - q^{-1})^{4m - 1} [m + 1]_q (e_0 e_2 e_1 e_3)^m \quad (m \geq 1).
\end{align*}
Thus the images of imaginary root vectors are
\begin{align*}
\pi_B^+ T_{\Check{w}} (I_{1, m}) &= \pi_B^+ T_{\Check{w}} (I_{3, m}) = \frac{(-1)^{m - 1}}{m(q - q^{-1})} D^m,\\
\pi_B^+ T_{\Check{w}} (I_{2, m}) &=\frac{q^m + q^{-m}}{m(q - q^{-1})} D^m \quad (m \geq 1),
\end{align*}
where $D := (q - q^{-1})^4 e_0 e_2 e_1 e_3$.

Our last task is to compute the image of $S'_m := (T_{\Check{w}} \otimes T_{\Check{w}}) (S_m)$ $(m \geq 1)$. The definition (\ref{eq:bijn}) reads
\begin{align*}
(b_{i, j; m})_{i, j = 1}^3 &= \frac{1}{m(q^{-1} - q)} \begin{bmatrix} M_2 & (-1)^{m - 1}M_1 & 0 \\ (-1)^{m - 1}M_1 & M_2 & (-1)^{m - 1}M_1 \\ 0 & (-1)^{m - 1}M_1 & M_2\end{bmatrix}, \\
(c_{i, j; m})_{i, j = 1}^3 &= \frac{m(q^{-1} - q)}{M_2^3 - 2M_1^2 M_2} \begin{bmatrix} M_2^2 - M_1^2 & (-1)^{m}M_1 M_2 & M_1^2 \\ (-1)^{m}M_1 M_2 & M_2^2 & (-1)^{m}M_1 M_2 \\ M_1^2 & (-1)^{m}M_1 M_2 & M_2^2 - M_1^2 \end{bmatrix},
\end{align*}
where $M_k := [km]_q$ for $k = 1, 2$. Thus the image of $S'_m$ is computed as follows.
\begin{align*}
(\pi_B^+ \otimes \pi_B^-) (S'_m) &= \frac{m(q^{-1} - q)}{[2m]_q^3 - 2[2m]_q[m]_q^2}  \\
&\quad \times \Bigl(2([2m]_q^2 - [m]_q^2) \frac{1}{m (q - q^{-1}) m (q ^{-1} - q))} \\
&\quad \quad + 4(-1)^m [2m]_q [m]_q \frac{(-1)^{m - 1}(q^m + q^{-m})}{m (q - q^{-1}) m (q ^{-1} - q))} \\
&\quad \quad + 2[m]_q^2 \frac{1}{m (q - q^{-1}) m (q ^{-1} - q))} \\
&\quad \quad + [2m]_q^2 \frac{(q^m + q^{-m})^2}{m (q - q^{-1}) m (q ^{-1} - q))}\Bigr) \\
&\quad \times D^m \otimes D'^m \\
&= -\frac{1}{m} \frac{q^m (q^m + q^{-m})}{1 - q^{2m}} (D \otimes D')^m,
\end{align*}
where $D' := \overline{\Omega_B} (D)$. This is identical with the previous examples and we conclude $\pi_B^+ \widehat{\otimes} \pi_B^- (\Theta_{\mathrm{im}}) = \mathbb{E} (-q D\otimes D')^{-1} \mathbb{E} (-q^{-1} D\otimes D')^{-1}$. As a result,
\begin{thm}
The following identity holds in skew formal power series algebra $\widehat{\mathcal{S}_3} := \widehat{\mathcal{S}_B} \cong \Q(q)[[y_0, y_1, y_2, y_3]]$.
\begin{equation} \label{eq:A3_result}
\begin{aligned}
& \mathbb{E}(\normord{y_1}) \mathbb{E}(\normord{y_3}) \mathbb{E}(\normord{y_2}) \mathbb{E}(\normord{y_0}) \\
&= \left\{\prod_{m \geq 0}^{\rightarrow} X_m \right\} \mathbb{E}(\normord{y_1 y_2}) \mathbb{E}(\normord{y_2 y_3}) \\
& \times \mathbb{E} (-q \normord{y_0 y_1 y_2 y_3})^{-1} \mathbb{E} (-q^{-1} \normord{y_0 y_1 y_2 y_3})^{-1} \\
& \times \mathbb{E}(\normord{y_0 y_1}) \mathbb{E}(\normord{y_0 y_3}) \left\{\prod_{m \geq 0}^{\leftarrow} Y_m \right\},
\end{aligned}
\end{equation}
where
\begin{align*}
X_m &= \mathbb{E}(\normord{y_0^{m + 1} y_1^{m} y_2^{m} y_3^{m}}) \mathbb{E}(\normord{y_0^{m} y_1^{m} y_2^{m + 1} y_3^{m}}) \\
&\quad \times \mathbb{E}(\normord{y_0^{m + 1} y_1^{m + 1} y_2^{m + 1} y_3^{m}}) \mathbb{E}(\normord{y_0^{m + 1} y_1^{m} y_2^{m + 1} y_3^{m + 1}}), \\
Y_m &= \mathbb{E}(\normord{y_0^{m + 1} y_1^{m + 1} y_2^{m} y_3^{m + 1}}) \mathbb{E}(\normord{y_0^{m} y_1^{m + 1} y_2^{m + 1} y_3^{m + 1}}) \\
&\quad \times \mathbb{E}(\normord{y_0^{m} y_1^{m} y_2^{m} y_3^{m + 1}}) \mathbb{E}(\normord{y_0^{m} y_1^{m + 1} y_2^{m} y_3^{m}}).
\end{align*}
\end{thm}
Set $R := \bigl(\begin{smallmatrix} 0 & 1 & 0 & 1 \\ 1 & 1 & 1 & 1 \end{smallmatrix}\bigr)$. Then the condition of Proposition \ref{prop:S_B_hom} holds for the same $B'$. Thus we have continuous algebra homomorphism $\widehat{\psi_3}: \widehat{\mathcal{S}_3} \rightarrow \widehat{\mathcal{S}}$ which satisfies 
\begin{equation}
\widehat{\psi_3}(y_0) = \widehat{\psi_3}(y_2) = x_2, \quad \widehat{\psi_3}(y_1) = \widehat{\psi_3}(y_3) = \normord{x_1 x_2}.
\end{equation}
Applying $\widehat{\psi_3}$ on (\ref{eq:A3_result}) and transforming $x_1, x_2$ by $\psi_S$ yields (\ref{eq:A3}). This proves the identity (\ref{eq:A3}).

\subsection{Type $D_4^{(1)}$}

Let $\mathfrak{g} = \widehat{\mathfrak{so}_8}$ be affine algebra of type $D_4^{(1)}$. Let
\begin{equation}
B := \begin{pmatrix}
0 & 0 & -1 & 0 & 0 \\
0 & 0 & -1 & 0 & 0 \\
1 & 1 &  0 & 1 & 1 \\
0 & 0 & -1 & 0 & 0 \\
0 & 0 & -1 & 0 & 0
\end{pmatrix} \quad (\text{corresponding quiver:}
\begin{xy}
(0, 8) *+[Fo]{0}="0", (0, -8) *+[Fo]{1}="1", (16, -8) *+[Fo]{4}="4", (16, 8) *+[Fo]{3}="3", (8, 0) *+[Fo]{2}="2",
\ar "0" ; "2",
\ar "1" ; "2",
\ar "3" ; "2",
\ar "4" ; "2"
\end{xy}
).
\end{equation}
Using the assignment of indices in the Dynkin quiver above, we set a convex order by
\begin{align*}
w &:= s_1 s_3 s_4 s_2 s_1 s_3 s_4, \quad n := 4; \\
\mathring{I} &= J_0 \supsetneqq J_1 := \{1, 3, 4\} \supsetneqq J_2 := \{3, 4\} \supsetneqq J_3 := \{4\}  \supsetneqq J_4 = \emptyset, \\
{\bf s}_0 &:= (s_0 s_1 s_3 s_4 s_2)^{\infty}, \quad {\bf s}_1 := (s_1 s_{\delta - \alpha_1})^{\infty}, \\
{\bf s}_2 &:= (s_3 s_{\delta - \alpha_3})^{\infty}, \quad {\bf s}_3 := (s_4 s_{\delta - \alpha_4})^{\infty}, \\
& \\
\Check{w} &= ww_{\circ} = s_2 s_1 s_3 s_4 s_2; \\
\mathring{I} &= \Check{J}_0 \supsetneqq \Check{J}_1 := \{1, 3, 4\} \supsetneqq \Check{J}_2 := \{3, 4\} \supsetneqq \Check{J}_3 := \{3\} \supsetneqq \Check{J}_4 = \emptyset, \\
\Check{\bf s}_0 &:= (s_2 s_1 s_3 s_4 s_0)^{\infty}, \quad \Check{\bf s}_1 := (s_{\delta - \alpha_1} s_1)^{\infty}, \\
\Check{\bf s}_2 &:= (s_{\delta - \alpha_3} s_3)^{\infty}, \quad \Check{\bf s}_3 := (s_{\delta - \alpha_4} s_4)^{\infty}.
\end{align*}
Then the corresponding convex order $\leq$ is as follows.
\begin{equation*}
\begin{aligned}
&\delta - \alpha_1 - 2\alpha_2 - \alpha_3 - \alpha_4 < \alpha_1 < \alpha_3 < \alpha_4 \\
&< \delta - \alpha_2 < \alpha_1 + \alpha_2 + \alpha_3 + \alpha_4 < \delta - \alpha_1 - \alpha_2 \\
&< \delta - \alpha_2 - \alpha_3 < \delta - \alpha_2 - \alpha_4 < 2\delta - \alpha_2 \\
&\quad 2\delta - \alpha_1 - 2\alpha_2 - \alpha_3 - \alpha_4 < \delta + \alpha_1 < \delta + \alpha_3 < \delta + \alpha_4 \\
&\quad < 3\delta - \alpha_2 < \delta + \alpha_1 + \alpha_2 + \alpha_3 + \alpha_4 < 2\delta - \alpha_1 - \alpha_2 \\
&\quad < 2\delta - \alpha_2 - \alpha_3 < 2\delta - \alpha_2 - \alpha_4 < 4\delta - \alpha_2 \\
&\quad \dots \\
&\quad < \alpha_2 + \alpha_3 + \alpha_4 < \delta + \alpha_2 + \alpha_3 + \alpha_4 < 2\delta + \alpha_2 + \alpha_3 + \alpha_4 < \dots \\
&\quad < \alpha_1 + \alpha_2 + \alpha_4 < \delta + \alpha_1 + \alpha_2 + \alpha_4 < 2\delta + \alpha_1 + \alpha_2 + \alpha_4 < \dots \\
&\quad < \alpha_1 + \alpha_2 + \alpha_3 < \delta + \alpha_1 + \alpha_2 + \alpha_3 < 2\delta + \alpha_1 + \alpha_2 + \alpha_3 < \dots \\
&\quad \quad < \delta < 2\delta < 3\delta < 4\delta < \dots \\
&\dots < 3\delta - \alpha_1 - \alpha_2 - \alpha_3 < 2\delta - \alpha_1 - \alpha_2 - \alpha_3 < \delta - \alpha_1 - \alpha_2 - \alpha_3 \\
&\dots < 3\delta - \alpha_1 - \alpha_2 - \alpha_4 < 2\delta - \alpha_1 - \alpha_2 - \alpha_4 < \delta - \alpha_1 - \alpha_2 - \alpha_4 \\
&\dots < 3\delta - \alpha_2 - \alpha_3 - \alpha_4 < 2\delta - \alpha_2 - \alpha_3 - \alpha_4 < \delta - \alpha_2 - \alpha_3 - \alpha_4 \\
&\quad \quad \dots \\
&\quad \quad < \delta + \alpha_1 + 2\alpha_2 + \alpha_3 + \alpha_4 < 2\delta - \alpha_4 < 2\delta - \alpha_3 < 2\delta - \alpha_1 \\
&\quad \quad < 3\delta + \alpha_2 < 2\delta - \alpha_1 - \alpha_2 - \alpha_3 - \alpha_4 < \delta + \alpha_2 + \alpha_4 \\
&\quad \quad < \delta + \alpha_2 + \alpha_3 < \delta + \alpha_1 + \alpha_2 < 2\delta + \alpha_2 \\
&\quad \quad \quad < \alpha_1 + 2\alpha_2 + \alpha_3 + \alpha_4 < \delta - \alpha_4 < \delta - \alpha_3 < \delta - \alpha_1 \\
&\quad \quad \quad < \delta + \alpha_2 < \delta - \alpha_1 - \alpha_2 - \alpha_3 - \alpha_4 < \alpha_2 + \alpha_4 \\
&\quad \quad \quad < \alpha_2 + \alpha_3 < \alpha_1 + \alpha_2 < \alpha_2,
\end{aligned}
\end{equation*}
where the null root $\delta = \alpha_0 + \alpha_1 + 2\alpha_2 + \alpha_3 + \alpha_4$.

To compute root vectors, reduced expressions of fundamental translations $t_{\varepsilon_i} \in \widehat{W}$ $(i = 1, 2, 3, 4)$ are required. Let $\tau := \bigl(\begin{smallmatrix}
0 & 1 & 2 & 3 & 4 \\
1 & 0 & 2 & 4 & 3
\end{smallmatrix}\bigr)$, $\tau' := \bigl(\begin{smallmatrix}
0 & 1 & 2 & 3 & 4 \\
3 & 4 & 2 & 0 & 1
\end{smallmatrix}\bigr)$ be Dynkin automorphisms. In the same way of the proof of Proposition \ref{prop:translation_reduced_rep}, one can verify that
\begin{equation}
\begin{aligned}
t_{\varepsilon_1} &= \tau s_1 s_2 s_3 s_4 s_2 s_1, \\
t_{\varepsilon_2} &= s_0 s_2 s_3 s_4 s_2 s_1 s_2 s_3 s_4 s_2, \\
t_{\varepsilon_3} &= \tau' s_3 s_2 s_1 s_4 s_2 s_3, \\
t_{\varepsilon_4} &= \tau \tau' s_4 s_2 s_1 s_3 s_2 s_4
\end{aligned}
\end{equation}
are reduced expressions in $\widehat{W}$. Using these formulas, one can compute all the root vectors and verify that all of them except for simple root vectors vanish by $\pi_B^+$. For the reversed order $\leq'$, real root vectors in the first row satisfy the condition of Proposition \ref{prop:image_real_root_vec}. On the other hand, real root vectors in the second, third, and fourth row, namely,
\begin{align*}
&E_{\leq', m\delta + \alpha_1 + \alpha_2 + \alpha_3}, E_{\leq', m\delta + \alpha_1 + \alpha_2 + \alpha_4}, E_{\leq', m\delta + \alpha_2 + \alpha_3 + \alpha_4}, \\
&E_{\leq', (m + 1) \delta - \alpha_1 - \alpha_2 - \alpha_3}, E_{\leq', (m + 1) \delta - \alpha_1 - \alpha_2 - \alpha_4}, E_{\leq', (m + 1) \delta - \alpha_2 - \alpha_3 - \alpha_4} 
\end{align*}
vanish for $m > 0$, and satisfy the condition of Proposition \ref{prop:image_real_root_vec} when $m = 0$.

The image of $T_{\Check{w}} (\varphi_{i, m})$ are computed as follows.
\begin{equation}
\begin{aligned}
\pi_B^+ T_{\Check{w}} (\varphi_{1, m}) &{=} \pi_B^+ T_{\Check{w}} (\varphi_{3, m}) {=} \pi_B^+ T_{\Check{w}} (\varphi_{4, m}) = \begin{cases}
q (q - q^{-1})^5 e_0 e_1 e_3 e_4 e_2^2 & m = 1 \\
0 & m > 1
\end{cases}, \\
\pi_B^+ T_{\Check{w}} (\varphi_{2, m}) &= q^m (q - q^{-1})^{6m - 1} [m + 1]_q (e_0 e_1 e_3 e_4 e_2^2)^m.
\end{aligned}
\end{equation}
Note that $e_i e _j = e_j e_i$ and $e_2 e_i = q e_i e_2$ in $\mathcal{P}_B^+$ for $i, j = 0, 1, 3, 4$. Thus
\begin{equation}
\begin{aligned}
\pi_B^+ T_{\Check{w}} (I_{i, m}) &= \frac{(-1)^{m - 1}}{m(q - q^{-1})} (qD)^m \quad (i = 0, 1, 3, 4), \\
\pi_B^+ T_{\Check{w}} (I_{2, m}) &= \frac{q^m + q^{-m}}{m(q - q^{-1})} (qD)^m \quad (m \geq 1),
\end{aligned}
\end{equation}
where $D := (q - q^{-1})^6 e_0 e_1 e_3 e_4 e_2^2$, $D' := \overline{\Omega_B} (D)$. The definition (\ref{eq:bijn}) reads
\begin{equation}
\begin{aligned}
b_{i, j; m} &= \frac{1}{m(q^{-1} - q)} \times \begin{cases}
s & i = j \\
t & i \neq j; \, 2 \in \{i, j\} \\
0 & i \neq j; \, i, j \neq 2
\end{cases}, \\
c_{i, j; m} &= \frac{m(q^{-1} - q)}{s^2 (s^2 - 3t^2)} \times \begin{cases}
s^3 & i = j = 2 \\
s(s^2 - 2t^2) & i = j \neq 2 \\
-s^2 t & i \neq j; \, 2 \in \{i, j\} \\
s t^2 & i \neq j; \, i, j \neq 2
\end{cases},
\end{aligned}
\end{equation}
where $s := [2m]_q$, $t := (-1)^{m - 1} [m]_q$. Finally, the image of $S'_m := (T_{\Check{w}} \otimes T_{\Check{w}}) (S_m)$ is computed as follows.
\begin{align*}
(\pi_B^+ \otimes \pi_B^-) (S'_m) &= \frac{m(q^{-1} - q)}{s^2 (s^2 - 3t^2)} \frac{1}{m(q - q^{-1}) m(q^{-1} - q)} \\
&\quad \times \Bigl[\left\{3\cdot s\left(s^2 - 2 t^2\right) + 6 \cdot st^2\right\} \cdot (-1)^{m - 1} \cdot (-1)^{m - 1} \\
&\quad \quad + 6 \cdot (-s^2 t) \cdot (-1)^{m - 1} \cdot (q^m + q^{-m}) \\
&\quad \quad + (q^m + q^{-m})^2 s^3\Bigl] \\
&\quad \times (qD)^m \otimes (q^{-1}D')^m \\
&= \frac{1}{m(q - q^{-1}) s^2(s^2 - 3t^2)} \cdot (q^{2m} - 1 + q^{-2m}) s^3 \cdot D^m \otimes D'^m \\
&= -\frac{1}{m} \frac{q^m (q^m + q^{-m})}{1 - q^{2m}} (D \otimes D')^m.
\end{align*}
Therefore, surprisingly, the image of $S'_m$ is identical to that in the examples of type $A_{\ell}^{(1)}$, and we have $\pi_B^+ \widehat{\otimes} \pi_B^- (\Theta_{\mathrm{im}}) = \mathbb{E} (-q D\otimes D')^{-1} \mathbb{E} (-q^{-1} D\otimes D')^{-1}$. As a result,

\begin{thm}
The following identity holds in skew formal power series algebra $\widehat{\mathcal{S}_4} := \widehat{\mathcal{S}_B} \cong \Q(q)[[y_0, y_1, y_2, y_3, y_4]]$.
\begin{equation} \label{eq:D4_result}
\begin{aligned}
& \mathbb{E}(\normord{y_2}) \mathbb{E}(\normord{y_4}) \mathbb{E}(\normord{y_3}) \mathbb{E}(\normord{y_1}) \mathbb{E}(\normord{y_0}) \\
&= \left\{\prod_{m \geq 0}^{\rightarrow} X_m \right\} \mathbb{E}(\normord{y_2 y_3 y_4}) \mathbb{E}(\normord{y_1 y_2 y_4}) \mathbb{E}(\normord{y_1 y_2 y_3}) \\
& \times \mathbb{E} (-q \normord{y_0 y_1 y_2^2 y_3 y_4})^{-1} \mathbb{E} (-q^{-1} \normord{y_0 y_1 y_2^2 y_3 y_4})^{-1} \\
& \times \mathbb{E}(\normord{y_0 y_2 y_4}) \mathbb{E}(\normord{y_0 y_2 y_3}) \mathbb{E}(\normord{y_0 y_1 y_2}) \left\{\prod_{m \geq 0}^{\leftarrow} Y_m \right\},
\end{aligned}
\end{equation}
where
\begin{align*}
X_m &= \mathbb{E}(\normord{y_0^{m + 1} y_1^{m} y_2^{2m} y_3^{m} y_4^{m}}) \mathbb{E}(\normord{y_0^{m} y_1^{m + 1} y_2^{2m} y_3^{m} y_4^{m}}) \\
&\times \mathbb{E}(\normord{y_0^{m} y_1^{m} y_2^{2m} y_3^{m + 1} y_4^{m}}) \mathbb{E}(\normord{y_0^{m} y_1^{m} y_2^{2m} y_3^{m} y_4^{m + 1}}) \\
&\times  \mathbb{E}(\normord{y_0^{2m + 1} y_1^{2m + 1} y_2^{4m + 1} y_3^{2m + 1} y_4^{2m + 1}}) \mathbb{E}(\normord{y_0^{m} y_1^{m + 1} y_2^{2m + 1} y_3^{m + 1} y_4^{m + 1}}) \\
&\times \mathbb{E}(\normord{y_0^{m + 1} y_1^{m} y_2^{2m + 1} y_3^{m + 1} y_4^{m + 1}}) \mathbb{E}(\normord{y_0^{m + 1} y_1^{m + 1} y_2^{2m + 1} y_3^{m} y_4^{m + 1}}) \\
&\times \mathbb{E}(\normord{y_0^{m + 1} y_1^{m + 1} y_2^{2m + 1} y_3^{m + 1} y_4^{m}}) \mathbb{E}(\normord{y_0^{2m + 2} y_1^{2m + 2} y_2^{4m + 3} y_3^{2m + 2} y_4^{2m + 2}}), \\
Y_m &= \mathbb{E}(\normord{y_0^{m} y_1^{m + 1} y_2^{2m + 2} y_3^{m + 1} y_4^{m + 1}}) \mathbb{E}(\normord{y_0^{m + 1} y_1^{m + 1} y_2^{2m + 2} y_3^{m + 1} y_4^{m}}) \\
&\times \mathbb{E}(\normord{y_0^{m + 1} y_1^{m + 1} y_2^{2m + 2} y_3^{m} y_4^{m + 1}}) \mathbb{E}(\normord{y_0^{m} y_1^{m + 1} y_2^{2m + 2} y_3^{m + 1} y_4^{m + 1}}) \\
&\times  \mathbb{E}(\normord{y_0^{2m + 1} y_1^{2m + 1} y_2^{4m + 3} y_3^{2m + 1} y_4^{2m + 1}}) \mathbb{E}(\normord{y_0^{m + 1} y_1^{m} y_2^{2m + 1} y_3^{m} y_4^{m}}) \\
&\times \mathbb{E}(\normord{y_0^{m} y_1^{m} y_2^{2m + 1} y_3^{m} y_4^{m + 1}}) \mathbb{E}(\normord{y_0^{m} y_1^{m} y_2^{2m + 1} y_3^{m + 1} y_4^{m}}) \\
&\times \mathbb{E}(\normord{y_0^{m} y_1^{m + 1} y_2^{2m + 1} y_3^{m} y_4^{m}}) \mathbb{E}(\normord{y_0^{2m} y_1^{2m} y_2^{4m + 1} y_3^{2m} y_4^{2m}}).
\end{align*}
\end{thm}
Set $R := \bigl(\begin{smallmatrix} 0 & 0 & 1 & 0 & 0 \\ 1 & 1 & 2 & 1 & 1\end{smallmatrix}\bigr)$. Then the condition of Proposition \ref{prop:S_B_hom} holds. Thus we have continuous algebra homomorphism $\widehat{\psi_4}: \widehat{\mathcal{S}_4} \rightarrow \widehat{\mathcal{S}}$ which satisfies 
\begin{equation}
\widehat{\psi_4}(y_0) = \widehat{\psi_4}(y_1) = \widehat{\psi_4}(y_3) = \widehat{\psi_4}(y_4) = x_2, \quad \widehat{\psi_4}(y_2) = \normord{x_1 x_2^2}.
\end{equation}
Applying $\psi_4$ on (\ref{eq:D4_result}) and transforming variables using $\psi_S$ yields (\ref{eq:D4}). This proves the identity (\ref{eq:D4}).

\begin{mainthm}
The four identities (\ref{eq:A1}), (\ref{eq:A2}), (\ref{eq:A3}), (\ref{eq:D4}) holds in the skew formal power series algebra $\Q(q)[[\frac{x_1}{x_2^2}, x_2]]$.
\end{mainthm}

\subsection*{Acknowledgements}
The author would like to thank Koji Hasegawa, Gen Kuroki for stimulating discussions and helpful advice. The author also appreciate Akihiro Tsuchiya, Yuji Terashima for valuable comments. This work was supported by the Research Institute for Mathematical Sciences, an International Joint Usage/Research Center located in Kyoto University.

\end{document}